\documentclass[12pt,abstracton]{scrreprt}
\usepackage[arrow, matrix, curve]{xy}
\usepackage{graphics}
\usepackage{setspace} 
\usepackage[usenames,dvipsnames,x11names]{xcolor}
\usepackage{tikz}
\usetikzlibrary{arrows,decorations.pathmorphing}
\usetikzlibrary{intersections,shapes.arrows}
\usepackage{amsmath}
\usepackage{upgreek}
\usepackage{epstopdf}
\usepackage{textcomp}
\usepackage{caption}
\usepackage{hyperref}
\usepackage{picinpar}
\usepackage{subcaption}
\usepackage{stmaryrd}
\usepackage{amssymb}
\usepackage{amsthm}
\usepackage{wrapfig}
\usepackage{float}
\usepackage{etex}
\usepackage{mathrsfs} 
\usepackage{booktabs}
\usepackage{longtable}
\usepackage{color}
\usepackage{titling}
\usepackage[figuresright]{rotating}
\usepackage{enumitem, color, amssymb}
\tikzset{
LL/.style={
  draw=black,decorate,
  decoration={snake, segment length=3mm, amplitude=1mm,post length=2mm}
  }
}  
\theoremstyle{definition}
\newtheorem{lem}{Lemma}[section]
\newtheorem{Prop}[lem]{Proposition}
\newtheorem{Theo}[lem]{Theorem}
\newtheorem{Rem}[lem]{Remark}
\newtheorem{Cor}[lem]{Corollary}
\newtheorem{Exp}[lem]{Example}
\newtheorem{Def}[lem]{Definition}
\newtheorem{Con}[lem]{Conjecture}
\newtheorem{Quest}[lem]{Question}

\newcommand\DistTo{\xrightarrow{
   \,\smash{\raisebox{-0.5ex}{\ensuremath{\scriptstyle\sim}}}\,}}
\newcommand{\bew}{\hfill $\boxempty$     
}
\begin{document}
\onehalfspacing
\huge
\begin{center}
Higher degree S-lemma and the stability of quadratic modules
\end{center}
\normalsize
\vspace{0.2cm}
\begin{center}
Philipp Jukic
\end{center}
\vspace{1cm}
\begin{center}
\textbf{Abstract}
\end{center}
In this work we will investigate a certain generalization of the so called S-lemma in higher degrees. 
The importance of this generalization is, that it is closely related to Hilbert's 1888 theorem about tenary quartics. In fact, if such a generalization exits, then one can state a Hilbert-like theorem, where positivity is only demanded on some semi-algebraic set. 
We will show that such a generalization is not possible, at least not without additional conditions. To prove this, we will use and generalize certain tools developed in \cite{z12}. In fact, these new tools will allow us to conclude that this generalization of the S-lemma is not possible because of geometric reasons. Furthermore, we are able to establish a link between geometric reasons and algebraic reasons. This will be accomplished within the framework of quadratic modules.    

\thispagestyle{empty}
\tableofcontents
\thispagestyle{empty} 
\chapter{Introduction}
\pagestyle{plain}
\setcounter{page}{1}
First of all, let us talk about the motivation of this article. In 1888 Hilbert showed in his work \cite{z14} that a ternary quartic $f$, that is a 4-form in three variables, can be written as a sum of three squares of quadratic forms if and only if $f$ is non-negative on $\mathbb{R}^3$. The question is: Can we find a Hilbert-like theorem in a more general setting? What does a more general setting mean in this context? Instead of considering non-negative ternary quartics, we consider a ternary quartic that needs to be non-negative on a semi-algebraic set $S\subseteq \mathbb{R}^3$. Furthermore, the semi-algebraic set $S$ should also satisfy the following two conditions: First, there exists a quadratic form $g$ in three variables such that $S=S(g):=\{x\in \mathbb{R}^3:g(x)\geq 0\}$. Second, the set $S$ has a non-empty interior. 
Of course, if $S\neq \mathbb{R}^3$ then $f$ can, in general, not be written as a sum of three squares of quadratic forms. In this case we need a sort of a correcting term. This correcting term should also satisfy some conditions. First, we demand that this term is of the form $-tg$, where $t$ is a non-negative quadratic form. Second, $f-tg$ should be a non-negative ternary quartic.
Thus a generalization of Hilbert's theorem could look like the following: 
Let $g$ be a quadratic form such that there exists a point $x'\in \mathbb{R}^3$ with $g(x')>0$. 
A ternary quartic $f$ is non-negative on the set $S(g)$ if and only if there exists a non-negative quadratic form $t$ such that $f-tg$ can be written as a sum of three squares of quadratic forms.\\\\
The interpretation of this statement is simple. If $g$ is non-negative, then this statement is equal to Hilbert's statement. If $g$ is not non-negative, then $-tg$ measures 'how far away' $f$ is from being a sum of three squares of quadratic forms.\\
Let us illustrate this statement by considering the two polynomials $g=\mathrm{x}_1^2-\mathrm{x}_2^2$ and $f=\mathrm{x}_1^4-\mathrm{x}_2^4$. It is easy to see that $S(g)$ has a non-empty interior and that $f$ is non-negative on $S(g)$. Furthermore, we have $f-2\mathrm{x}_2^2g=\mathrm{x}_1^4-2\mathrm{x}_1^2\mathrm{x}_2^2+\mathrm{x}_2^4=\left(\mathrm{x}_1^2-\mathrm{x}_2^2\right)^2$. Thus $f-2\mathrm{x}_2^2g$ is a sum of three squares of quadratic forms: One quadratic form is given by $\mathrm{x}_1^2-\mathrm{x}_2^2$, the other two are $0$.\\\\      
We are looking to clarify the following question: Can such a generalization of Hilbert's theorem be made? It turns out that this question is closely related to the so called S-lemma resp. to a certain generalization of the S-lemma.
Hence the \textbf{first Chapter} is all about the introduction and the proof of the S-lemma. The machinery presented in this chapter relies heavily on the work of \cite{z10} and \cite{z1}. 
The results in the first chapter are all well known. Therefore there is nothing new in this part of the article.
In this chapter and throughout the whole article no fancy knowledge will be required. One should be familiar with basic linear algebra, convex geometry, and real algebraic geometry. 
In the \textbf{second Chapter} we will formulate a generalization of the S-lemma. For the sake of simplicity we will refer to the generalization as the S4-conjecture. 
The importance of this S4-conjecture is the following: If the conjecture is true, then the generalization of Hilbert's theorem is possible. If it is not true, then such a generalization is impossible. However, it turns out that it is impossible because we can find a counterexample for the S4-conjecture. 
Although we can find a counterexample, we will still refer to this mentioned generalization as the S4-conjecture.
Next, we do some geometric investigations and finally generalize the counterexample to higher degrees. In the \textbf{third and last chapter} we use and generalize the machinery developed in \cite{z12} to further investigate the counterexample.
It turns out that the tools presented in \cite{z12} are quite suitable in analyzing the S4-conjecture. 
In fact, by using these new methods we will see that the conjecture fails because of geometric reasons. 
Since \cite{z12} connects geometric properties and algebraic properties, we will see that there is an interesting link between the S4-conjecture and the stability of quadratic modules.
Finally, this article will be concluded by presenting some new questions that should serve as a motivation for further studies.\newpage
We will use the following \textbf{notation} throughout this article:
\small
\begin{itemize}
\item $\mathbb{R}$, $\mathbb{C}$, $\mathbb{Z}$, $\mathbb{N}, \mathbb{N}_0$: The real, complex, integer, natural numbers and the natural numbers with $0$.
\item $\mathbb{R}^n$, $\mathbb{C}^n$, $\mathbb{Z}^n$: The $0\leq n$ dimensional vector spaces $\mathbb{R}^n$ resp. $\mathbb{C}^n$ and the free $\mathbb{Z}$-module $\mathbb{Z}^n$. 
\item $\langle\cdot,\cdot\rangle$: The standard scalar product in $\mathbb{R}^n$. 
\item $\mathbb{K}\left[\mathrm{x}_1,\ldots,\mathrm{x}_n\right]$: The polynomial ring over a field $\mathbb{K}$ in $n\geq 1$ variables. Polynomial variables will always be denoted by upright letters $\mathrm{x},\mathrm{y},\upbeta,\uplambda$ etc.
\item $\mathbb{K}\left[\mathrm{x}_1,\ldots,\mathrm{x}_n\right]_d$: The set of all polynomials $f\in \mathbb{K}\left[\mathrm{x}_1,\ldots,\mathrm{x}_n\right]$ with $\deg(f)\leq d$.
\item A polynomial $f\in \mathbb{R}[\mathrm{x}_1,\ldots,\mathrm{x}_n]$ is called non-negative if $\forall x\in \mathbb{R}^n:f(x)\geq 0$. Negative, positive and non-positive polynomials are defined in the same manner.
\item The homogenization of a polynomial $f\in \mathbb{K}[\mathrm{x}_1,\ldots,\mathrm{x}_n]$ will be denoted by $\overline{f}$. The dehomogenization of a homogeneous polynomial $g$ with $\tilde{g}$.
\item $\mathbb{A}^n$, $\mathbb{P}^n$: The $n$-dimensional affine space and the $n$-dimensional projective space.  
\item $\mathcal{V}(f_1,\ldots,f_s)$: For polynomials $f_1,\ldots,f_s\in \mathbb{K}[\mathrm{x}_1,\ldots,\mathrm{x}_n]$ the set $\mathcal{V}(f_1,\ldots,f_s)$ is defined to be the set of all solutions $x\in \overline{\mathbb{K}}^n$ ($\overline{\mathbb{K}}$ denotes the algebraic closure of $\mathbb{K}$) of the polynomial equalities $f_1(x)=0,\ldots,f_s(x)=0$. If $\mathbb{K}=\mathbb{R}$ then we will fix $\mathbb{C}$ as the algebraic closure of $\mathbb{R}$. If $f_1,\ldots,f_s$ are homogeneous, then we can interpret $\mathcal{V}(f_1,\ldots,f_s)$ as the set of all solutions in the projective space $\mathbb{P}^n$. 
\item Let $V$ be a variety defined over a field $\mathbb{K}$ and $\mathbb{L}|\mathbb{K}$ an algebraic extension of $\mathbb{K}$. The $\mathbb{L}$-rational points of $V$ are denoted by $V(\mathbb{L})$. 
\item $S(f_1,\ldots,f_s)$: The basic closed semi-algebraic set $$S(f_1,\ldots,f_s)=\left\{x\in \mathbb{R}^n:f_1(x)\geq 0,\ldots,f_s(x)\geq 0\right\}$$ 
defined by the polynomials $f_1,\ldots,f_s\in \mathbb{R}[\mathrm{x}_1,\ldots,\mathrm{x}_n]$.  
\item $\mathrm{GL}_n$, $\mathrm{O}_n$: The general linear group over $\mathbb{R}$ and its orthogonal subgroup over $\mathbb{R}$.
\item Let $X$ be a topological space and $A$ a subset. The interior of $A$ is denoted with $\mathrm{int}(A)$ and the closure with $\overline{A}$. 
\end{itemize} \normalsize\mdseries                              
\newpage
\chapter{The S-lemma}\label{k1}
In this chapter we will formulate and prove the so called S-lemma. Before doing this, however, it shall be noted that the S-lemma has many variations in the literature. While all versions are in fact equivalent, we will use a version that is closer to real algebraic geometry. Thus the original statement of the S-lemma made by Yakubovich, that can be found in the work of Polik and Terlaky \cite{z10}, will not be used. In the sense of real algebraic geometry the S-lemma is formulated in the following way:
\begin{Theo}\label{t1}\textbf{S-lemma}:
Let $f,g$ be polynomials in $\mathbb{R}[\mathrm{x}_1,\ldots,\mathrm{x}_n]_2$. If there exists a point $x'\in \mathbb{R}^n$ with $g(x')>0$, then the following statements are equivalent:
\begin{enumerate}[label=(\alph*)]
\item The inclusion $S(g)\subseteq S(f)$ holds.
\item There exists a non-negative real number $t$ such that $f(x)-tg(x)\geq 0$ for all $x\in \mathbb{R}^n$
\end{enumerate}
\end{Theo}   
\noindent The aim of this chapter is to provide a proof for Theorem \ref{t1}. Simultaneously, it should serve as an introduction in what is to come later. Before we are ready to prove Theorem \ref{t1}, we need some preparatory results, which will be bundled together in the following section.
\section{Preliminaries}
First of all, it is worth mentioning that one could prove Theorem \ref{t1} directly, without any notable machinery. One such proof can be found in \cite[pp. 376-378]{z10}. The disadvantage is, however, that it  needs quite a lot of computations. As already pointed out, we will use a different approach. For the proof of theorem \ref{t1} we will need the following definitions, lemmas, and propositions:
\begin{Def}
Let $f=\sum_{i=1}^n\sum_{j=1}^na_{ij}\mathrm{x}_i\mathrm{x}_j$ be a quadratic form in $\mathbb{R}[\mathrm{x}_1,\ldots,\mathrm{x}_n]$, where all coefficients $a_{ij}$ of $f$ lie in $\mathbb{R}$.
The matrix that corresponds to $f$ is defined to be
the symmetric matrix $A_f=\left(\frac {1} {2}(a_{ij}+a_{ji})\right)_{1\leq i\leq n,1\leq j\leq n}$.
If $f$ is an arbitrary form of degree $d\geq 0$ in $\mathbb{R}[\mathrm{x}_1,\ldots,\mathrm{x}_n]$, then $f$ is said to be positive semi-definite resp. positive definite if $\forall x\in \mathbb{R}^n:f(x)\geq 0$ resp. $\forall x\in \mathbb{R}^n\backslash \{0\}:f(x)>0$.
\end{Def} 
\begin{Rem}
A quadratic form $f$ is positive (semi-) definite if and only if the corresponding matrix $A_f$ is positive (semi-) definite.
\end{Rem}
\noindent In the following we will assume that the coefficients of a quadratic form in $\mathbb{R}[\mathrm{x}_1,\ldots,\mathrm{x}_n]$ lie in $\mathbb{R}$.
\begin{Def}\label{d3}
Let $P_{d,n}$ be the set of all forms of even degree $d>0$ in $\mathbb{R}[\mathrm{x}_1,\ldots,\mathrm{x}_n]$. With $P_{d,n}^+$ we denote the subset of $P_{d,n}$ that consist of all positive semi-definite forms in $P_{d,n}$. 
\end{Def}
\begin{Def}
Let $V$ be finite dimensional $\mathbb{R}$-vector space. A (convex) cone $C\subseteq V$ is a subset of $V$ that satisfies the following two conditions:
\begin{itemize}
\item The set $C$ is not empty.
\item For any real number $\lambda\geq 0$ and any element $g\in C$ we have $\lambda g\in C$.
\end{itemize}
We say that a cone $C\subseteq V$ is pointed, if the identity $C\cap -C=\{0\}$ holds. 
\end{Def}
\begin{Rem}\label{rr1}
One can easily see that $P_{2,n}$ is a finite dimensional $\mathbb{R}$-vector space. To be more precise, there is a vector-space isomorphism $P_{2,n}\DistTo\mathbb{R}^{\frac {n(n+1)} {2}}$. For $f,g\in P_{2,n}$ the dot product on $P_{2,n}$ is defined by $\langle f,g\rangle=\mathrm{tr}(A_fA_g)$, which is just the pullback of the dot product in $\mathbb{R}^{\frac {n(n+1)} {2}}$. Thus $P_{2,n}$ is an euclidean space. The same is also true for $P_{d,n}$, where $d\geq 0$.\par\smallskip   
Finally, it should be noted that this vector space $P_{2,n}$ has a more or less surprising upcoming in algebraic geometry: See \cite[Example 3, p. 44]{z2} about determinantal varieties.      
\end{Rem}
\begin{Def}
Let $V$ be a finite dimensional real vector space and $C\subseteq V$ a convex subset. A convex subset $F$ of $C$ is called a face of $C$ if the following statement holds: 
Suppose $u$ and $v$ are two points in $C$. If there exists a $\lambda\in (0,1)$ such that $\lambda u+(1-\lambda)v\in F$, then $u$ and $v$ lie already in $F$.\par\smallskip  
A face $F$ of $C$ is called proper if $\varnothing\subsetneqq F\subsetneqq C$ holds. If there is a point $u\in C$ such that $\{u\}$ is a face of $C$, then the point $u$ is called an extremal point. With $\mathrm{ex}(C)$ we will denote the set of all extremal points of $C$.
\end{Def}
\begin{Def}
Let $V$ be a finite dimensional real vector space and $C\subseteq V$ a convex subset. Let $H$ be a hyper plane given by $H=\left\{x\in V:\ell(x)=0\right\}$, where $\ell:V\rightarrow \mathbb{R}$ is a linear form. Set $\overline{H}_+=\left\{x\in V:\ell(x)\geq 0\right\}$ and $\overline{H}_-=\left\{x\in V:\ell(x)\leq 0\right\}$. A face $F$ of $C$ is called exposed if there exists a linear form $\ell:V\rightarrow \mathbb{R}$ such that $C$ is contained in $\overline{H}_+$ or $\overline{H}_-$ and $F=C\cap H$.  
\end{Def}
\begin{Def}
Let $V$ be a finite dimensional real vector space and $S$ an arbitrary subset of $V$. The affine hull $\mathrm{aff}(S)$ of $S$ is defined by $$\mathrm{aff}(S)=\left\{\sum_{i=1}^n\lambda_is_i:s_1,\ldots,s_n\in S,\lambda_1,\ldots,\lambda_n\in \mathbb{R},\sum_{i=1}^n\lambda_i=1,n\in \mathbb{N}\right\}.$$ 
\end{Def}
\begin{Def}
Let $V$ be a finite dimensional real vector space. The dimension of a convex set $C\subseteq V$ is defined to be the dimension of its affine hull. In short $\dim(C)=\dim(\mathrm{aff}(C))$. 
\end{Def}
\begin{Def}
Let $V$ be a finite dimensional real vector space and $C$ a convex subset of $V$. Let $B_\varepsilon(x)$ denote the open ball in $x$ with radius $\varepsilon>0$. The relative interior $\mathrm{relint}(C)$ of $C$ in $V$ is defined by $$\mathrm{relint}(C)=\left\{x\in C:\exists \varepsilon>0:B_\varepsilon(x)\cap \mathrm{aff}(C)\subseteq C\right\}.$$ 
\end{Def}
\begin{lem}\label{l10}
(a): For every $x\in \mathbb{R}^n\backslash \{0\}$ the symmetric $n\times n$-matrix $xx^T$ matrix is positive semi-definite of rank $1$.\par\smallskip  
(b): Let $A$ be a positive semi-definite matrix. The rank of $A$ is the smallest natural number $r$ such that $A$ can be written as $A=\sum_{i=1}^rx_ix_i^T$ for some $x_1,\ldots,x_r\in \mathbb{R}^n$.  
\end{lem}
\textbf{Proof}: (a): Trivial.\par\smallskip  
(b): First of all, note that statement (b) is independent with respect to transformations $S^TAS$, where $S\in \mathrm{O}_n$. Indeed, set $D=S^TAS$ and assume that $A=\sum_{i=1}^rx_ix_i^T$, where $r$ is minimal. Then we have $D=S^T\sum_{i=1}^rx_ix_i^TS=\sum_{i=1}^rS^Tx_ix_i^TS=\sum_{i=1}^rS^Tx_i(S^Tx_i)^T$. 
It is clear that $r$ is also the minimal length of the sum for $D$: Otherwise, $A$ could be written as a sum of smaller length, which would be a contradiction. Choose $S\in \mathrm{O}_n$ such that $D$ is a diagonal matrix. The diagonal of $D$ consists of the eigenvalues of $A$, which are all non-negative. Thus it is easy to see that $D$ can be written as a sum $\sum_{i=1}^rx_ix_i^T$ for some $x_1,\ldots,x_n\in \mathbb{R}^n$ and $r=\mathrm{rk}(A)$. It remains to verify that $r$ is minimal. But this follows from $\mathrm{rk}\left(\sum_{i=1}^rx_ix_i^T\right)\leq r\mathrm{rk}(x_ix_i^T)=r$.\bew 
\begin{Prop}\label{p11}
Let $d\geq 0$ be an even number.
\begin{enumerate}[label=(\alph*)]
\item The set $P_{d,n}^+$ is a closed cone in $P_{d,n}$.
\item The cone $P_{d,n}^+$ is pointed.
\item Let $L\subseteq\mathbb{R}^n$ be a subspace and $F_L:=\left\{f\in P_{2,n}^+:\forall x\in L:f(x)=0\right\}$. The set $F_L$ is an exposed face with $\dim(F_L)=\frac {r(r+1)} {2}$ and $r=n-\dim(L)$. If $f\in P_{2,n}^+$ and $L=\mathrm{ker}\left(A_f\right)$, then $f$ is in $\mathrm{relint}(F_L)$.
\end{enumerate}
\end{Prop}
\textbf{Proof}: (a): We will show that $P:=P_{d,n}\backslash P_{d,n}^+$ is open. Take an element $f\in P$.
Since $f\in P$, there exists a point $x\in \mathbb{R}^n$ such that $f(x)<0$. Consider the evaluation homomorphism $\mathrm{ev}_x:P_{d,n}\rightarrow \mathbb{R},p\mapsto p(x)$. Furthermore, $P_{d,n}$ is, as already stated in Remark \ref{rr1}, an euclidean space. Thus $\mathrm{ev}_x$ is continuous. Let $U\subseteq \mathbb{R}$ be an open neighborhood of $f(x)$ such that all elements of $U$ are negative real numbers.
The set $U'=\mathrm{ev}_x^{-1}\left(U\right)$ is an open neighborhood of $f$ that satisfies $U'\subseteq P$. Thus we proved that $P$ is open resp. that $P_{d,n}^+$ is closed. The second assertion that $P_{d,n}^+$ is a cone is trivial.\par\smallskip  
(b): Trivial.\par\smallskip  
(c): 
In the following we will just omit the trivial parts of the proof.\footnote{Pay attention to statements that begin with 'It is easy to see'.}
Let us begin with the easiest part, verifying that $\dim(F_L)=\frac {r(r+1)} {2}$, where $r=n-\dim(L)$. This can be done by proving $\mathrm{aff}(F_L)\cong \mathbb{R}^{\frac {r(r+1)} {2}}$. 
Let $L^\perp$ be the orthogonal complement of $L$ in $\mathbb{R}^n$. Without loss of generality we can identify $L$ with $\mathbb{R}^{n-r}$ and $L^\perp$ with $\mathbb{R}^r$.
Consider the cone $P_{2,r}^+$. It is easy to see that $F_L$ can be identified with $P_{2,r}^+$. Furthermore, $P_{2,r}^+$ has a non-empty interior. A well known result in convex geometry states that a cone with non-empty interior is full. This means that $\mathrm{aff}\left(P_{2,r}^+\right)=P_{2,r}$. Thus $\mathrm{aff}\left(P_{2,r}^+\right)\cong \mathbb{R}^{\frac {r(r+1)} {2}}$. Identifying $\mathrm{aff}(F_L)$ with $\mathrm{aff}\left(P_{2,r}^+\right)$ proves the assertion.
Next, we show that $F_L$ is an exposed face. A quadratic form $h_1\in P_{2,n-r}^+$ and a quadratic form $h_2\in P_{2,r}^+$ give rise to a quadratic form $h\in P_{2,n}^+$ in an obvious manner. In fact, the corresponding matrix $A_h$ of $h$ is given by $A_h=\left(\begin{array}{cc}
A_{h_1} & 0\\
0& A_{h_2}
\end{array}\right)$. Fix $h_1=\mathrm{x}_1^2+\cdots+\mathrm{x}_{n-r}^2$ and define $\tilde{A}_{h_1}=\left(\begin{array}{cc}
A_{h_1} & 0\\
0& 0
\end{array}\right)$, $\tilde{A}_{h_2}=\left(\begin{array}{cc}
0 & 0\\
0& A_{h_2}
\end{array}\right)$. Then we can identify $F_L$ with $\left\{h\in P_{2,n}^+:\exists h_2\in P_{2,r}^+:A_h=\tilde{A}_{h_2}\right\}$. It is easy to see that $F_L$ consists of all quadratic forms $h\in P_{2,n}^+$ that satisfy $\mathrm{tr}\left(A_h\tilde{A}_{h_1}\right)=0$. Thus it is convenient to consider the linear form $\ell:P_{2,n}\rightarrow \mathbb{R},p\mapsto \mathrm{tr}\left(A_p\tilde{A}_{h_1}\right)$ and the hyper plane $H=\left\{p\in P_{2,n}:\ell(p)=0\right\}$. So far, we know that $F_L=P_{2,n}^+\cap H$.
Finally, we just have to deal with the inclusion $P_{2,n}^+\subseteq \overline{H}_+$. Let $h=\sum_{i,j}a_{ij}\mathrm{x}_i\mathrm{x}_j$ be a quadratic form in $P_{2,n}^+$. Since $h$ is non-negative, the coefficients $a_{ii}$ must be non-negative for all $i=1,\ldots,n$. Thus the diagonal of $A_h$ consists of non-negative real numbers. This implies $\ell(h)\geq 0$. Altogether we proved that $F_L$ is an exposed face.
Let us deal with the last statement in (c). Suppose $f\in P_{2,n}^+$ and $L=\mathrm{ker}(A_f)$. As before set $r=n-\dim(L)$. Again we identify $F_L$ with $P_{2,r}^+$ and interpret $f$ as a quadratic form in $P_{2,r}^+$. 
Then $f$ does only vanish at the origin in $\mathbb{R}^r$. 
Hence $f$ lies in the interior of the cone $P_{2,r}^+$ resp. in $\mathrm{relint}\left(F_L\right)$, which proves the assertion.\bew

\begin{lem}\label{l12}
Let $C\subset \mathbb{R}^n$ be a non-empty closed convex set which does not contain any straight line. Then $\mathrm{ex}(C)$ is a non-empty set.
\end{lem}
\textbf{Proof}: See \cite[Lemma 3.5, p. 53]{z1}.\bew\\\\
Because the next result is very important, it will be proven, although there exists a suitable reference.
\begin{Prop}\label{p13}
Let $L$ be an affine subspace of $P_{2,n}$ such that $S=L\cap P_{2,n}^+$ is not empty. Suppose the inequality $\mathrm{codim}_{P_{2,n}}(L)<\frac {(r+2)(r+1)} {2}$ holds for some $r\in \mathbb{N}_0$. Then there exists a quadratic form $f\in S$ such that the rank of $A_f$ is bounded by $r$.
\end{Prop}
\textbf{Proof}\cite[Proposition 13.1, p. 83]{z1}: According to Proposition \ref{p11} the cone $P_{2,n}^+$ is pointed and closed. This means that there is no way that the cone $P_{2,n}^+$ contains a straight line. If $P_{2,n}^+$ does not contain such a line so does not the subset $S$ of $P_{2,n}^+$. By using Lemma \ref{l12} we get $\mathrm{ex}(S)\neq \varnothing$. Choose an arbitrary $f\in S$ and let $A_f$ be its corresponding matrix of rank $m$. Consider $W=\mathrm{ker}(A_f)$ and the exposed face $F_W$ (Proposition \ref{p11}). We want to show that $f$ is an element of the set $\mathrm{relint}(L\cap F_W)$. Since $f$ is an element of $\mathrm{relint}(F_W)$ and $\mathrm{relint}(L)$, it is enough to verify the inclusion $\mathrm{relint}(F_W)\cap \mathrm{relint}(L)\subset \mathrm{relint}(L\cap F_W)$.
Take a point $g\in \mathrm{relint}(F_W)\cap \mathrm{relint}(L)$. There exist $\varepsilon_1,\varepsilon_2>0$ such that $B_{\varepsilon_1}(g)\cap \mathrm{aff}(F_W)\subset F_W$ and $B_{\varepsilon_2}(g)\cap \mathrm{aff}(L)\subset L$. By setting $\varepsilon=\min\{\varepsilon_1,\varepsilon_2\}$ we get $B_\varepsilon(g)\cap \mathrm{aff}(L)\cap \mathrm{aff}(F_W)\subset L\cap F_W$. Since $\mathrm{aff}(L\cap F_W)\subset \mathrm{aff}(L)\cap \mathrm{aff}(F_W)$, we have $B_\varepsilon(g)\cap \mathrm{aff}(L\cap F_W)\subset L\cap F_W$ which implies the assertion.\par\smallskip  
We know that $f$ lies in both sets, $\mathrm{relint}(L\cap F_W)$ and $\mathrm{ex}(L\cap F_W)$. This can only work if $\dim(L\cap F_W)=0$ holds: Suppose $\dim(L\cap F_W)>0$. For every $\varepsilon>0$ we can find two different points $\delta_1,\delta_2\in B_\varepsilon(f)$ such that $\delta_1,\delta_2\neq f$ and $\delta_1,\delta_2,f\in \mathrm{aff}(L\cap F_W)$. Choose $\varepsilon>0$ such that $B_\varepsilon(f)\cap \mathrm{aff}(L\cap F_W)\subset L\cap F_W$ holds. Now, $f$ is some point on the line segment that connects the two points $\delta_1$ and $\delta_2$. But both points lie in $L\cap F_W$. This contradicts the fact that $f$ is an extremal point of $L\cap F_W$.\par\smallskip  
Since $\dim(L\cap F_W)=0$, we get $\dim(L)+\dim(F_W)=\dim(L+F_W)\leq \dim(P_{2,n})$. This and Proposition \ref{p11} imply $\frac {m(m+1)} {2}=\dim(F_W)\leq \dim(P_{2,n})-\dim(L)=\mathrm{codim}_{P_{2,n}}(L)<\frac {(r+1)(r+2)} {2}$ and thus $m<r+1$.\bew  
\begin{Cor}\label{c14}
Let $f_1,f_2\in P_{2,n}$ be two quadratic forms. The two equations $f_1(x)=\alpha_1$ and $f_2(x)=\alpha_2$ have a simultaneous solution $x\in \mathbb{R}^n$ if and only if there exists a quadratic form $q\in P_{2,n}^+$ such that $\mathrm{tr}\left(A_{f_1}A_q\right)=\alpha_1$ and $\mathrm{tr}\left(A_{f_2}A_q\right)=\alpha_2$. 
\end{Cor}
\textbf{Proof}\cite[Corollary 13.2, p. 84]{z1}: $\Rightarrow$: Choose $x\in \mathbb{R}^n$ such that $f_1(x)=\alpha_1$ and $f_2(x)=\alpha_2$. Define $X=xx^T$. Then $\mathrm{tr}\left(A_{f_i}X\right)=\mathrm{tr}\left(A_{f_i}xx^T\right)=x^TA_{f_i}x=f_i(x)=\alpha_i$ holds for $i=1,2$.\par\smallskip  
$\Leftarrow$: The map $\ell_i:P_{2,n}\rightarrow \mathbb{R},p\mapsto \mathrm{tr}\left(A_{f_i}A_p\right)$ is obviously a vector space homomorphism for $i=1,2$. 
It is easy to see that $\ell_i^{-1}(\alpha_i)=\mathrm{ker}(\ell_i)+q$ holds for $i=1,2$. Hence $\dim(\ell_i^{-1}(\alpha_i))=\dim(P_{2,n})-1$ for $i=1,2$.
This implies $\mathrm{codim}_{P_{2,n}}(L)<3=\frac {(1+1)(2+1)} {2}$, where $L=\ell_1^{-1}(\alpha_1)\cap \ell_2^{-1}(\alpha_2)$. According to Proposition \ref{p13} we can find a quadratic form $h$ in $L\cap P_{2,n}^+$ such that $\mathrm{tr}\left(A_{f_i}A_h\right)=\alpha_i$ and $\mathrm{rk}(A_h)\leq 1$ for $i=1,2$. Now, Lemma \ref{l10} tells us that there exists a point $x\in \mathbb{R}^n$ with $A_h=xx^T$. Substituting $A_h$ through $xx^T$ in $\mathrm{tr}\left(A_{f_i}A_h\right)$ results in $f_i(x)=\alpha_i$ for $i=1,2$.\bew   
\begin{Cor}\label{c1.15}
Consider the quadratic forms $f_1,f_2\in P_{2,n}$ and the map $\varphi:\mathbb{R}^n\rightarrow \mathbb{R}^2,x\mapsto (f_1(x),f_2(x))^T$. Then the set $M=\varphi(\mathbb{R}^n)$ is a convex subset of $\mathbb{R}^2$.
\end{Cor} 
\textbf{Proof}\cite[Corollary 13.3, p. 84]{z1}: Consider the map $$\psi:P_{2,n}\rightarrow \mathbb{R}^2,h\mapsto \left(\mathrm{tr}(A_{f_1}A_h),\mathrm{tr}(A_{f_2}A_h)\right).$$ Since $\psi$ is linear and $P_{2,n}^+$ is a cone (Proposition \ref{p11}), the image of $P_{2,n}^+$ under $\psi$ is a cone. Finally, Corollary \ref{c14} implies the equality $\psi(P_{2,n}^+)=M$.\bew
\begin{Prop}\label{p1.16}
Let $A$ and $B$ be two non-empty convex subsets of $\mathbb{R}^n$ such that $A\cap B=\varnothing$. Then there exists a linear form $\ell:\mathbb{R}^n\rightarrow \mathbb{R}$ such that $\ell(x)\leq \ell(y)$ for all $x\in A$ and $y\in B$.
\end{Prop}
\textbf{Proof}: See \cite[Proposition 1.2, p. 106]{z1}.\bew
\section{Proof of the S-lemma}
Before actually proving theorem \ref{t1} we will consider a special case, from which the Theorem \ref{t1} will easily follow. 
\begin{Prop}\label{p2.1}
\textbf{Homogeneous S-lemma}:
Let $f,g$ be quadratic forms in\\
$\mathbb{R}[\mathrm{x}_1,\ldots,\mathrm{x}_n]$. If there exists a point $x'\in \mathbb{R}^n$ with $g(x')>0$, then following statements are equivalent:
\begin{enumerate}[label=(\alph*)]
\item The inclusion $S(g)\subseteq S(f)$ holds.
\item There exists a non-negative real number $t\geq 0$ such that $f(x)-tg(x)\geq 0$ for all $x\in \mathbb{R}^n$
\end{enumerate}
\end{Prop} 
\textbf{Proof}\cite[Proposition 2.3, p. 377]{z10}: (b)$\Rightarrow$(a): This implication is quite trivial: Suppose we could find a non-negative real number $t$ such that $f(x)-tg(x)\geq 0$ holds for all $x\in \mathbb{R}^n$ and a point $y\in \mathbb{R}^n$ with $g(y)\geq 0$ and $f(y)<0$. It is clear that $f(y)-tg(y)$ would be negative, contradicting $f(x)-tg(x)\geq 0$ for all $x\in \mathbb{R}^n$.\\
(a)$\Rightarrow$(b): According to Corollary \ref{c1.15} the set $M=\left\{(f(x),g(x)):x\in \mathbb{R}^n\right\}$ is a convex subset of $\mathbb{R}^2$. Define $C=\left\{(u_1,u_2):u_1<0,u_2\geq 0\right\}$. Because $S(g)\subseteq S(f)$ holds, the intersection between $M$ and the convex set $C$ is empty. According to Proposition \ref{p1.16} there exists a linear form $\ell:\mathbb{R}^2\rightarrow \mathbb{R}$ such that $\ell(x)\leq \ell(y)$ for all $x\in M$ and $y\in C$. Since $0\in M$ and $0\in \overline{C}$, we have $\ell(x)\leq 0\leq \ell(y)$ for all $x\in M$ and $y\in C$. 
Choose $\alpha_1,\alpha_2\in \mathbb{R}$ such that $\ell=\alpha_1\mathrm{x}_1+\alpha_2\mathrm{x}_2$.
Consider the two statements 
\begin{eqnarray*}
\begin{split}
\forall (u_1,u_2)\in C:\alpha_1u_1+\alpha_2u_2\geq 0\\
\forall x\in \mathbb{R}^n:\alpha_1f(x)+\alpha_2g(x)\leq 0.
\end{split}
\end{eqnarray*}
If we take the point $(-1,0)\in C$ and evaluate $\ell$ at $(-1,0)$, we get $\ell(-1,0)=-\alpha_1\geq 0$. Thus $\alpha_1$ must be a non-positive real number. 
For an arbitrary $\varepsilon>0$ the point $(-\varepsilon,1)$ lies in $C$. Evaluating $\ell$ at the point $(-\varepsilon,1)$ leads to $\ell(-\varepsilon,-1)=-\alpha_1\varepsilon+\alpha_2\geq 0$. If $\alpha_2\neq 0$ then $\alpha_2$ must be positive. Otherwise, we could choose $\varepsilon>0$ so small that the inequality $-\alpha_1\varepsilon+\alpha_2<0$ would hold.  
We know that there exists a point $x'\in \mathbb{R}^n$ such that $g(x')>0$. From the inequality $\alpha_1f(x')+\alpha_2g(x')\leq 0$ we conclude that $\alpha_1$ cannot vanish. Set $\alpha=\frac {\alpha_2} {\alpha_1}$. 
Without loss generality we can assume that $\alpha\neq 0$. Otherwise, $f$ would be non-negative and therefore we could set $t=0$.
Since $\alpha<0$, we have the inequality $\alpha \alpha_1f(x)+\alpha\alpha_2g(x)=\alpha_2(f(x)+\alpha g(x))\geq 0$ for all $x\in \mathbb{R}^n$. Hence $f+\alpha g$ is a non-negative polynomial. Set $t=-\alpha$ and the assertion follows.\bew\\\\
\textbf{Proof of theorem \ref{t1}}: (b)$\Rightarrow$(a): This implication is trivial.\par\smallskip  
(a)$\Rightarrow$(b): Without loss of generality we can assume that $x'=0$. 
Set $\mathrm{x}=(\mathrm{x}_1,\ldots,\mathrm{x}_n)$.
Let
\begin{eqnarray*}
\begin{split}
f=f_1(\mathrm{x})+f_2(\mathrm{x})+c_f\\
g=g_1(\mathrm{x})+g_2(\mathrm{x})+c_g
\end{split}
\end{eqnarray*}
be the decompositions of $f$ and $g$ with respect to the standard grading in $\mathbb{R}[\mathrm{x}]$,
where $f_1$ resp. $g_1$ denotes the component of $f$ resp. $g$ that has degree $2$, $f_2$ resp. $g_2$ denotes the component of $f$ resp. $g$ that has degree $1$ and finally, $c_f$ resp. $c_g$ denotes the constant component. Let $\overline{f},\overline{g}\in \mathbb{R}[\mathrm{x},\mathrm{y}]$ be given by:
 \begin{eqnarray*}
\begin{split}
\overline{f}=f_1(\mathrm{x})+\mathrm{y}f_2(\mathrm{x})+\mathrm{y}^2c_f\\
\overline{g}=g_1(\mathrm{x})+\mathrm{y}g_2(\mathrm{x})+\mathrm{y}^2c_g.
\end{split}
\end{eqnarray*}
In fact, we just need to prove that $\overline{f}$ and $\overline{g}$ satisfy the condition (a) in Proposition \ref{p2.1}: Since then, there would exist a non-negative real number $t$ such that $\overline{f}(x,y)-t\overline{g}(x,y)\geq 0$ for all $(x,y)\in \mathbb{R}^{n+1}$ and the assertion would follow from the dehomogenization of $\overline{f}(\mathrm{x},\mathrm{y})-t\overline{g}(\mathrm{x},\mathrm{y})$.\\
Suppose we could find a point $(x,y)\in \mathbb{R}^{n+1}$ with $y\neq 0$ such that    
\begin{eqnarray*}
\begin{split}
\overline{f}(x,y)<0\\
\overline{g}(x,y)\geq 0.
\end{split}
\end{eqnarray*}
Then we would get a contradiction, since the two identities 
\begin{eqnarray*}
\begin{split}
f\left(\frac {x} {y}\right)=\frac {\overline{f}(x,y)} {y^2}\\
g\left(\frac {x} {y}\right)=\frac {\overline{g}(x,y)} {y^2}
\end{split}
\end{eqnarray*}
hold.
Suppose we could find a point $(x,0)\in \mathbb{R}^{n+1}$ with 
\begin{eqnarray*}
\begin{split}
\overline{f}(x,0)<0\\
\overline{g}(x,0)\geq 0.
\end{split}
\end{eqnarray*} 
This two inequalities imply that $f_1(x)<0$ and $g_1(x)\geq 0$. Consider $f(\uplambda x)=\uplambda^2f_1(x)+\uplambda f_2(x)+c_f$ as a polynomial in the new variable $\uplambda$. Since $f_1(x)<0$, we get that $\lambda^2f_1(x)+\lambda f_2(x)+c_f$  converges to $-\infty$ as $|\lambda|$ converges to $\infty$.\par\smallskip  
Acknowledging that $c_g>0$, $g_1(x)\geq 0$ and treating $g(\uplambda x)=\uplambda^2g_1(x)+\uplambda g_2(x)+c_g$ as a polynomial in $\uplambda$, leads us to the following distinctions
\begin{itemize}
\item $g(\lambda x)\rightarrow \infty$ for $\lambda\rightarrow \infty$ if $g_1(x)>0$
\item $g(\lambda x)\rightarrow \infty$ for $\lambda\rightarrow \mathrm{sign}(g_2(x))\infty$ if $g_1(x)=0$ and $g_2(x)\neq 0$.
\end{itemize}
This proves that no matter what, we can always find a suitable $\lambda\in \mathbb{R}$ such that $f(\lambda x)<0$ and $g(\lambda x)\geq 0$ are satisfied, which clearly contradicts our assumption.\bew 
\newpage
\chapter{Higher degree S-lemma}
\section{Counterexample}
Let us revisit Proposition \ref{p2.1}. The question here is, if there is a generalization of the mentioned proposition in higher degrees. To be more precise, can we give up the restriction that the degree of the two homogeneous polynomials in Proposition \ref{p2.1} is bounded by $2$? The answer is no, as the following simple example illustrates it.
\begin{Exp}\label{l3.1}
Consider $g=\mathrm{x}_1^2-\mathrm{x}_2^2$ and $f=\mathrm{x}_1^2(\mathrm{x}_1^2-\mathrm{x}_2^2)$. Then there is no non-negative real number $t$ such that $f(x)-tg(x)\geq 0$ holds for all $x\in \mathbb{R}^2$. It is easy to check that $g$ and $f$ satisfy the prerequisites and the condition (a) of Proposition \ref{p2.1}. Suppose we could find a non-negative real number $t$ such that $f(x)-tg(x)\geq 0$ holds for all $x\in \mathbb{R}^2$. Let $(k_n)_n\subset \mathbb{R}$ be a sequence such that $k_n\rightarrow 0$ for $n\rightarrow \infty$. The inequality $f(k_n,0)-tg(k_n,0)\geq 0$ implies $f(k_n,0)\geq tg(k_n,0)$ for all $n\in \mathbb{N}$. But this is not true. Choose a natural number $N\in \mathbb{N}$ such that $k_N^2<t$. Then $tg(k_N,0)$ would be greater than $f(k_N,0)$, which would contradict
our assumption.
\end{Exp} 
\begin{window}[3, r, \includegraphics[scale=0.2]{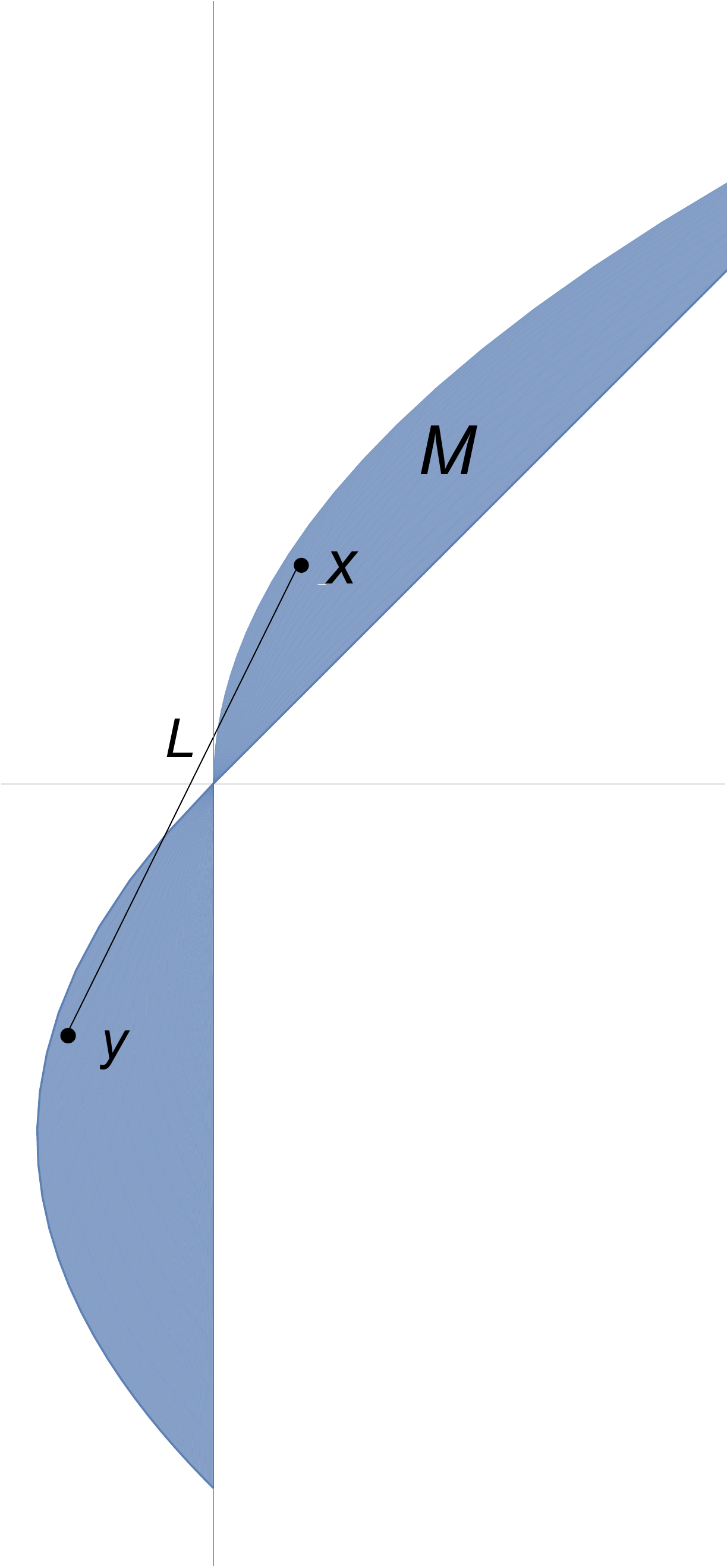},{}]
\begin{Rem}
Since we cannot generalize Proposition \ref{p2.1} to higher degrees (Example \ref{l3.1}), the method we used to prove Proposition \ref{p2.1} should also fail in higher degrees. The interesting question is, where does it fail? In case of Lemma \ref{l3.1} it turns out that the set $M=\left\{(f(x),g(x)):x\in \mathbb{R}^2\right\}$ is not a convex subset of $\mathbb{R}^2$ anymore. We reprise that $g=\mathrm{x}_1^2-\mathrm{x}_2^2$ and $f=\mathrm{x}_1^2(\mathrm{x}_1^2-\mathrm{x}_2^2)$.
Consider the half-line $H=\left\{\lambda(-1,1):\lambda\in \mathbb{R}_{\geq 0}\right\}$.
Since $S(g)\subseteq S(f)$ holds, the intersection $M\cap H$ can only consist of one point and this point is the origin $(0,0)$.
On the other hand, we can find two points $(x_1,x_2),(y_1,y_2)\in M$ such that $x_1>0$, $x_2>0$, $y_1<0$, $y_2<0$ and the line segment $L$ connecting the points $(x_1,x_2)$, $(y_1,y_2)$ does not go through the origin. This means $L$ intersects $H$ in some other point than the origin.
But this intersection point cannot be in $M$. Thus $M$ is not convex.
\end{Rem}
\end{window}   
\begin{window}[3, r, \includegraphics[scale=0.3]{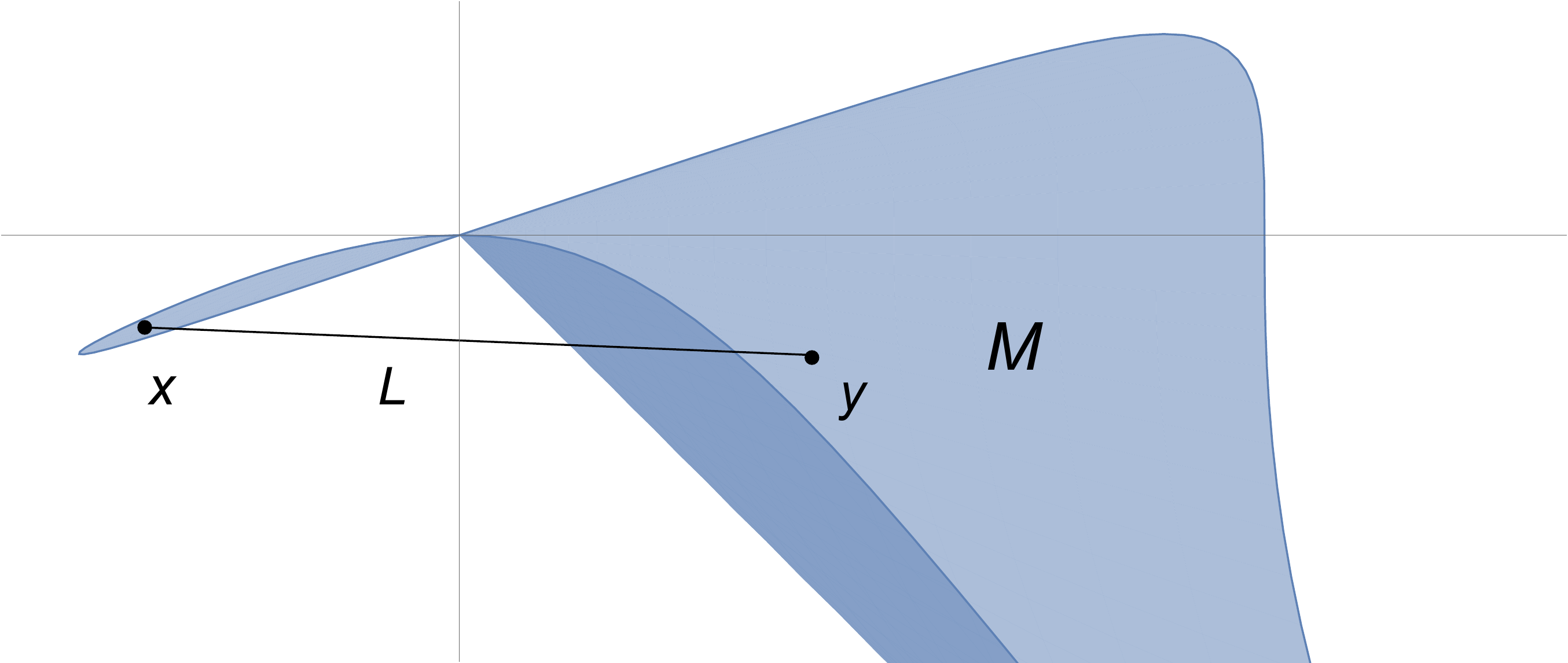},{}]
\begin{Exp}
Consider the two homogeneous polynomials $q=\mathrm{x}_1^3\mathrm{x}_2-\mathrm{x}_1^2\mathrm{x}_2^2$ and $p=\mathrm{x}_1^4-\mathrm{x}_1\mathrm{x}_2^3$.
We want to show that $M=\left\{(p(x),q(x)):x\in \mathbb{R}^2\right\}$ is not convex, giving an example where $M$ fails to be convex even if both homogeneous polynomials have the same degree. Set $H=\left\{(0,x_2):x_2<0\right\}$.
There are points $(x_1,x_2),(y_1,y_2)\in M$ with $x_1<0,x_2<0,y_1>0,y_2<0$ such that the segment line $L$, connecting both points, doesn't go through the origin. So the intersection point of $L$ and $H$ is not in $M$. Thus $M$ is not convex.
\end{Exp}
\end{window}
\section{Formulating a higher degree S-lemma}
In the past section we have seen that there is no way to increase the degree in the homogeneous S-lemma (Proposition \ref{p2.1}) and keep all other statement as they are. 
There are two ways one can proceed. One way would be to left the conditions (a) and (b) in the homogeneous S-lemma unchanged, and find some additional statements that might be plugged in into the homogeneous S-lemma such that the S-lemma remains true. Results in this sense can be found in \cite{z11}.\\\\
But there is another way. Instead of keeping the statements (a) and (b) as they are, we could simply modify the condition (b) by giving up that $t$ should be a non-negative real number. We demand that $t$ should be a non-negative homogeneous polynomial. The advantage is that we do not need to make up some new statements. The disadvantage, however, is that we have less information about $t$. This philosophy motivates to formulate 
\begin{Con}\label{k3.c1}
\textbf{S4-Conjecture}: Let $f$ be a tenary quartic, that is a $4$-form in $\mathbb{R}[\mathrm{x}_1,\mathrm{x}_2,\mathrm{x}_3]$, and let $g$ be a quadratic form in $\mathbb{R}[\mathrm{x}_1,\mathrm{x}_2,\mathrm{x}_3]$. Suppose there exists a point $x'\in \mathbb{R}^3$ such that $g(x')>0$. Then the following statements are equivalent:
\begin{enumerate}[label=(\alph*)]
\item The inclusion $S(g)\subseteq S(f)$ holds.
\item There exists a non-negative homogeneous polynomial $t\in \mathbb{R}[\mathrm{x}_1,\mathrm{x}_2,\mathrm{x}_3]_2$ such that $f(x)-t(x)g(x)\geq 0$ for all $x\in \mathbb{R}^3$.
\end{enumerate}
\end{Con}
\noindent The statement \ref{k3.c1} is originally a question posted in \textit{Mathoverflow} and was answered by the author. See \cite{z16}. 
\begin{Rem}
The homogeneous polynomial $t$ in the S4-conjecture is either a quadratic form or the zero-polynomial, which can be interpreted as a homogeneous polynomial of negative-infinite degree. The case $t=0$ can only occur if $f$ is a non-negative form. If $f$ is not a non-negative form, then $t$ must be of degree 2. This becomes clear if we take a point $x\notin S(f)$. 
Then we get $f(\lambda x)-t(\lambda x)g(\lambda x)=\lambda^4f(x)-\lambda^2t(x)g(x)\geq 0$ for all $\lambda\geq 0$ (*). Note that if $t$ is not of degree $2$, then it must be a constant.
But no matter what sign the constant $t$ has, we will always get $\lambda^4f(x)-\lambda^2t(x)g(x)\rightarrow -\infty$ for $\lambda\rightarrow \infty$, which contradicts (*).\par\smallskip  
Since $g$ is a quadratic form, we can take a look at the signature of $g$. If the S4-conjecture fails, then the next obvious question is: Is there a counterexample for all non-trivial signatures $-1,0$ and $1$. Note that if the signature of $g$ is $-2$ resp. $2$, then $g$ is non-positive resp. non-negative. But it is obvious how to deal with the S4-conjecture if $g$ is non-positive or non-negative.  
\end{Rem}
\section{S4-conjecture in two variables}
In this section we are going to prove the S4-conjecture in just two variables aka. 
\begin{Theo}\label{k4.t1}
\textbf{S4-conjecture in 2-variables}: Let $f$ be a $4$-form in $\mathbb{R}[\mathrm{x}_1,\mathrm{x}_2]$ and let $g$ be a quadratic form in $\mathbb{R}[\mathrm{x}_1,\mathrm{x}_2]$. Suppose there exists a point $x'\in \mathbb{R}^2$ such that $g(x')>0$. Then the following statements are equivalent:
\begin{enumerate}[label=(\alph*)]
\item The inclusion $S(g)\subseteq S(f)$ holds.
\item There exists a non-negative homogeneous polynomial $t\in \mathbb{R}[\mathrm{x}_1,\mathrm{x}_2]_2$ such that $f(x)-t(x)g(x)\geq 0$ for all $x\in \mathbb{R}^2$.
\end{enumerate}
\end{Theo}
\begin{Rem}
Theorem \ref{k4.t1} is not true, without the condition that there exists a point $x'\in \mathbb{R}^2$ such that $g(x')>0$! Consider for example the polynomials $g=-\mathrm{x}_1^2$ and $f=\mathrm{x}_1^4+\mathrm{x}_1^3\mathrm{x}_2+\mathrm{x}_1\mathrm{x}_2^3$. Obviously $g$ is a non-positive quadratic form and the equation $g(x_1,x_2)=0$ holds if and only if $x_1=0$. But for $x_1=0$ we have $f(0,x_2)=0$, which implies $S(g)\subseteq S(f)$. Take a homogeneous polynomial $t=a_1\mathrm{x}_1^2+a_2\mathrm{x}_2^2+a_3\mathrm{x}_1\mathrm{x}_2$, where $a_1,a_2,a_3\in \mathbb{R}$. Then $f(\mathrm{x}_1,1)-t(\mathrm{x}_1,1)g(\mathrm{x}_1,1)=\mathrm{x}_1^4+\mathrm{x}_1^3+\mathrm{x}_1+a_1\mathrm{x}_1^4+a_2\mathrm{x}_1^2+a_3\mathrm{x}_1^3$. But no matter what the coefficients $a_1,a_2,a_3$ are, the polynomial $f(\mathrm{x}_1,1)-t(\mathrm{x}_1,1)g(\mathrm{x}_1,1)$ has a sign change at the origin. Thus the implication (a)$\Rightarrow$(b) in Theorem \ref{k4.t1} would be false.   
\end{Rem} 
\begin{lem}\label{k5.l1}
Let $p$ and $q$ be polynomials in $\mathbb{R}[\mathrm{x}]$ such that $\deg(p)=4$, $\deg(q)=2$ and $S(q)\subseteq S(p)$. Suppose that there exists a point $x'\in \mathbb{R}$ with $q(x')>0$. Then there exists a non-negative polynomial $t\in \mathbb{R}[\mathrm{x}]$ of degree at most $2$ such that $p(x)-t(x)q(x)\geq 0$ for all $x\in \mathbb{R}$.
\end{lem}
\textbf{Proof}: Without loss of generality we can assume that neither $p$ nor $q$ are non-negative polynomials. Otherwise we could just take $t=0$.
For the sake of simplicity and oversight we will devide the proof in several meaningful cases:\par\smallskip  
\underline{Case I: The polynomial $p$ has a real double root $y\in \mathbb{R}$}: In this situation $p$ is divided by $s=(\mathrm{x}-y)^2\in \mathbb{R}[\mathrm{x}]$. Set $h=\frac {p} {s}\in \mathbb{R}[\mathrm{x}]$. It is clear that $h$ is a polynomial of degree $2$ and that the inclusion $S(q)\subseteq S(h)$ holds: 
The only situation in which $S(q)\subseteq S(h)$ might fail is the one, where $y$ is an isolated point of $S(q)$. But this would imply that $s|q$, which means that $q$ is either a non-negative or a non-positive polynomial. Obviously non of these two cases can occur.
According to the S-lemma (Theorem \ref{t1}) there is a non-negative real number $t'\geq 0$ such that $h(x)-t'q(x)\geq 0$ for all $x\in \mathbb{R}$. This implies that the inequality $0\leq s(x)(h(x)-t'q(x))=p(x)-t's(x)q(x)$ holds for all $x\in \mathbb{R}$. Set $t=t's$ and the assertion follows.\par\smallskip  
\underline{Case II: The polynomial $p$ has a complex root $c\in \mathbb{C}\backslash \mathbb{R}$:} If $c$ is a complex root of $p$, then $p$ is divided by the polynomial $s=(\mathrm{x}-c)(\mathrm{x}-\overline{c})\in \mathbb{R}[\mathrm{x}]$. Since $s$ has only complex roots, we conclude that $s$ is nowhere changing its sign. It is easy to see that $\lim_{x\rightarrow \infty}s(x)=\infty$. Thus $s$ is non-negative. By defining $h=\frac {p} {s}$ and repeating the procedure in case I, we get the desired result.\par\smallskip  
\underline{Case III: The polynomial $p$ has four distinct real roots $x_1<x_2<x_3<x_4$}: Here we     
have two possibilities how $q$ may actually look like. If there exists a point $c\in [x_1,x_2]$ such that $p(c)<0$, then $q=\alpha(\mathrm{x}-\tilde{x}_1)(\mathrm{x}-\tilde{x}_4)$ for $\alpha>0$ and $\tilde{x}_1\leq x_1$, $\tilde{x}_4\geq x_4$. If there exists a point $c\in [x_1,x_2]$ such that $p(c)>0$, then $q=\alpha(\mathrm{x}-\tilde{x}_1)(\mathrm{x}-\tilde{x}_4)$ for $\alpha<0$ and  $\tilde{x}_1\geq x_1$, $\tilde{x}_4\leq x_4$.
\begin{figure}[h]
\begin{minipage}[t]{0.5\textwidth}
\caption{First possibility}
\begin{tikzpicture}[scale=0.5]
\draw[blue ,line width=1.8pt,shorten >= 3pt,shorten <= 3pt]
(-4,0) to[out=-80, in=180] (-2,-2) to[bend right] (0,0)to[out=80, in=-180] (2,2) to [bend left] (4,0) to[out=-80, in=-180] (6,-2) to[bend right] (8,0);
\draw[dashed ,line width=1.8pt,shorten >= 3pt,shorten <= 3pt]
(-2,0) to (-2,-2);
\draw[dashed ,line width=1.8pt,shorten >= 3pt,shorten <= 3pt]
(-5,0) to (0,0);
\draw[dashed,blue ,line width=1.8pt,shorten >= 3pt,shorten <= 3pt]
(0,0) to (4,0);
\draw[dashed ,line width=1.8pt,shorten >= 3pt,shorten <= 3pt]
(4,0) to (9,0);
\filldraw
(-4,0) circle (4pt) node[align=left, above] {$x_1$};
\filldraw
(-2,0) circle (4pt) node[align=left, above] {$c$};
\filldraw
(0,0) circle (4pt) node[above right] {$x_2$};
\filldraw
(4,0) circle (4pt) node[above right] {$x_3$};
\filldraw
(8,0) circle (4pt) node[above] {$x_4$};
\draw[red ,line width=1.8pt,shorten >= 3pt,shorten <= 3pt]
(-5,0) to[bend right] (2,-4) to[bend right] (9,0);
\filldraw
(-5,0) circle (4pt) node[align=left, above] {$\tilde{x}_1$};
\filldraw
(9,0) circle (4pt) node[align=left, above] {$\tilde{x}_4$};
\filldraw
(2,2) circle (0pt) node[above right] {$\mathrm{graph}(p)$};
\filldraw
(2,-4) circle (0pt) node[below right] {$\mathrm{graph}(q)$};
\end{tikzpicture}
\end{minipage}
\begin{minipage}[t]{0.5\textwidth}
\caption{Second possibility}
\begin{tikzpicture}[scale=0.5]
\draw[blue ,line width=1.8pt,shorten >= 3pt,shorten <= 3pt]
(-4,0) to[out=80, in=180] (-2,2) to[bend left] (0,0)to[out=-80, in=-180] (2,-2) to [bend right] (4,0) to[out=80, in=-180] (6,2) to[bend left] (8,0);
\draw[dashed,line width=1.8pt,shorten >= 3pt,shorten <= 3pt]
(-2,0) to (-2,2);
\draw[dashed,blue ,line width=1.8pt,shorten >= 3pt,shorten <= 3pt]
(-4,0) to (0,0);
\draw[dashed,line width=1.8pt,shorten >= 3pt,shorten <= 3pt]
(0,0) to (4,0);
\draw[dashed,blue ,line width=1.8pt,shorten >= 3pt,shorten <= 3pt]
(4,0) to (8,0);
\filldraw
(-4,0) circle (4pt) node[align=left, below] {$x_1$};
\filldraw
(-2,0) circle (4pt) node[align=left, below] {$c$};
\filldraw
(0,0) circle (4pt) node[below right] {$x_2$};
\filldraw
(4,0) circle (4pt) node[below left] {$x_3$};
\filldraw
(8,0) circle (4pt) node[below] {$x_4$};
\draw[red ,line width=1.8pt,shorten >= 3pt,shorten <= 3pt]
(-3,0) to[bend left] (-2,1) to[bend left] (-1,0);
\draw[red ,line width=1.8pt,shorten >= 3pt,shorten <= 3pt]
(5,0) to[bend left] (6,1) to[bend left] (7,0);
\filldraw
(-3,0) circle (4pt) node[align=left, below] {$\tilde{x}_1$};
\filldraw
(-1,0) circle (4pt) node[align=left, below] {$\tilde{x}_4$};
\filldraw
(5,0) circle (4pt) node[align=left, below] {$\tilde{x}_1$};
\filldraw
(7,0) circle (4pt) node[align=left, below] {$\tilde{x}_4$};
\end{tikzpicture}
\end{minipage}
\end{figure}
Consider the first possibility. Define $h=\alpha(\mathrm{x}-x_1)(\mathrm{x}-x_4)$ and $s_1=\frac {p'(x_1)} {h'(x_1)}$, $s_4=\frac {p'(x_4)} {h'(x_4)}$. Note that $h'$ does not vanish at $x_1$ resp. $x_4$. Thus $s_1$ and $s_4$ are well defined positive real numbers. 
Let $v\in \mathbb{R}[\mathrm{x}]$ be a positive polynomial of degree 2 such that $v(x_1)=s_1$ and $v(x_4)=s_4$. For example, consider the polynomial $v=a(\mathrm{x}-x_1)^2+s_1$ with $a=\frac {s_4-s_1} {(x_4-x_1)^2}$ if $s_1\leq s_4$ resp. the polynomial $v=a(\mathrm{x}-x_4)^2+s_4$ with $a=\frac {s_1-s_4} {(x_1-x_4)^2}$ if $s_1>s_4$.  
The polynomial $w=p-vh$ has two double roots, namely one double root at $x_1$ and the other at $x_4$. This proves that $w$ is either non-negative or non-positive. Since $w(x_3)>0$, the polynomial $w$ is indeed non-negative.
It is easy to see that the inclusion $S(q)\subseteq S(h)$ holds. According to the S-lemma there exists a real non-negative number $t'$ such that $h(x)-t'q(x)\geq 0$ holds for all $x\in \mathbb{R}$. This implies $-h(x)\leq -t'q(x)$ for all $x\in \mathbb{R}$ and therefore $w(x)\leq p(x)-t'v(x)q(x)$. Since $w$ is non-negative, we are done. The second possibility is considered nearly analogous:
Without loss of generality we can assume that $x_1\leq \tilde{x}_1$, $x_2\geq \tilde{x}_4$ and $h=\alpha(\mathrm{x}-x_1)(\mathrm{x}-x_2)$.
As before we define $s_1=\frac {p'(x_1)} {h'(x_1)}$ resp. $s_2=\frac {p'(x_2)} {h'(x_2)}$. Let $v$ be a positive polynomial of degree $2$ such that $v(x_1)=s_1$ and $v(x_2)=s_2$ are satisfied. The polynomial $w=p-vh$ has two double roots at $x_1$ and $x_2$. Since $w(x_4)>0$, we see that $w$ is non-negative. As before we can deduce from $S(q)\subseteq S(h)$ that there exists a positive real number $t'$ with $h(x)-t'q(x)\geq 0$ for all $x\in \mathbb{R}$, implying that $w(x)\leq p(x)-t'v(x)q(x)$ for all $x\in \mathbb{R}$. Set $t=t'v$ and we are done. It is clear that there are no cases left that can occur. Thus the lemma is proven.\bew
\begin{lem}\label{k5.l2x}
Let $q\in \mathbb{R}[\mathrm{x}_1,\ldots,\mathrm{x}_n]$ be a quadratic form such that there exists a point $x'\in \mathbb{R}^n$ with $q(x')>0$. For every point $x\in S(q)$ and for every $\varepsilon>0$ the intersection between $B_{\varepsilon}(x)$ and $\mathrm{int}(S(q))$ is non-empty.
\end{lem}      
\textbf{Proof}: By using an appropriate change of coordinates we can rewrite $q$ as $q=a_1\mathrm{x}_1^2+\cdots+a_n\mathrm{x}_n^2$, where $a_1,\ldots,a_n\in \mathbb{R}$.
Define $I=\left\{i\in \{1,\ldots,n\}:a_i>0\right\}$. Note that $I\neq \varnothing$ because of $q(x')>0$. Take a point $x\in \mathbb{R}^n$ with $q(x)\geq 0$.
Choose an index $j\in I$, a real positive number $\varepsilon>0$ and consider $q(x_1,\ldots,x_j+\varepsilon,\ldots,x_n)$ if $x_j\geq 0$ resp. $q(x_1,\ldots,x_j-\varepsilon,\ldots,x_n)$ if $x_j<0$. It is clear that $q(x_1,\ldots,x_j+\varepsilon,\ldots,x_n)$ resp. $q(x_1,\ldots,x_j-\varepsilon,\ldots,x_n)$ is positive.   
This means that the point $y=\begin{cases}
(x_1,\ldots,x_j+\varepsilon,\ldots,x_n),\,\text{if}\,x_j\geq 0\\
(x_1,\ldots,x_j-\varepsilon,\ldots,x_n),\,\text{if}\,x_j<0\\
\end{cases}$ is lying in both, the interior of $S(q)$ and the ball $B_\varepsilon(x)$. In other words, the assertion is proven.\bew  
\begin{lem}\label{k5.l2}
Let $p\in \mathbb{R}[\mathrm{x}_1,\ldots,\mathrm{x}_n]$ be a polynomial of even degree $m\in \mathbb{N}_0$ and $q\in \mathbb{R}[\mathrm{x}_1,\ldots,\mathrm{x}_n]$ a polynomial of degree $2$. Suppose further that there is a point $x'\in \mathbb{R}^n$ such that $q(x')>0$. Let $\overline{p},\overline{q}\in \mathbb{R}[\mathrm{x}_1,\ldots,\mathrm{x}_{n+1}]$ denote their homogenizations. If $S(q)\subseteq S(p)$ then $S(\overline{q})\subseteq S(\overline{p})$.  
\end{lem}
\textbf{Proof}: 
We have to show that $S(\overline{q}(\mathrm{x},0))\subseteq S(\overline{p}(\mathrm{x},0))$ holds. Without loss of generality we can assume that $S(\overline{q}(\mathrm{x},0))\neq \varnothing$.
Note that $q(x')>0$ implies $\overline{q}(x',1)>0$. Thus we can use Lemma \ref{k5.l2x}. 
Let $c$ be an arbitrary point in $S(\overline{q}(\mathrm{x},0))$. The task is to verify $c\in S(\overline{p}(\mathrm{x},0))$. Lemma \ref{k5.l2x} tells us that $B_{\varepsilon}(c)\cap \mathrm{int}(S(\overline{q}))\neq \varnothing$ for $\varepsilon>0$. For every $\varepsilon>0$ we can find a point $y_\varepsilon\in B_{\varepsilon}(c)\cap \mathrm{int}(S(\overline{q}))$ such that the $n+1$-th component of $y_\varepsilon$ does not vanish. We have $\overline{q}(y_{\varepsilon,1},\ldots,y_{\varepsilon,n+1})=\frac {1} {y_{\varepsilon,n+1}^2}\overline{q}\left(\frac {y_{\varepsilon,1}} {y_{\varepsilon,n+1}},\ldots,\frac {y_{\varepsilon,n}} {y_{\varepsilon,n+1}},1\right)$, where $\overline{q}\left(\frac {y_{\varepsilon,1}} {y_{\varepsilon,n+1}},\ldots,\frac {y_{\varepsilon,n}} {y_{\varepsilon,n+1}},1\right)>0$. Therefore $\overline{p}\left(\frac {y_{\varepsilon,1}} {y_{\varepsilon,n+1}},\ldots,\frac {y_{\varepsilon,n}} {y_{\varepsilon,n+1}},1\right)\geq 0$. Since $m$ is even, we get $\overline{p}(y_{\varepsilon,1},\ldots,y_{\varepsilon,n+1})=\frac {1} {y_{\varepsilon,n+1}^m}\overline{p}\left(\frac {y_{\varepsilon,1}} {y_{\varepsilon,n+1}},\ldots,\frac {y_{\varepsilon,n}} {y_{\varepsilon,n+1}},1\right)\geq 0$.  
This implies that $\mathrm{dist}(c,S(\overline{p}))<\varepsilon$. By making $\varepsilon>0$ arbitrary small and using the fact that $S(\overline{p})$ is a closed set, we conclude that $c\in S(\overline{p})$.\bew\\\\
\textbf{Proof of Theorem \ref{k4.t1} aka. S4-conjecture in 2 variables}: Without loss of generality we can restrict ourselves to the case where $f$ and $g$ are both not non-negative.\par\smallskip    
(b)$\Rightarrow$(a): Trivial.\par\smallskip  
(a)$\Rightarrow$(b): Since we are talking about quadratic forms, and $g$ is a quadratic form, it helps to take a look at its diagonal-form. Let us take a matrix $A\in \mathrm{O}_2$ and consider the induced map $\psi:\mathbb{R}[\mathrm{x}_1,\mathrm{x}_2]\rightarrow \mathbb{R}[\mathrm{x}_1,\mathrm{x}_2],p\mapsto p\circ A$. We can choose $A$ in such a way that $\psi(g)$ is in diagonal form. Applying $\psi$ on $f$ and $g$ does not mess up the prerequisites of Theorem \ref{k4.t1}. Thus we can assume that $g$ is already in diagonal from. Since $g$ is neither non-negative nor non-positive, two real numbers $a_{11},a_{22}\neq 0$ with $\mathrm{sgn}(a_{11})\neq \mathrm{sgn}(a_{22})$ can be found such that $g=a_{11}\mathrm{x}_1^2+a_{22}\mathrm{x}_2^2$. Furthermore, we can demand that $a_{11}>0$ and $a_{22}<0$ because otherwise we could apply the coordinate transformation $\mathbb{R}^2\rightarrow \mathbb{R}^2,(x_1,x_2)\mapsto (x_2,x_1)$. The proof is devided into two cases:\par\smallskip  
\underline{Case I: $\deg_{\mathrm{x}_1}(f)$ and $\deg_{\mathrm{x}_2}(f)<4$}: First of all, we show that the monomial $\mathrm{x}_1^3\mathrm{x}_2$ cannot appear in $f$, while the monomials $\mathrm{x}_1^2\mathrm{x}_2^2$ and $\mathrm{x}_1\mathrm{x}_2^3$ must appear in $f$. Suppose the monomial $\mathrm{x}_1^3\mathrm{x}_2$ would appear in $f$. Consider the polynomial $f(\mathrm{x}_1,1)\in \mathbb{R}[\mathrm{x}_1]$. Since $S(g)\subseteq S(f)$, we get the conditions $\lim_{x_1\rightarrow \infty}f(x_1,1)=\infty$ and $\lim_{x_1\rightarrow -\infty}f(x_1,1)=\infty$. But the leading term of $f(\mathrm{x}_1,1)$ is of the form $\alpha \mathrm{x}_1^3$ for $\alpha\neq 0$.\\
\begin{tikzpicture}[scale=0.5]
\draw[line width=2pt] (-8,4) -- (8,-4);
\draw[line width=2pt] (-8,-4) -- (8,4) ;
\draw[red, line width=2pt, dashed] (-8,1) -- (-4,1) node[anchor=north] {$S(g(\mathrm{x}_1,1))$} --(-2,1);
\draw[red, line width=2pt, dashed] (2,1) -- (8,1);
\draw[blue, line width=2pt, dashed] (2,1) -- (2,-1) node[anchor=west] {$S(g(1,\mathrm{x}_2))$};
\draw[fill=gray, opacity=0.3] (-8,4) to (-8,-4) to (0,0);
\draw[fill=gray, opacity=0.3] (8,-4) to (8,4) to (0,0);
\fill (5,0) circle (0pt) node[anchor=west] {$S(g)$};
\end{tikzpicture}\\   
Thus $f_1(\mathrm{x}_1,1)$ cannot fulfill the two conditions and therefore we get a contradiction. On the other side, we cannot exclude the monomial $\mathrm{x}_1\mathrm{x}_2^3$. In case of the monomial $\mathrm{x}_2\mathrm{x}_1^3$ we exploited that the set $S(g(\mathrm{x}_1,1))$ is symmetric and unbounded. This is not the case when we consider $S(g(1,\mathrm{x}_2))$. Instead we can state that since $f$ is non-negative on the set $S(g)$, the polynomial $f$ cannot consist of just one monomial $\mathrm{x}_1\mathrm{x}_2^3$ or $\mathrm{x}_1^2\mathrm{x}_2^2$ alone. So, we can rewrite $f$ as $f=\gamma \mathrm{x}_1^2\mathrm{x}_2^2+\beta \mathrm{x}_1\mathrm{x}_2^3$ with $\gamma>0$ and $\beta\in \mathbb{R}\backslash \{0\}$.
Define $s=-\frac {1} {2}ba_{22}\mathrm{x}_2^4+\frac {1} {2}ba_{11}\mathrm{x}_1^2\mathrm{x}_2^2+\beta \mathrm{x}_1\mathrm{x}_2^3$ and $t=\frac {1} {2}b\mathrm{x}_2^2$ where $a_{11}b=\gamma>0$. A simple computation shows $f=tg+s$. We are done, if we can show that $s$ is a non-negative polynomial. Since $s$ is divided by $\mathrm{x}_2^2$, it is sufficient to prove that $s'=\frac {s} {\mathrm{x}_2^2}=-\frac {1} {2}ba_{22}\mathrm{x}_2^2+\beta \mathrm{x}_1\mathrm{x}_2+\frac {1} {2}ba_{11}\mathrm{x}_1^2$ is non-negative. The discriminant of $s'$ is given by $\mathrm{disc}(s')=\left(\beta^2+a_{22}a_{11}b^2\right)\mathrm{x}_1^2$. Then we have the equivalence $\mathrm{disc}(s')\leq 0$ for all $x_1\in \mathbb{R}$ $\Leftrightarrow$ $\beta^2+a_{22}a_{11}b^2\leq 0$.
It is sufficient to show that $\beta^2+a_{22}a_{11}b^2\leq 0$.    
To prove this, take a point $y\in \partial S(g)$ such that $y_1,y_2>0$. Furthermore we assume that $\beta< 0$. Now $y\in \partial S(g)$ implies $y_1=y_2\sqrt{\left|\frac {a_{22}} {a_{11}}\right|}$ and thus $f(y)=\left(\gamma\left|\frac {a_{22}} {a_{11}}\right|+\beta \sqrt{\left|\frac {a_{22}} {a_{11}}\right|}\right)y_2^4\geq 0$. This is only possible if $\gamma\left|\frac {a_{22}} {a_{11}}\right|+\beta \sqrt{\left|\frac {a_{22}} {a_{11}}\right|}\geq 0$. The inequality $\gamma\left|\frac {a_{22}} {a_{11}}\right|+\beta \sqrt{\left|\frac {a_{22}} {a_{11}}\right|}\geq 0$ is equivalent to $\gamma\sqrt{\left|\frac {a_{22}} {a_{11}}\right|}\geq |\beta|$. By substituting $\gamma$ through $a_{11}b$, we get $\sqrt{|a_{11}a_{22}|}b\geq |\beta|$ and finally $\beta^2+a_{11}a_{22}b^2\leq 0$, because $a_{11}a_{22}<0$.\\
If $\beta>0$ then we simply consider another point $y\in \partial S(q)$ with $y_1>0$, $y_2<0$, and repeat the arguments above.\par\smallskip  
\underline{Case II: $\deg_{\mathrm{x}_1}(f)=4$ or $\deg_{\mathrm{x}_2}(f)=4$}: We are only considering the case $\deg_{\mathrm{x}_1}(f)=4$. Let $\tilde{f}\in \mathbb{R}[\mathrm{x}_1]$ be the dehomogenization of the polynomial $f$ in the variable $\mathrm{x}_2$.
In the same manner, let $\tilde{g}\in \mathbb{R}[\mathrm{x}_1]$ be the dehomogenization of $g$. It is easy to see that $\tilde{f}$ and $\tilde{g}$ fulfill the prerequisites of Lemma \ref{k5.l1}. Thus there exists a non-negative polynomial $\tilde{t}\in \mathbb{R}[\mathrm{x}_1]$ with $\deg(\tilde{t})\leq 2$ such that $\tilde{f}(x_1)-\tilde{t}(x_1)\tilde{g}(x_1)\geq 0$ for all $x_1\in \mathbb{R}$. 
If $\tilde{f}(x_1)-\tilde{t}(x_1)\tilde{g}(x_1)=0$ holds for all $x_1\in \mathbb{R}$, then the assertion follows immediately. 
Otherwise, Lemma \ref{k5.l2} tells us that $S(\tilde{f}-\tilde{t}\tilde{g})=\mathbb{R}$ and $S(\tilde{t})=\mathbb{R}$ imply $S\left(\overline{\tilde{f}-\tilde{t}\tilde{g}}\right)=\mathbb{R}^2$ and $S(t)=\mathbb{R}^2$, where $t=\begin{cases}
\overline{\tilde{t}}\,\text{if}\,\deg(\tilde{t})=2\\
\mathrm{x}_2^2\tilde{t}\,\text{if}\,\deg(\tilde{t})=0
\end{cases}$.
If $\deg(\tilde{t})=2$ then we get that $f-tg=\mathrm{x}_2^4\left(\overline{\tilde{f}-\tilde{t}\tilde{g}}\right)$, $f-tg=\mathrm{x}_2^2\left(\overline{\tilde{f}-\tilde{t}\tilde{g}}\right)$ or
$f-tg=\overline{\tilde{f}-\tilde{t}\tilde{g}}$. This follows directly from the fact that $\overline{\tilde{f}-\tilde{t}\tilde{g}}$ is of even degree and that $\deg(f)=\deg(tg)$.
If $\deg(\tilde{t})=0$ then it is easy to see that $f-tg=\overline{\tilde{f}-\tilde{t}\tilde{g}}$. Thus $f-tg$ is non-negative.\bew\newpage
\section{The S4-conjecture: A counterexample}
It turns out that the S4-conjecture stated in \ref{k3.c1} is wrong. First of all, we are going to give a 'lucky' counterexample and afterwards, that means in the next chapter, we will investigate why the S4-conjecture cannot work. 
\begin{Exp}\label{k6.e1x}
\textbf{A counterexample for \ref{k3.c1}}: Consider the polynomials $f=\mathrm{x}_1^3\mathrm{x}_3+\mathrm{x}_1^3\mathrm{x}_2+\mathrm{x}_2^2\mathrm{x}_3^2$ and $g=\mathrm{x}_1\mathrm{x}_3+\mathrm{x}_2\mathrm{x}_3+\mathrm{x}_1\mathrm{x}_2$. One can show that the inclusion $S(g)\subseteq S(f)$ holds: For example, type in 
\footnotesize
\begin{verbatim}
Reduce[ForAll[{x1,x2,x3},Implies[x1*x3+x2*x3+x1*x2>=0,x1^3*x3+x1^3*x2+x2^2*x3^2>=0]]]
\end{verbatim}
\normalsize
in \cite{z15}.\\
But there is no non-negative homogeneous polynomial $t\in \mathbb{R}[\mathrm{x}_1,\mathrm{x}_2,\mathrm{x}_3]$ of degree at most $2$ such that $f(y)-t(y)g(y)\geq 0$ for all $y\in \mathbb{R}^3$: Suppose we could find such a polynomial $t=a_1\mathrm{x}_1^2+a_2\mathrm{x}_2^2+a_3\mathrm{x}_1\mathrm{x}_2+a_4\mathrm{x}_1\mathrm{x}_3+a_5\mathrm{x}_2\mathrm{x}_3+a_6\mathrm{x}_3^2$, where $a_1,\ldots,a_6\in \mathbb{R}$. A simple computation shows that $f-tg=\mathrm{x}_1^3\mathrm{x}_2-a_1\mathrm{x}_1^3\mathrm{x}_2-a_3\mathrm{x}_1^2\mathrm{x}_2^2-a_2\mathrm{x}_1\mathrm{x}_2^3+\mathrm{x}_1^3\mathrm{x}_3-a_1\mathrm{x}_1^3\mathrm{x}_3-a_1\mathrm{x}_1^2\mathrm{x}_2\mathrm{x}_3-a_3\mathrm{x}_1^2\mathrm{x}_2\mathrm{x}_3-a_4\mathrm{x}_1^2\mathrm{x}_2\mathrm{x}_3-a_2\mathrm{x}_1\mathrm{x}_2^2\mathrm{x}_3-a_3\mathrm{x}_1\mathrm{x}_2^2\mathrm{x}_3-a_5\mathrm{x}_1\mathrm{x}_2^2\mathrm{x}_3-a_2\mathrm{x}_2^3\mathrm{x}_3-a_4\mathrm{x}_1^2\mathrm{x}_3^2-a_4\mathrm{x}_1\mathrm{x}_2\mathrm{x}_3^2-a_5\mathrm{x}_1\mathrm{x}_2\mathrm{x}_3^2-a_6\mathrm{x}_1\mathrm{x}_2\mathrm{x}_3^2+\mathrm{x}_2^2\mathrm{x}_3^2-a_5\mathrm{x}_2^2\mathrm{x}_3^2-a_6\mathrm{x}_1\mathrm{x}_3^3-a_6\mathrm{x}_2\mathrm{x}_3^3$.\par\smallskip  
Thus $f(\mathrm{x}_1,0,1)-t(\mathrm{x}_1,0,1)g(\mathrm{x}_1,0,1)=\mathrm{x}_1^3-a_1\mathrm{x}_1^3-a_4\mathrm{x}_1^2-a_6\mathrm{x}_1$. We know that $f-tg$ is a non-negative polynomial. This is only possible if $a_1=1$, $a_6=0$, and $a_4\leq 0$. Since $t$ is non-negative and $a_6=0$, the two coefficients $a_4$, $a_5$ must also vanish. The leading term of $f(0,\mathrm{x}_2,1)-t(0,\mathrm{x}_2,1)g(0,\mathrm{x}_2,1)$ is $-a_2\mathrm{x}_2^3$, and therefore $a_2=0$. Now only $a_3$ is not determined. But it is easy to see that $a_3$ must also vanish. Thus $f-tg$ reduces to $f-tg=-\mathrm{x}_1^2\mathrm{x}_2\mathrm{x}_3+\mathrm{x}_2^2\mathrm{x}_3^2$, which is obviously not non-negative.    
\end{Exp}
\begin{Rem}
Under an appropriate linear change of coordinates, we can rewrite $g$ as $g=\mathrm{x}_1^2-\frac {1} {2}\mathrm{x}_2^2-\frac {1} {2}\mathrm{x}_3^2$. Thus under this new coordinates $S(g)$ is a double cone and every slice with a plane, parallel to the $x_2,x_3$-plane, is compact. Fix $c>0$ and consider $S:=S\left(c^2-\frac {1} {2}\mathrm{x}_2^2-\frac {1} {2}\mathrm{x}_3^2\right)$, which is a compact subset of $\mathbb{R}^2$. A simple computation shows that $f(c,x_2,x_3)>0$ for all $x_2,x_3\in S$. Using the Positivstellensatz of Schm\"udgen, we see that $f(c,\mathrm{x}_2,\mathrm{x}_3)\in T(g(c,\mathrm{x}_2,\mathrm{x}_3))$, where $T(g(c,\mathrm{x}_2,\mathrm{x}_3))$ denotes the preordering generated by $g(c,\mathrm{x}_2,\mathrm{x}_3)$. This, however, is not true for $f$ and $g$, i.e $f\notin T(g)$. See Lemma \ref{l.count} and Remark \ref{remT}.\newpage      
\end{Rem}
\section{Geometric analysis}
In this section we are going to take a closer look at the counterexample. In particular we are interested in the geometric properties of $V_1=\mathcal{V}(f)$ and $V_2=\mathcal{V}(g)$ and what they have to do with the counterexample. In this section we will fix $f=\mathrm{x}_1^3\mathrm{x}_3+\mathrm{x}_1^3\mathrm{x}_2+\mathrm{x}_2^2\mathrm{x}_3^2$ and $g=\mathrm{x}_1\mathrm{x}_3+\mathrm{x}_2\mathrm{x}_3+\mathrm{x}_1\mathrm{x}_2$. 
\begin{lem}\label{l.count}
Let $f$ and $g$ be as in Example \ref{k6.e1}. Then there is no non-negative homogeneous polynomial $t\in \mathbb{R}[\mathrm{x}_1,\mathrm{x}_2,\mathrm{x}_3]$ of even degree $n$ such that $f(y)-t(y)g(y)\geq 0$ for all $y\in \mathbb{R}^3$. 
\end{lem}
\textbf{Proof}:
Without loss of generality we can assume that $n>2$. 
Let $t$ be a non-negative homogeneous polynomial in $\mathbb{R}[\mathrm{x}_1,\mathrm{x}_2,\mathrm{x}_3]$ of even degree $n>2$. We are going to show that $f(\mathrm{x}_1,\mathrm{x}_2,1)-t(\mathrm{x}_1,\mathrm{x}_2,1)g(\mathrm{x}_1,\mathrm{x}_2,1)$ is not a non-negative polynomial in $\mathbb{R}[\mathrm{x}_1,\mathrm{x}_2]$, which is a stronger statement than that in the lemma. Because $f(\mathrm{x}_1,\mathrm{x}_2,0)=\mathrm{x}_1^3\mathrm{x}_2$ is not non-negative, we can assume that $t(\mathrm{x}_1,\mathrm{x}_2,0)\neq 0$. Thus we have $\deg(t(\mathrm{x}_1,\mathrm{x}_2,0))=n$.\\  
Write $t(\mathrm{x}_1,\mathrm{x}_2,1)=\sum_{\alpha\in \mathbb{N}^2,|\alpha|\leq n}c_\alpha \mathrm{x}^\alpha$ and define $I=\left\{\alpha\in \mathbb{N}^2:|\alpha|=n,c_\alpha\neq 0\right\}$. Note that $I$ is not empty, since $\deg(t(\mathrm{x}_1,\mathrm{x}_2,0))=n$. 
Without loss of generality we can assume that there is an element $\alpha\in I$ such that $\alpha_2>0$. Otherwise, we could simply interchange the variables $\mathrm{x}_1$ and $\mathrm{x}_2$.
Let $\alpha '$ be the uniquely determined element of $I$ that satisfies $\alpha_2'>\alpha_2$ for all $\alpha\in I\backslash \{\alpha '\}$. There exists a real number $\beta>0$ such that $c_{\alpha '}\beta^{\alpha_2 '}+\sum_{\alpha\in I\backslash \{\alpha'\}}c_\alpha\beta^{\alpha_2}\neq 0$: Indeed, $\sum_{\alpha\in I}c_\alpha\upbeta^{\alpha_2}$ is a polynomial in $\mathbb{R}[\upbeta]$ of degree $\alpha_2'>0$ and therefore does not vanish.  
In fact, $\sum_{\alpha\in I}c_\alpha\beta^{\alpha_2}$ is the leading coefficient and $\sum_{\alpha\in I}c_\alpha\beta^{\alpha_2}\mathrm{x}_1^{\alpha_1+\alpha_2}=\sum_{\alpha\in I}c_\alpha\beta^{\alpha_2}\mathrm{x}_1^n$ the leading term of the polynomial $t(\mathrm{x}_1,\beta \mathrm{x}_1,1)$. 
So, the real positive number $\beta$ is needed to make sure that this leading term does not vanish.
Because $t(\mathrm{x}_1,\beta \mathrm{x}_2,1)$ is non-negative, the coefficient $\sum_{\alpha\in I}c_\alpha\beta^{\alpha_2}$ is positive. The leading term of $tg$ is $\left(\sum_{\alpha\in I}c_\alpha \beta^{\alpha_2+1}\right)\mathrm{x}_1^{n+2}$ with a positive coefficient $\sum_{\alpha\in I}c_\alpha \beta^{\alpha_2+1}$. This implies that the leading term of $f-tg$ is $-\left(\sum_{\alpha\in I}c_\alpha \beta^{\alpha_2+1}\right)\mathrm{x}_1^{n+2}$. Therefore we get
$$\lim_{x_1\rightarrow \infty}\left(f(x_1,\beta x_1,1)-t(x_1,\beta x_1,1)g(x_1,\beta x_1,1)\right)=-\infty,$$ which proves the lemma.\bew
\begin{Rem}\label{remT}
Lemma \ref{l.count} combined with the Counterexample \ref{k6.e1x} proves that $f\notin T(g)$. Suppose we could find sums of squares $\sigma_1$ and $\sigma_2$ in $\mathbb{R}[\mathrm{x}_1,\mathrm{x}_2,\mathrm{x}_3]$ such that $f=\sigma_1+\sigma_2 g$. Then we have $\tilde{f}=\tilde{\sigma}_1+\tilde{\sigma}_2\tilde{g}$, where we dehomogenize with respect to $\mathrm{x}_3$. Without loss of generality we can assume that $\deg(\tilde{\sigma}_2)\geq 2$. We distinguish between two cases:
\begin{itemize}
\item We have $\deg\left(\tilde{\sigma}_2\tilde{g}\right)=\deg(\tilde{f})=4$: Under this condition, we have $f-\overline{\tilde{\sigma}_2}g=\mathrm{x}_3^n\left(\overline{\tilde{f}-\tilde{\sigma}_2\tilde{g}}\right)$, where $n\leq 4$ is a even number.
By using Lemma \ref{k5.l2} we see that $\overline{\tilde{f}-\tilde{\sigma}_2\tilde{g}}$ is non-negative. 
Thus $f-\overline{\tilde{\sigma}_2}g$ is non-negative. Since $\overline{\tilde{\sigma}_2}$ is a non-negative polynomial of degree $2$, we get a contradiction. 
\item  We have $\deg\left(\tilde{\sigma}_2\tilde{g}\right)>\deg(\tilde{f})=4$: In this case we proceed as in Lemma \ref{l.count}: Choose a real positive number $\beta$ such that the leading monomial of $f(\mathrm{x}_1,\beta\mathrm{x}_1,1)-\sigma_2(\mathrm{x}_1,\beta\mathrm{x}_1,1)g(\mathrm{x}_1,\beta\mathrm{x}_1,1)$ is $-L\beta\mathrm{x}_1^2$, where $L$ denotes the leading monomial of $\sigma_2(\mathrm{x}_1,\beta\mathrm{x}_1,1)$. Since $\sigma_2(\mathrm{x}_1,\beta\mathrm{x}_1,1)$ is a sum of squares in $\mathbb{R}[\mathrm{x}_1]$, the polynomial $L$ is a sum of squares and therefore $L\beta\mathrm{x}_1^2$ is non-negative.
Thus $$\lim_{x_1\rightarrow \infty}\left(f(x_1,\beta x_1,1)-\sigma_2(x_1,\beta x_1,1)g(x_1,\beta x_1,1)\right)\rightarrow -\infty,$$ which contradicts our assumption. 
\end{itemize}
\end{Rem}
\begin{lem}\label{kx.l0}
The two $\mathbb{R}$-varieties $\mathcal{V}(g)$ and $\mathcal{V}(f)$ are both geometrically irreducible. Furthermore, the set $H=\left\{(0,x_2,0):x_2\in \mathbb{R}\right\}$ is the set of all $\mathbb{R}$-rational singularities of $\mathcal{V}(g)$ and $\mathcal{V}(f)$.  
\end{lem}
\textbf{Proof}: First, $\mathcal{V}(f)$ resp. $\mathcal{V}(g)$ is irreducible if and only if $\mathcal{V}(\tilde{f})$ resp. $\mathcal{V}(\tilde{g})$ is irreducible, where $\tilde{f}=\mathrm{x}_1^3+\mathrm{x}_1^3\mathrm{x}_2+\mathrm{x}_2^2$ and $\tilde{g}=\mathrm{x}_1+\mathrm{x}_2+\mathrm{x}_1\mathrm{x}_2$. This statement is a well known fact.
Consider the polynomial $\tilde{g}$ as an element of the ring $\mathbb{C}[\mathrm{x}_1][\mathrm{x}_2]$. Furthermore, $\tilde{g}$ is a primitive polynomial. According to \cite[Satz 2, p. 68]{z18} we know that $\tilde{g}$ is irreducible in $\mathbb{C}[\mathrm{x}_1][\mathrm{x}_2]$ if the image of $\tilde{g}$ in $\left(\mathbb{C}[\mathrm{x}_1]/\mathrm{x}_1\mathbb{C}[\mathrm{x}_1]\right)[\mathrm{x}_2]$ is irreducible, which is easy to verify. 
Hence $\mathcal{V}(\tilde{g})$ is an irreducible $\mathbb{C}$-variety.
Let us consider $\tilde{f}$ as an element of the polynomial ring $\mathbb{C}[\mathrm{x}_1][\mathrm{x}_2]$. 
The polynomial $\tilde{f}$ cannot be divided by any irreducible polynomial in $\mathbb{C}[\mathrm{x}_1]$: Indeed, $\mathbb{C}[\mathrm{x}_1][\mathrm{x}_2]/\mathrm{x}_1\mathbb{C}[\mathrm{x}_1][\mathrm{x}_2]\cong \left(\mathbb{C}[\mathrm{x}_1]/\mathrm{x}_1\mathbb{C}[\mathrm{x}_1]\right)[\mathrm{x}_2]$ and the image of $\tilde{f}$ in $\left(\mathbb{C}[\mathrm{x}_1]/\mathrm{x}_1\mathbb{C}[\mathrm{x}_1]\right)[\mathrm{x}_2]$ is not a unit.
Suppose an irreducible factor $h$ of $\tilde{f}$ has the same degree in $\mathrm{x}_2$ as $\tilde{f}$. Then $h$ must coincide with $\tilde{f}$: If there would exist another irreducible factor $v$, then $v$ must lie in $\mathbb{C}[\mathrm{x}_1]$. But this would be a contradiction, since $v\nmid \tilde{f}$. It is impossible that $\tilde{f}$ factors into more than one component in $\mathbb{C}[\mathrm{x}_1][\mathrm{x}_2]$:
Suppose we could write $\tilde{f}=h_1h_2$, where $h_1,h_2\in \mathbb{C}[\mathrm{x}_1][\mathrm{x}_2]$ are polynomials of degree $1$ in $\mathrm{x}_2$. The polynomials $h_1$ and $h_2$ can be written as $h_1=\mathrm{x}_2-v_1$ and $h_2=\mathrm{x}_2-v_2$, where $v_1,v_2\in \mathbb{C}[\mathrm{x}_1]$. Then $\tilde{f}=\mathrm{x}_2^2-\mathrm{x}_2v_2-\mathrm{x}_2v_1+v_1v_2$. Therefore $v_1v_2=\mathrm{x}_1^3$ and $-\mathrm{x}_2(v_1+v_2)=\mathrm{x}_2\mathrm{x}_1^3$, which is utterly impossible.   
Hence $\mathcal{V}(\tilde{f})$ is an irreducible $\mathbb{C}$-variety.
It remains to verify the statement about the singularities. Consider a point $x\in \mathcal{V}(\tilde{f})$ such that $\nabla \tilde{f}(x)=0$. Then we have $\nabla f(x')=0$ for $x'=(x,1)$. On the other hand, suppose there is a point $x'\in V_1$  with $x_3'\neq 0$ such that $\nabla f(x')=0$. Then the point $x=\left(\frac {x_1'} {x_3'},\frac {x_2'} {x_3'}\right)$ will satisfy $\nabla \tilde{f}(x)=0$.
This shows that the singular points $x'\in V_1$ with $x_3'\neq 0$ 'come from' the singular points of $\mathcal{V}(\tilde{f})$. Thus it suffices to show that $\mathcal{V}(\tilde{f})$ has only one $\mathbb{R}$-rational singularity at the origin and that all other $\mathbb{R}$-rational singularities of $V_1$ are in $H$.\\
The equation $\nabla\tilde{f}(x)=\left(\begin{array}{c}
3x_1^2+3x_1^2x_2\\
x_1^3+2x_2
\end{array}\right)=0$ has only one real solution $x=(0,0)$ in $\mathcal{V}(\tilde{f})$, proving the first assertion of the last statement.\\ 
Finally, consider the equation $\nabla f(x_1,x_2,0)=\left(\begin{array}{c}
3x_1^2x_2\\
x_1^3\\
x_1^3\end{array}\right)=0$, where the set of all solutions in $\mathbb{R}^3$ is exactly $H$. Note that $H$ is a subset of $V(f)(\mathbb{R})$ and $V(g)(\mathbb{R})$. The same argumentation applied to $g$ gives the same result. Thus the lemma is proven.\bew
\begin{Rem}
A standard theorem in algebraic geometry states that a irreducible variety $V$ over $\mathbb{C}$ is connected with respect to the norm topology. For a proof see \cite[Theorem 1, p. 126]{z3}. In \ref{kx.l0} we proved that $V_2=V(g)$ is irreducible. Thus $V_2\subseteq \mathbb{C}^2$ is a connected set. But it is easy to see that $V_2(\mathbb{R})$ is not connected. This implies that \cite[Theorem 1, p. 126]{z3} is not true if we just consider the $\mathbb{R}$-rational points.   
\end{Rem}
\begin{lem}\label{kx.l1}
Let $(q_n)_n$ be a convergent sequence of quadratic forms in $\mathbb{R}[\mathrm{x}_1,\ldots,\mathrm{x}_n]$ and $(p_n)_n$ a convergent sequence of forms in $\mathbb{R}[\mathrm{x}_1,\ldots,\mathrm{x}_n]$ of degree $d$ such that $S(q_n)\subseteq S(p_n)$ for every $n\in \mathbb{N}$. If $q$ and $p$ are the limits of the sequences $(q_n)_n$ and $(p_n)_n$, then $S(q)$ is a subset of $S(p)$ if $q\neq 0$.
\end{lem}          
\textbf{Proof}: 
First of all we are going to prove the assertion under the assumption that $\mathrm{int}(S(q))\neq \varnothing$. 
Under an appropriate change of coordinates we can assume that $q=a_1\mathrm{x}_1^2+\cdots+a_n\mathrm{x}_n^2$, where $a_1,\ldots,a_n\in \mathbb{R}$.
For every point $y\in S(q)$ there exists a sequence $(y_n)_n\subset \mathrm{int}(S(q))$ such that $\lim_{n\rightarrow \infty}y_n=y$.
Since $\mathrm{int}(S(q))$ is not empty, we can handle this statement with Lemma \ref{k5.l2x}.
Consider a point $x$ lying in the interior of $S(q)$. Then there exists a number $N\in \mathbb{N}$ such that $q_n(x)>0$ for all $n\geq N$, implying that $\lim_{n\rightarrow \infty}q_n(x)\geq 0$. Since $S(q_n)$ is a subset of $S(p_n)$, we get $\lim_{n\rightarrow \infty}p_n(x)\geq 0$ resp. $x\in S(p)$.
If $x$ lies in $\partial S(q)$, then there exists a sequence $(x_n)\subset \mathrm{int}(S(q))$ such that $\lim_{n\rightarrow\infty}x_n=x$. But we showed above that this sequence also lies in $S(p)$. Thus $x$ lies in $S(p)$, since $S(p)$ is closed.\par\smallskip  
Finally, the case $\mathrm{int}(S(q))=\varnothing$ must be considered. This is only possible if the coefficients of $q=a_1\mathrm{x}_1^2+\cdots+a_n\mathrm{x}_n^2$ are all non-positive and at least one of them is negative. If all coefficients are negative, then we are done, since any form in $\mathbb{R}[\mathrm{x}_1,\ldots,\mathrm{x}_n]$ is non-negative at the origin. Therefore we can assume that not all coefficients are negative. Thus the set $I'=\left\{i\in \{1,\ldots,n\}:a_i=0\right\}$ is not empty. Set $H=\prod_{i\notin I'}\mathbb{R}\times \prod_{i\in I'}\{0\}$, $H'=\prod_{i\in I'}\mathbb{R}\times \prod_{i\notin I'}\{0\}$ and consider $q'=q|_H$, $p'=p|_H$. Then we have $S(q')=\{0\}$ by assumption and therefore $S(q')\subseteq S(p')$. Since $S(q)=S(q')\cup H'$, we have to make sure that $p''=p|_{H'}$ is non-negative. Consider $q_n''=q_n|_{H'}$, $q''=q|_{H'}$ and $p_n''=p_n|_{H'}$ as polynomials in $\mathbb{R}\left[\mathrm{x}_i:i\in I'\right]$. By using the facts that $S(q'')$ has a non-empty interior and $S(q_n'')\subseteq S(p_n'')$ for all $n\in \mathbb{N}$, we can apply the result made in the first part to deduce that $S(p'')\supseteq S(q'')$. Thus the lemma is proven.\bew 
\begin{Prop}\label{kx.p1}
Let \\
Set $S=\left\{(q,p)\in \mathbb{R}[\mathrm{x}_1,\mathrm{x}_2,\mathrm{x}_3]_2\times \mathbb{R}[\mathrm{x}_1,\mathrm{x}_2,\mathrm{x}_3]_4,q\,\text{quadratic form},\,p\,\text{4-form}\right\}$ and let $S_4$ be the set of all $(q,p)\in S$ that satisfy the following condition:
\begin{itemize}
\item There exists a non-negative homogeneous polynomial $t\in \mathbb{R}[x_1,x_2,x_3]_2$ such that $p(y)-t(y)q(y)\geq 0$ for all $y\in \mathbb{R}^3$. 
\end{itemize} 
Then the set $S_4$ is a closed subset of $S$.
\end{Prop}          
\textbf{Proof}: Let $P_{4,3}\subset \mathbb{R}[\mathrm{x}_1,\mathrm{x}_2,\mathrm{x}_3]_4$ be the set of all non-negative 4-forms and $P_{2,3}\subset \mathbb{R}[\mathrm{x}_1,\mathrm{x}_2,\mathrm{x}_3]_2$ the set of all non-negative quadratic forms. It is well known that $P_{4,3}$ and $P_{2,3}$ are closed cones (see Proposition \ref{p11}).\par\smallskip  
Let $(q_n,p_n)_n\subset S_4$ be a convergent sequence in $S$. For every $n\in \mathbb{N}$ there is a $t_n\in P_{2,3}$ such that $p_n-t_nq_n\in P_{4,3}$. Or in other words, there exists a sequence $(t_n)_n\subset P_{2,3}$ such that $(p_n-t_nq_n)_n\subset P_{4,3}$. Since $P_{2,3}$ and $P_{4,3}$ are closed, we get $t=\lim_{n\rightarrow \infty}t_n\in P_{2,3}$ and $\lim_{n\rightarrow \infty}(p_n-t_nq_n)=p-tq\in P_{4,3}$, where $p=\lim_{n\rightarrow \infty}p_n\in P_{4,3}$ and $q=\lim_{n\rightarrow\infty}q_n\in P_{2,3}$.\par\smallskip  
Lemma \ref{kx.l1} tells us that $S(q)$ is a subset of $S(p)$ if $q\neq 0$. 
If $q=0$ then $f-tq\in P_{4,3}$ implies that $f\in P_{4,3}$, which leads straight to $S(q)=S(f)=\mathbb{R}^3$.
Altogether we have that $(q,p)\in S_4$ and therefore $S_4$ is a closed subset of $S$.\bew\\\\
Let us consider the following statement:    
\begin{Con}\label{k6.c1}
\textbf{Dehomogenized S4-Conjecture}: Let $f$ be a polynomial fo degree $4$ in $\mathbb{R}[\mathrm{x}_1,\mathrm{x}_2]$ and let $g$ be a polynomial of degree $2$ in $\mathbb{R}[\mathrm{x}_1,\mathrm{x}_2]$. Suppose there exists a point $x'\in \mathbb{R}^2$ such that $g(x')>0$. Then the following statements are equivalent:
\begin{enumerate}[label=(\alph*)]
\item The inclusion $S(g)\subseteq S(f)$ holds.
\item There exists a non-negative polynomial $t\in \mathbb{R}[\mathrm{x}_1,\mathrm{x}_2]_2$ such that $f(x)-t(x)g(x)\geq 0$ for all $x\in \mathbb{R}^2$.
\end{enumerate}
\end{Con}
Note that Lemma \ref{k5.l2} makes sure that if we find a counterexample for \ref{k6.c1} we have a counterexample for the original S4-conjecture by homogenization:\\ 
If $p,q\in \mathbb{R}[\mathrm{x}_1,\mathrm{x}_2,\mathrm{x}_3]$ satisfy the S4-conjecture, then $\tilde{p}\in \mathbb{R}[\mathrm{x}_1,\mathrm{x}_2]$ and $\tilde{q}\in \mathbb{R}[\mathrm{x}_1,\mathrm{x}_2]$ will satisfy the dehomogenized S4-conjecture. Suppose $\tilde{p},\tilde{q}\in \mathbb{R}[\mathrm{x}_1,\mathrm{x}_2]$ satisfy the dehomogenized S4-conjecture. Then $\tilde{p}-\tilde{t}\tilde{q}$ is non-negative. Lemma \ref{k5.l2} tells us that $\overline{\tilde{p}-\tilde{t}\tilde{q}}$ and $t:=\begin{cases}
\overline{\tilde{t}},\,\text{if}\,\deg(\tilde{t})=2\\
\mathrm{x}_3^2\tilde{t},\,\text{if}\,\deg(\tilde{t})=0
\end{cases}$ are non-negative. Let $p$ and $q$ be the homogenizations of $\tilde{p}$ and $\tilde{q}$. Then we have $f-tg=\mathrm{x}_3^4\left(\overline{\tilde{p}-\tilde{t}\tilde{q}}\right)$, $f-tg=\mathrm{x}_3^2\left(\overline{\tilde{p}-\tilde{t}\tilde{q}}\right)$ or $f-tg=\overline{\tilde{p}-\tilde{t}\tilde{q}}$, which implies the non-negativity of $p-tq$. 
It is easy to see that the polynomials $\tilde{f}=\mathrm{x}_1^3+\mathrm{x}_1^3\mathrm{x}_2+\mathrm{x}_2^2$ and $\tilde{g}=\mathrm{x}_1+\mathrm{x}_2+\mathrm{x}_1\mathrm{x}_2$ form a counterexample for \ref{k6.c1}. What is the point with the dehomogenized S4-conjecture?\\\\
We want to use Proposition \ref{kx.p1} to prove that for a small $\varepsilon>0$ the two homogeneous polynomials $f_\varepsilon=f+\varepsilon \mathrm{x}_3^4$ and $g_\varepsilon=g+\varepsilon \mathrm{x}_3^2$ form still a counterexample to the S4-conjecture. First of all, we must make sure that $S(g_\varepsilon)\subseteq S(f_\varepsilon)$. But it is easy to see that $S\left(\tilde{g}_\varepsilon\right)\subseteq S\left(\tilde{f}_\varepsilon\right)$ holds for small $\varepsilon>0$. Using Lemma \ref{k5.l2} we can deduce that $S(g_\varepsilon)\subseteq S(f_\varepsilon)$.\par\smallskip  
By using Proposition \ref{kx.p1} and shrinking $\varepsilon>0$ further if necessary, we can achieve $(g_\varepsilon,f_\varepsilon)\notin S_4$. Thus we are getting a counterexample for the S4-conjecture.\par\smallskip  
But if we consider the geometry of $\mathcal{V}(g_\varepsilon)$ and $\mathcal{V}(f_\varepsilon)$, then not much has changed. A simple computation shows that these two varieties have the same $\mathbb{R}$-rational singularities as $\mathcal{V}(g)$ and $\mathcal{V}(f)$. But $\mathcal{V}(\tilde{g}_\varepsilon)$ and $\mathcal{V}(\tilde{f}_\varepsilon)$ have not just no $\mathbb{R}$-rational singularities, they are indeed non-singular varieties. While the geometric situation has not change in the homogeneous situation, the situation concerning the dehomogenized S4-conjecture is obviously different. But as we mentioned at the beginning of this investigation, $\tilde{f}$ and $\tilde{g}$ resp. $\tilde{f}_\varepsilon$ and $\tilde{g}_\varepsilon$ will fail the dehomogenized S4-conjecture. The point with dehomogenized S4-conjecture is, that it reflects the geometric differences between $\tilde{f},\tilde{g}$ and $\tilde{f}_\varepsilon,\tilde{g}_\varepsilon$ in a better way than its homogeneous counterpart.
In fact this result tells us, that the reason $f,g$ and $f_\varepsilon,g_\varepsilon$ are failing the S4-conjecture must lie in some other geometric properties.        
\section{A generalization of the counterexample}
Before continuing our investigation, it is worth to prove a generalization of the counterexample \label{k6.e1}. We can extend this result by using some simple results in algebraic geometry.
Let us consider the polynomials $f=\mathrm{x}_1^3+\mathrm{x}_1^3\mathrm{x}_2+\mathrm{x}_2^2$ and $g=\mathrm{x}_1+\mathrm{x}_2+\mathrm{x}_1\mathrm{x}_2$ which are the dehomogenizations of the polynomials in Example \ref{k6.e1}. Instead of introducing new symbols for the dehomogenization of the polynomials in Example \ref{k6.e1}, we will refer to them with the same symbols in this section.   
\begin{Def}\label{k6.d1}
(\textbf{Blow-up of} $\mathbb{A}^2$): Let $Y$ be the variety that is defined by the equation $x_1z_2=x_2z_1$, $(x_1,x_2;z_1:z_2)\in \mathbb{A}^2\times \mathbb{P}^1$. The restriction $\sigma:Y\rightarrow \mathbb{A}^2$ of the projection $\mathbb{A}^2\times \mathbb{P}^1\rightarrow \mathbb{A}^2$ onto $\mathbb{A}^2$ is called the blow-up of $\mathbb{A}^2$ centered at the origin. 
\end{Def}
\begin{Rem}
Definition \ref{k6.d1} might suggest that a blow-up centered at a regular point $x$ of a variety is unique in nature. But this is not the fact. In case of Definition \ref{k6.d1} it is true. But in case of an arbitrary quasi-projective variety $X$, where $X$ is not projective, we have just uniqueness up to 'isomorphism' \cite[Lemma, p. 117]{z2}. 
\end{Rem}
\begin{Prop}\label{pblow}
Let $V$ be a irreducible curve in $\mathbb{A}^2$ and $\sigma:Y\rightarrow \mathbb{A}^2$ the blow-up of $\mathbb{A}^2$ centered at the origin. Consider the curve $V'=\overline{\sigma^{-1}\left(V\backslash \{0\}\right)}$, where the bar denotes the Zariski-closure of $\sigma^{-1}\left(V\backslash \{0\}\right)$. Then we have the following statements about $V$ and $V'$:
\begin{enumerate}[label=(\alph*)]
\item If $0\notin V$ then there is an isomorphism $V\DistTo V'$.
\item If $0\in V$ then $\sigma^{-1}(V)$ decomposes into two irreducible components $E=\{0\}\times \mathbb{P}^1$ (exceptional curve) and $V'$ (birational transform).
\end{enumerate} 
\end{Prop}
\textbf{Proof}: See \cite[Theorem 1, p.118]{z2}, though the statement is much more general.\bew\\\\
Let us return to our polynomials $f$ and $g$. Set $V_1=\mathcal{V}(f)$ and $V_2=\mathcal{V}(g)$. If we want to know more about the behaviour of $V_1$ resp. $V_2$ under the blow-up of $\mathbb{A}^2$ centered at the origin, we have to make sure that $f$ and $g$ are irreducible polynomials in $\mathbb{R}[\mathrm{x}_1,\mathrm{x}_2]$. But this has been done in Lemma \ref{kx.l0}.\par\smallskip  
Consider the variety $Y$ given by the equation $x_1z_2=x_2z_1$, $(x_1,x_2;z_1:z_2)\in \mathbb{A}^2\times \mathbb{P}^1$. Suppose $z_1\neq 0$. Then we can choose $z_1=1$ and therefore we are getting $x_2=x_1z_2$. In other words, we are considering $Y$ on the affine piece $\mathbb{A}^2\times \mathbb{A}^1$. Substituting $\mathrm{x}_2$ through $\mathrm{x}_1\mathrm{z}_2$ leads to 
$$f_1=\mathrm{z}_2^2\mathrm{x}_1^2+\mathrm{x}_1^3+\mathrm{z}_2\mathrm{x}_1^4=\underbrace{\mathrm{x}_1^2}_{E\text{-component}}\underbrace{\left(\mathrm{z}_2^2+\mathrm{x}_1+\mathrm{z}_2\mathrm{x}_1^2\right)}_{V_1'\text{-component}}$$ 
resp.
$$g_1=\mathrm{x}_1+\mathrm{z}_2\mathrm{x}_1+\mathrm{z}_2\mathrm{x}_1^2=\underbrace{\mathrm{x}_1}_{E\text{-component}}\underbrace{\left(1+\mathrm{z}_2+\mathrm{z}_2\mathrm{x}_1\right)}_{V_2'\text{-component}}.$$
We can make the following statement about $f_1$ and $g_1$:
\begin{Prop}\label{kx.pg}
There is no non-negative polynomial $t\in \mathbb{R}[\mathrm{x}_1,\mathrm{z}_2]$ such that $f_1(y)-t(y)g_1(y)\geq 0$ for all $y\in \mathbb{R}^2$.
\end{Prop}
\textbf{Proof}: Suppose there is such a polynomial $t$. Consider the homomorphism $\phi:\mathbb{R}[\mathrm{x}_1,\mathrm{z}_2]\rightarrow \mathbb{R}\left[\mathrm{x}_1,\frac {\mathrm{x}_2} {\mathrm{x}_1}\right],p(\mathrm{x}_1,\mathrm{z}_2)\mapsto p\left(\mathrm{x}_1,\frac {\mathrm{x}_2} {\mathrm{x}_1}\right)$. The next step is to verify that $\phi(t)$ is a polynomial in $\mathbb{R}[\mathrm{x}_1,\mathrm{x}_2]$. Write $\phi(t)=\sum_{i,j}a_{ij}\mathrm{x}_1^i\left(\frac {\mathrm{x}_2} {\mathrm{x}_1}\right)^j=\sum_{i<j}a_{ij}\mathrm{x}_1^i\left(\frac {\mathrm{x}_2} {\mathrm{x}_1}\right)^j+\sum_{j\leq i}a_{ij}\mathrm{x}_1^i\left(\frac {\mathrm{x}_2} {\mathrm{x}_1}\right)^j$. Suppose $\phi(t)\notin \mathbb{R}[\mathrm{x}_1,\mathrm{x}_2]$, i.e $\sum_{i<j}a_{ij}\mathrm{x}_1^i\left(\frac {\mathrm{x}_2} {\mathrm{x}_1}\right)^j\neq 0$. 
Choose $(i',j')\in \mathbb{N}^2_0$ such that the following conditions are satisfied:
\begin{itemize}
\item We have $j'>i'$ and $a_{i'j'}\neq 0$.
\item Define $\Delta_{ij}=(j'-i')-(j-i)$, where $j,i\in \mathbb{N}_0$.
The inequality $\Delta_{ij}\geq 0$ holds for all $(i,j)\in \mathbb{N}^2_0$ with $a_{ij}\neq 0$ and $j>i$. 
\item We have $j'>j$ for all $(i,j)\in \mathbb{N}^2_0$ with $j>i$, $a_{ij}\neq 0$ and $\Delta_{ij}=0$.
\end{itemize}
By using $\mathrm{x}_1^{j'-i'}$ as a common denominator for $\sum_{i<j}a_{ij}\mathrm{x}_1^i\left(\frac {\mathrm{x}_2} {\mathrm{x}_1}\right)^j$ we get  
$$\sum_{i<j}a_{ij}\mathrm{x}_1^i\left(\frac {\mathrm{x}_2} {\mathrm{x}_1}\right)^j=\frac {a_{i'j'}\mathrm{x}_2^{j'}+\sum_{i<j,(i,j)\neq (i',j')}a_{ij}\mathrm{x}_2^j\mathrm{x}_1^{\Delta_{ij}}} {\mathrm{x}_1^{j'-i'}}.$$
Note that the polynomial 
$$h(0,\mathrm{x}_2):=a_{i'j'}\mathrm{x}_2^{j'}+\sum_{i<j,(i,j)\neq (i',j'),\Delta_{ij}=0}a_{ij}\mathrm{x}_2^j$$ 
has degree $j'>0$ in the variable $\mathrm{x}_2$.
Take a point $x'=(0,c)\in \mathbb{R}^2$ that satisfies $(0,c)\in \mathrm{int}(S(g))$. All points on the $x_2$-axis, but the origin, are inner points of $S(g)$. This means that $c$ is just a positive real number. Since $\deg_{\mathrm{x}_2}(h(0,\mathrm{x}_2))>0$, we can choose $c>0$ such that $h(x')\neq 0$. 
Consider the sequence $(x_n)_{n\in \mathbb{N}}=\left(\left(\frac {1} {n},c\right)\right)_{n\in \mathbb{N}}$.
Then the limit $\lim_{n\rightarrow \infty}\phi(t)(x_n)$ is not finite. While the nominator tends to a finite value $h(x')\neq 0$, the denominator tends to zero. Thus the limit cannot be finite. Since $t\in \mathbb{R}[\mathrm{x}_1,\mathrm{z}_2]$ is non-negative, we see that $\lim_{n\rightarrow \infty}\phi(t)(x_n)=\infty$. But then $$\lim_{n\rightarrow \infty}\left(\phi(f_1)(x_n)-\phi(t)(x_n)\phi(g_1)(x_n)\right)=\lim_{n\rightarrow \infty}\left(f(x_n)-\phi(t)(x_n)g(x_n)\right)=-\infty.$$ This contradicts $f_1(y)-t(y)g_1(y)\geq 0$ for all $y\in \mathbb{R}^2$.
We just proved $\phi(t)\in \mathbb{R}[\mathrm{x}_1,\mathrm{x}_2]$. We can deduce that $f(y)-\phi(t)(y)g(y)\geq 0$ for all $y\in \mathbb{R}^2$, since $f_1(y)-t(y)g_1(y)\geq 0$ for all $y\in \mathbb{R}^2$. But this contradicts Lemma \ref{l.count}. Thus the proposition is proven.\bew\\\\
Set $f_1'=\mathrm{z}_2^2+\mathrm{x}_1+\mathrm{z}_2\mathrm{x}_1^2$ and $g_1'=1+\mathrm{z}_2+\mathrm{z}_2\mathrm{x}_1$. Repeat the same procedure done to $V_1$ and $V_2$ with $\mathcal{V}(f_1')$ and $\mathcal{V}(g_1')$.
Therefore we get two new polynomials $f_2$ and $g_2$. 
Substituting $f_1'$ resp. $g_1'$ in $f_1=\mathrm{x}_1^2f_1'$ resp. $g_1=\mathrm{x}_1g_1'$ by $f_2$ and $g_2$ gives us the polynomials $\mathrm{x}_1^2f_2$ and $\mathrm{x}_1g_2$.  
So, what's the point in doing that? The degree of $f_1$ resp. $g_1$ compared to the degree of $f$ resp. $g$ has increased by one. The same is also true for $\mathrm{x}_1^2f_2$ and $\mathrm{x}_1g_2$ with respect to $f_1$ and $g_1$. 
Let $f_i'$ and $g_i'$ denote the equations of the birational transformations of $\mathcal{V}(f_{i-1}')$ and $\mathcal{V}(f_{i-1}')$ on the affine piece $\mathbb{A}^2\times \mathbb{A}^1$ for $i\geq 2$.
By repeating this procedure we get the polynomials $\mathrm{x}_1^2f_i=\mathrm{x}_1^3f_i'$, where $f_i'=\mathrm{x}_1^{2(i-2)+1}\mathrm{z}_i^2+1+\mathrm{z}_i\mathrm{x}_1^{2+i}$ and $\mathrm{z}_i=\mathrm{x}_1\mathrm{z}_{i-1}$ for $i\geq 2$. By using Proposition \ref{pblow} we can see immediately that the polynomial $f_i'$ is irreducible for $i\geq 2$. On the other hand, we get $\mathrm{x}_1g_i=\mathrm{x}_1g_i'$, where $g_i'=1+\mathrm{z}_i+\mathrm{z}_i\mathrm{x}_1^i$ for $i\geq 2$. Hence $\deg(\mathrm{x}_1^2f_i)=\deg(\mathrm{x}_1^2f_{i-1})+2$ and $\deg(\mathrm{x}_1g_i)=\deg(\mathrm{x}_1g_{i-1})+1$ for $i>2$. The only specifics we used about the polynomials $f_1$ and $g_1$ in the proof of Proposition \ref{kx.pg}, was that $f_1$ resp. $g_1$ emerged from $f$ resp. $g$ by blowing up $V_1$ resp. $V_2$ and that there is no non-negative polynomial $t\in \mathbb{R}[\mathrm{x}_1,\mathrm{x}_2]$ such that $f-tg$ is non-negative. Hence we can make the same statement with respect to the polynomials $\mathrm{x}_1^2f_i$ and $\mathrm{x}_1g_i$, where $i\geq 2$. Thus we get a counterexample for the dehomogenized S4-conjecture in higher degrees.
Finally, we can state: 
\begin{Prop}\label{blonk}
For any natural number $d\in \{2n:n\geq 4\}\cup \{4,5,6\}$ there are polynomials $f$ and $g$ in $\mathbb{R}[\mathrm{x}_1,\mathrm{x}_2]$ that satisfy the following statements:
\begin{itemize}
\item The degree of $f$ is $d$ and the degree of $g$ is $\nu(d)$, where $$\nu(d)=\begin{cases}
d-2,\,\text{if}\,d\leq 6\\
d-\frac {d-6} {2}-2,\,\text{if}\,d>6
\end{cases}.$$
\item There is a point $x'\in \mathbb{R}^2$ such that $g(x')>0$.
\item The inclusion $S(g)\subseteq S(f)$ holds.
\item There is no non-negative polynomial $t\in \mathbb{R}[\mathrm{x}_1,\mathrm{x}_2]$ such that $f(y)-t(y)g(y)\geq 0$ for all $y\in \mathbb{R}^2$. 
\end{itemize} 
\end{Prop}    
\noindent Of course, by applying the blow-up procedure to other counterexamples the result in Proposition \ref{blonk} can be refined. As a hint one could start with the polynomials $\mathrm{x}_1^2\mathrm{x}_2-\mathrm{x}_1^2+1$ and $-\mathrm{x}_1^2+\mathrm{x}_2$. But since a refinement of Proposition \ref{blonk} is not our aim, we will not further pursue it. 
\newpage     
\chapter{Quadratic modules and stability}
The aim of this chapter is to clarify the reasons why Example \ref{k6.e1x} does form a counterexample for the S4-conjecture. In the first part of this chapter we will introduce the necessary tools to answer this question. This tools will be based on the article \cite{z12}. Finally, the second part is meant to deal with the question mentioned at the beginning, by answering it through a geometric criterion.   
\section{Preliminaries}
The following definitions and theorems can be found in \cite{z12}. The aim is to provide a list of basic tools for later needs.  
\begin{Def}
A subset $M\subseteq \mathbb{R}[\mathrm{x}_1,\ldots,\mathrm{x}_n]=:A$ is called a quadratic module, if  
$1\in M$, $M+M\subseteq M$, and $A^2\cdot M\subseteq M$ holds, where $A^2$ denotes the set of squares in $A$ and $\Sigma A^2$ denotes the sum of squares in $A$. Furthermore, $\mathrm{QM}(f_1,\ldots,f_s)=\left\{\sigma_0+\sigma_1f_1+\cdots+\sigma_sf_s:\sigma_0,\ldots,\sigma_s\in \Sigma A^2\right\}$ is called the quadratic module generated by $f_1,\ldots,f_s\in A$.  
\end{Def}
\noindent Throughout this chapter $A$ will denote the polynomial ring $\mathbb{R}[\mathrm{x}_1,\ldots,\mathrm{x}_n]$.
\begin{Def}
Let $A=\bigoplus_{\gamma\in \Gamma}A_\gamma$ be a grading and let $M\subseteq A$ be a finitely generated quadratic module. $M$ is totally stable with respect to the grading if $\deg(f)\leq \deg(f+g)$ holds for all $f,g\in M$. This is equivalent to the fact that there are generators $f_1,\ldots,f_s$ of $M$ such that 
$$\deg(\sigma_jf_j)\leq \deg\left(\sum_i\sigma_if_i\right)$$
holds for all $\sigma_j\in \Sigma A^2$. Any finite set of generators of $M$ fulfills this condition then. 
\end{Def}
\begin{Def}
For $z\in \mathbb{Z}$ and $d\in \mathbb{Z}^n$ we define $$A^{(z)}_d:=\left\{\sum_{\delta\in \mathbb{N}^n,\langle z,\delta\rangle=d}c_\delta \mathrm{x}_1^{\delta_1}\cdots \mathrm{x}_n^{\delta_n}:c_\delta\in \mathbb{R}\right\}.$$
Then $$A=\bigoplus_{d\in \mathbb{Z}}A_d^{(z)}$$
is a grading that we will call the $z$-grading of $A$. For an element $f\in A$ we define $L_z(f)$ to be the degree component (component with the highest degree) of $f$ with respect to the $z$-grading of $A$. 
\end{Def}
\begin{Rem}
In the literature the polynomials that lie in $A^{(z)}_d$ are called quasi-homogeneous polynomials of type $z$ and degree $d$.
\end{Rem}
\begin{Def}\label{kxd1x}
For a compact set $K\subseteq \mathbb{R}^n$ with non-empty interior, we define the tentacle of $K$ in direction of $z\in \mathbb{Z}^n$ in the following way:
$$T_{K,z}:=\left\{\left(\lambda^{z_1}x_1,\ldots,\lambda^{z_n}x_n:\lambda \geq 1, x=(x_1,\ldots,x_n)\in K\right)\right\}.$$
\end{Def}
\begin{Theo}\label{kxtt11}
Let $f_1,\ldots,f_s$ be polynomials in the graded polynomial algebra $A=\bigoplus_{d\in \mathbb{Z}}A_d^{(z)}$, where $z\in \mathbb{Z}^n$. If the set $S(f_1,\ldots,f_s)\subseteq \mathbb{R}^n$ contains a tentacle $T_{K,z}$, then the quadratic module $M=\mathrm{QM}(f_1,\ldots,f_s)$ is totally stable with respect to the $z$-grading. If $M$ is closed under multiplication, then $S(f_1,\ldots,f_s)$ must contain such a tentacle for $M$ to be totally stable.   
\end{Theo}
\textbf{Proof}: See \cite{z12}[Theorem 5.2].\bew
\section{Stability and tentacles}
\begin{Def}
Let $q=\sum_{i,j}a_{ij}\mathrm{x}_i\mathrm{x}_j$ be a quadratic form in $A$. The diagonal part $\mathcal{D}(q)$ of $q$ is defined by $$\mathcal{D}(q)=\sum_ia_{ii}\mathrm{x}_i^2.$$
\end{Def}
\begin{Def}
Let $f$ be a polynomial in $A$. The set $\mathcal{T}_0(f)$ is defined to be the set of all $z\in \mathbb{Z}^n$, under which the quadratic module $\mathrm{QM}(f)$ is totally stable.  
\end{Def}
\begin{Prop}\label{rooters}
Let $f=\sum_{i=0}^na_i\mathrm{x}^i$ be a polynomial of degree $n>0$ in $\mathbb{C}[\mathrm{x}]$ with distinct roots $x_1,\ldots,x_r\in \mathbb{C}$, where $r\leq n$. Furthermore we define $f(\mathrm{x},y)=\sum_{i=0}^n(a_i+y_i)\mathrm{x}^i$ for a point $y\in \mathbb{C}^{n+1}$. For every $0<\varepsilon$ there exists a $\delta>0$ such that all distinct roots of $f(\mathrm{x},y)$ lie in $\bigcup_{i=1}^rB_{\varepsilon}(x_i)$ for all $\left\|y\right\|_2<\delta$.  
\end{Prop}
\textbf{Proof}: Define $h(\mathrm{x},y)=\frac {f'(\mathrm{x},y)} {f(\mathrm{x},y)}$ and $d=\begin{cases}
\frac {1} {2},\,\text{if}\,r=1\\
\frac {1} {2}\min_{i>j}|x_i-x_j|,\,\text{otherwise}
\end{cases}$. Without loss of generality we can assume that $0<\varepsilon<d$.  It is therefore easy to see that we can find a simple closed, null-homologous path $\Gamma_{i,y}$ in $B_\varepsilon(x_i)$ such that $x_i$ is in the interior of $\Gamma_{i,y}$ and $f(\mathrm{x},y)$ does not vanish on $\Gamma_{i,y}$. For example, choose $\Gamma_{i,y}$ to be a circle around $x_i$ such that $f(\mathrm{x},y)$ does not vanish on this circle. Take an arbitrary $i=1,\ldots,r$. According to a consequence of the residual theorem \cite{z13}[Proposition 13.2.3, p. 350] we have $$\frac {1} {2\pi i}\oint_{\Gamma_{i,y}}h(x,y)\mathrm{d}x=N_i(y),$$
where $N_i(y)$ denotes the number of roots with multiplicity of $f(\mathrm{x},y)$ in $B_{\varepsilon}(x_i)$.
It is easy to see that 
$$\left|\frac {1} {2\pi i}\oint_{\Gamma_{i,y}}h(x,y)\mathrm{d}x-\frac {1} {2\pi i}\oint_{\Gamma_{i,y}}h(x,0)\mathrm{d}x\right|\rightarrow 0$$
as $\left\|y\right\|_2\rightarrow 0$. Since both integrals are integer numbers, there exists a real positive number $\delta_i$ such that $\frac {1} {2\pi i}\oint_{\Gamma_{i,y}}h(x,y)\mathrm{d}x=\frac {1} {2\pi i}\oint_{\Gamma_{i,y}}h(x,0)\mathrm{d}x$ for all $\left\|y\right\|_2<\delta_i$. This implies that $N_i(y)=N_i(0)$, where $N_i(0)$ is simply the multiplicity of the root $x_i$ of $f$. Setting $\delta=\min\left\{\delta_i:i=1,\ldots,r\right\}$ concludes the proof.\bew
\begin{Cor}\label{crooters}
Suppose $f$ and $g$ are polynomials of degree $n_1>0$ resp. $n_2>0$ in $\mathbb{C}[\mathrm{x}]$ with distinct roots $x_1,\ldots,x_{r_1}\in \mathbb{C}$, where $r_1\leq n_1$ resp. $x_1',\ldots,x_{r_2}'\in \mathbb{C}$, where $r_2\leq n_2$. For every $y\in \mathbb{C}^{n+1}$ let $f(\mathrm{x},y)$ resp. $g(\mathrm{x},y)$ be defined as in the preceding proposition. Then for every $0<\varepsilon$ there exists a $\delta>0$ such that all roots and poles of $\frac {f(\mathrm{x},y)} {g(\mathrm{x},y)}$ lie in $\bigcup_{i=1}^{r_1}B_{\varepsilon}(x_i)\cup \bigcup_{i=1}^{r_2}B_{\varepsilon}(x_i')$ for all $\left\|\delta\right\|_2<\delta$. 
\end{Cor}    
\begin{Theo}\label{kxt1}
For a quadratic form $q\in A\backslash \{0\}$ the following statements hold:
\begin{enumerate}[label=(\alph*)]
\item Suppose $\mathcal{D}(q)$ is negative-definite. Then $\mathcal{T}_0(q)\subsetneqq \mathbb{Z}^n$.
\item If $\mathcal{D}(q)$ is non-negative, then $\mathcal{T}_0(q)=\mathbb{Z}^n$.
\end{enumerate}
\end{Theo}
\textbf{Proof}: 
We will write $q=\sum_{i,j}a_{ij}\mathrm{x}_i\mathrm{x}_j$ in this proof.\par\smallskip   
(a): We have to show that $\mathcal{T}_0(q)\subsetneqq \mathbb{Z}^n$. But this is quiet easy: Because $\mathcal{D}(q)$ is negative-definite there is a coefficient $a_{ii}$ for some $i\in \{1,\ldots,n\}$ such that $a_{ii}<0$. Take $z\in \mathbb{Z}^n$ such that all components but the $i$-th vanish and let the $i$-th component be a large positive number. It is clear that $z\notin \mathcal{T}_0(q)$.\par\smallskip  
(b): Suppose that $\mathcal{D}(q)=0$. We have to show that $\mathcal{T}_0(q)=\mathbb{Z}^n$. Because $\mathcal{T}_0(q)$ cannot be any bigger that $\mathbb{Z}^n$, it is enough to prove that the inclusion $\mathcal{T}_0(q)\supseteq \mathbb{Z}^n$ holds. 
For a given $z\in \mathbb{Z}^n$ let $I$ be the set of all $(i_1,j_1)\in \mathbb{N}^2$ with $a_{i_1j_1}+a_{j_1i_1}\neq 0$ such that there is no $(i_2,j_2)\in \mathbb{N}^2$ with $a_{i_2j_2}+a_{j_2i_2}\neq 0$ and $z_{i_2}+z_{j_2}>z_{i_1}+z_{j_1}$.
Take an arbitrary $z\in \mathbb{Z}^n\backslash \{0\}$ and take $(i,j)\in I$.
Without loss generality we can demand that $a_{ij}+a_{ji}>0$. Otherwise substitute the variable $\mathrm{x}_i$ through $-\mathrm{x}_i$ and define $\tilde{a}_{ij}=-a_{ij}$ resp. $\tilde{a}_{ji}=-a_{ji}$ as the new coefficient of $\mathrm{x}_i\mathrm{x}_j$ resp. $\mathrm{x}_j\mathrm{x}_i$. In the next step we prove that there exists a point $x'\in \mathbb{R}^n$ such that 
\begin{itemize}
\item $(a_{ij}+a_{ji})x_i'x_j'>0$
\item $q(x')>0$
\item $x_1',\ldots,x_n'\neq 0$.
\item $(a_{ij}+a_{ji})x_i'x_j'>2\left|\sum_{(i',j')\in I\backslash \{(i,j),(j,i)\}}a_{i'j'}x_{i'}'x_{j'}'\right|$
\end{itemize}    
Let us start with a point $x\in \mathbb{R}$ that satisfies $x_1,\ldots,x_n\neq 0$ and $\mathrm{sign}(x_i)=\mathrm{sign}(x_j)\neq 0$. Without loss of generality we can assume that $q(x)\leq 0$. Let us modify the point $x$. Consider $q(x_1,\ldots,\uplambda x_i,\ldots,\uplambda x_j,\ldots,x_n)\in \mathbb{R}[\uplambda]$. The leading term of this polynomial in $\uplambda$ is $(a_{ij}+a_{ji})x_ix_j\uplambda^2$.
Since $(a_{ij}+a_{ji})x_ix_j>0$, we have $q(x_1,\ldots,\lambda x_i,\ldots,\lambda x_j,\ldots,x_n)\rightarrow \infty$ for $\lambda\rightarrow \infty$. Choosing a large $\lambda\in \mathbb{R}$ and a new point $x'\in \mathbb{R}$ with $x'_k=x_k$ for $k\neq i,j$ and $x_i'=\lambda x_i$, $x_j'=\lambda x_j$ leads to $q(x')>0$. Finally, we can achieve $(a_{ij}+a_{ji})x_i'x_j'>2\left|\sum_{(i',j')\in I\backslash \{(i,j),(j,i)\}}a_{i'j'}x_{i'}'x_{j'}'\right|$ by enlarging $\lambda$ further if necessary.
Next, we want to find an appropriate point $x'\in \mathbb{R}^n$ and a neighborhood $U$ of $x'$ such that $T_{\overline{U},z}\subseteq S(q)$. Take $(i,j)\in I$ and a point $x'\in \mathbb{R}^n$ satisfying the four conditions mentioned above. Then we have $\sum_{(r,s)\in I}a_{rs} x_r'x_s'>0$.
Furthermore, $\sum_{(r,s)\in I}a_{rs}x_r'x_s'$ is the leading coefficient of 
the polynomial $\hat{q}(\uplambda,x')=q\left(\uplambda^{z_1}x_1',\ldots,\uplambda^{z_n}x_n'\right)\in \mathbb{R}[\uplambda]$. Consider $\hat{q}(\lambda,x')$ for $\lambda\geq 1$. So far, we have shown that the leading coefficient of $\hat{q}(\lambda,x')$ is positive. Therefore $\hat{q}(\lambda,x')\rightarrow \infty$ for $\lambda\rightarrow \infty$. This implies the existence of a $\lambda'\geq 1$ such that $\hat{q}(\lambda,x'')=q\left(\lambda^{z_1}x_1'',\ldots,\lambda^{z_n}x_n''\right)>0$ for all $\lambda\geq 1$, where $x''\in \mathbb{R}^n$ is defined by $x_i''=\lambda '^{z_i}x_i'$ for $i=1,\ldots,n$.
Let $U$ be a 'small' neighborhood of $x''$ such that $\sum_{\alpha\in I}a_\alpha y^\alpha>0$ for all $y\in U$. Interpreting $\hat{q}(\uplambda,x'')$ as a polynomial in $\mathbb{C}[\uplambda]$ and using Proposition \ref{rooters}, we see that no real root of $\hat{q}(\uplambda,y)$ can be greater than $1$, if $U$ is small enough. This implies that $\hat{q}(\lambda,y)$ is positive for all $\lambda\geq 1$ and all $y\in U$. Hence $T_{\overline{U},z}\subseteq S(q)$.
Let us assume that $\mathcal{D}(q)$ is positive definite. 
We need to verify that $\mathcal{T}_0(q)=\mathbb{Z}^n$. Consider the quadratic form $\tilde{q}=q-\mathcal{D}(q)$. Suppose that $\tilde{q}=0$. Since $q(x)\geq \tilde{q}(x)$ holds for all $x\in \mathbb{R}^n$, we see that $S(\tilde{q})$ is a subset of $S(q)$. This implies $S(q)=\mathbb{R}^n$ and $\mathcal{T}_0(q)=\mathbb{Z}^n$. If $\tilde{q}\neq 0$ then we get $\mathcal{T}_0(\tilde{q})=\mathbb{Z}^n$, since $\mathcal{D}(\tilde{q})=0$. 
The inclusion $S(\tilde{q})\subseteq S(q)$ leads straight to $\mathcal{T}_0(q)=\mathbb{Z}^n$.\bew   
\begin{Def}
Let $\varphi=(\varphi_1,\ldots,\varphi_n)\in \mathbb{R}(\uplambda)^n$ be a tuple of rational fractions $\varphi_1,\ldots,\varphi_n\neq 0$. Under a rational tentacle we understand the set 
$$T_{K,\varphi}=\left\{\left(\varphi_1(\lambda)x_1,\ldots,\varphi_n(\lambda)x_n\right):\lambda\geq 1,\varphi_1(\lambda),\ldots,\varphi_n(\lambda)\,\text{is defined},x\in K\right\}$$
where $K\subseteq \mathbb{R}^n$ is an compact set with non-empty interior. Furthermore, we denote by $\mathcal{T}_0(S)$ resp. $\mathcal{T}(S)$ the set of all tentacles resp. rational tentacles that are contained in a semi-algebraic set $S\subseteq \mathbb{R}^n$.  
\end{Def}
\begin{Rem}
We want to establish a link between rational tentacles and the $z$-gradings of $A$. Let $f_1,\ldots,f_s$ be polynomials in $A$ and $\mathcal{T}$ the set of all rational tentacles $T$ such that $T\subseteq S(f_1,\ldots,f_s)$. We say that a tentacle $T\in \mathcal{T}$ is of degree $z\in \mathbb{Z}^n$ if there exists compact set $K\subseteq \mathbb{R}^n$ with non-empty interior and a tuple of rational fractions $\varphi\in \mathbb{R}(\uplambda)^n$ such that $T=T_{K,\varphi}$ and $(\deg(\varphi_1),\ldots,\deg(\varphi_n))=z$, where $\deg$ is defined to be the negative degree valuation $-v_\infty$. Therefore we can assign each $T$ a tuple in $\mathbb{Z}^n$ by $D(\varphi)=(\deg(\varphi_1),\ldots,\deg(\varphi_n))$. However, this assignment is not unique. Thus a rational tentacle may have more degrees than merely just one.
\end{Rem}
\noindent We are now able to generalize Proposition \cite[Proposition 5.1]{z12} and Theorem \ref{kxtt11}:
\begin{Prop}\label{kxpp1}
Let $f_1,\ldots,f_s$ be polynomials in the graded polynomial algebra $A=\bigoplus_{d\in \mathbb{Z}}A_d^{(z)}$, where $z\in \mathbb{Z}^n$. Then the set $$S(L_z(f_1),\ldots,L_z(f_s))\subseteq \mathbb{R}^n$$ 
is Zariski-dense in $\mathbb{R}^n$ if and only if the set $S(f_1,\ldots,f_s)\subseteq \mathbb{R}^n$ contains a rational tentacle $T_{K,\varphi}$ of degree $z$ for some compact set $K\subseteq \mathbb{R}^n$ with non-empty interior.
\end{Prop}
\textbf{Proof}: $\Rightarrow$: The same proof as in \cite[Proposition 5.1]{z12}.\par\smallskip  
$\Leftarrow$: 
For each $k=1,\ldots,s$ let $f_k$ be given by $f_k=\sum_{\alpha}a_{k,\alpha}\mathrm{x}_1^{\alpha_1}\cdots \mathrm{x}_n^{\alpha_n}$, where $a_{k,\alpha}\in \mathbb{R}$.
Suppose there exists a rational tentacle $T:=T_{K,\varphi}$ such that $T\subseteq S(f_1,\ldots,f_s)$ and $D(\varphi)=z$.\par\smallskip  
We are going to show that there exists a point $x\in \mathrm{int}(K)$ and an open neighborhood $U$ of $x$ such that the component of $\hat{f}_i(\uplambda,x)=f_i(\varphi_1x_1,\ldots,\varphi_nx_n)=\sum_{\alpha}c_{i,\alpha}x^\alpha\varphi_1^{\alpha_1}\cdots\varphi_n^{\alpha_n}$ with the highest degree is $\hat{L}_z(f_i)(\uplambda,x)=\sum_{\langle \alpha,z\rangle=\delta_i}c_{i,\alpha}x^\alpha\varphi_1^{\alpha_1}\cdots\varphi_n^{\alpha_n}$ for all $i=1,\ldots,s$ and $x\in U$, where 
$\delta_i=\max\left\{\langle \alpha,z\rangle:c_{i,\alpha}\neq 0\right\}$.
In other words, we have to show that $\deg_\uplambda\left(\hat{L}_z(f_i)(\uplambda,x)\right)=\delta_i$ for all $i=1,\ldots,s$ and all $x\in U$.
There are polynomials $h_{1,\alpha},h_2\in \mathbb{R}[\uplambda]$ such that $\hat{L}_z(f_i)(\uplambda,x)$ can be rewritten as $$\hat{L}_z(f_i)(\uplambda,x)=\frac {\sum_{\langle \alpha,z\rangle=\delta_i}c_{i,\alpha}x^\alpha h_{1,\alpha}} {h_2}.$$
Let $m_\alpha$ be the leading coefficient of $h_{1,\alpha}$ if no other $h_{1,\alpha'}$ ($c_{i,\alpha'}\neq 0$), appearing in the sum above, has a higher degree. Otherwise, we set $m_\alpha=0$. The only situation in which the degree of $\hat{L}_z(f_i)(\uplambda,x)$ in $\uplambda$ is smaller than $\delta_i$, is the one, where $\sum_\alpha c_{i,\alpha}m_\alpha x^\alpha=0$. Since not all $m_\alpha$ can vanish, the sum $\sum_\alpha c_{i,\alpha}m_\alpha \mathrm{x}_1^{\alpha_1}\cdots \mathrm{x}_n^{\alpha_n}$ interpreted as an element of $\mathbb{R}[\mathrm{x}_1,\ldots,\mathrm{x}_n]$ is not the zero polynomial. Since $\mathrm{int}(K)$ is not empty, we can find a point $x\in \mathrm{int}(K)$ such that $\sum_\alpha c_{i,\alpha}m_\alpha x^\alpha\neq 0$. Additionally, we can find a neighborhood $U_i$ of $x$, where $\sum_\alpha c_{i,\alpha}m_\alpha y^\alpha\neq 0$ for all $y\in U_i$. We just proved that for every $y\in U_i$, the degree of $\hat{L}_z(f_i)(\uplambda,y)$ is exactly $\delta_i$.\par\smallskip  
Let us construct the following subset $U$ of $\mathrm{int}(K)$: Start with a point $x_1\in \mathrm{int}(K)$ and an open neighborhood $U_1$ of $x_1$ such that $\deg_\uplambda\left(\hat{L}_z(f_1)(\uplambda,y)\right)=\delta_1$ for all $y\in U_1$. Since $U_1$ is open, we can find another point $x_2\in U_1$ and an open neighborhood $U_2\subseteq U_1$ of $x_2$ such that $\deg_\uplambda\left(\hat{L}_z(f_2)(\uplambda,y)\right)=\delta_2$ for all $y\in U_2$. By repeating this procedure for the remaining polynomials $f_3,\ldots,f_s$, we get the open neighborhoods $U_3,\ldots,U_s$. Set $U=\bigcap_{k=1}^sU_k$. Hence $\deg_\uplambda\left(\hat{L}_z(f_i)(\uplambda,x)\right)=\delta_i$ for all $i=1,\ldots,s$ and all $x\in U$.\par\smallskip  
Fix a point $x'\in U$ with $x'_k\neq 0$ for all $k=1,\ldots,n$ and consider again the rational fraction $\hat{L}_z(f_i)(\uplambda,x')$ and $\hat{f}_i(\uplambda,x')$.
Since $T\subseteq S(f_1,\ldots,f_s)$, there exists a $\lambda_i\geq  1$ such that 

\begin{itemize}
\item $\hat{f}_i(\lambda,x')>0$
\item $\hat{L}_z(f_i)(\lambda,x')>0$
\item $\varphi_1(\lambda),\ldots,\varphi_n(\lambda)\neq 0$
\item $\varphi_1(\lambda),\ldots,\varphi_n(\lambda)$ is defined
\end{itemize}
for all $\lambda\geq \lambda_i$. If we take a small neighborhood $\widetilde{U}_i\subseteq U$ of $x'$, the inequalities $\hat{f}_i(\lambda,y)>0$ and $\hat{L}_z(f_i)(\lambda,y)>0$ will still hold for all $\lambda\geq \lambda_i$ and all $y\in \widetilde{U}_i$: 
To be more precise, we take an open neighborhood $\widetilde{U}_i\subseteq U$ of $x'$, such that for every point $y\in \widetilde{U}_i$ no component of $y$ vanishes. 
According to Corollary \ref{crooters} we can choose $\widetilde{U}_i$ so small that $\hat{f}_i(\lambda,y)$ and $\hat{L}_z(f_i)(\lambda,y)$ have no poles or roots for all $\lambda\geq \lambda_i$ and all $y\in \widetilde{U}_i$. Thus if $\widetilde{U}_i$ is small enough, all real roots of $\hat{f}_i(\lambda,y)$ and $\hat{L}_z(f_i)(\lambda,y)$ will be smaller than $\lambda_i$ for $y\in \widetilde{U}_i$ and therefore $\hat{L}_z(f_i)(\lambda,y)$ resp. $\hat{f}_i(\lambda,y)$ will be positive for all $\lambda\geq \lambda_i$ and all $y\in \widetilde{U}_i$.
Thus if $\widetilde{U}_i$ is a small neighborhood of $x'$, we get $\hat{f}_i(\lambda,y)>0$ and $\hat{L}_z(f_i)(\lambda,y)>0$ for all $\lambda\geq \lambda_i$ and all $y\in \widetilde{U}_i$.      
Set $\widetilde{U}=\bigcap_{i=1}^s\widetilde{U}_i$, $\lambda '=\max_{i\in \{1,\ldots,s\}}\lambda_i$ and consider the map $\psi:\widetilde{U}\rightarrow \mathbb{R}^n,u\mapsto (\varphi_1(\lambda ')u_1,\ldots,\varphi_n(\lambda ')u_n)$. So far, we have shown that $\psi\left(\widetilde{U}\right)\subseteq S(f_1,\ldots,f_s)$ and $\psi\left(\widetilde{U}\right)\subseteq S(L_z(f_1),\ldots,L_z(f_s))$.
Since $\mathrm{int}\left(\psi\left(\widetilde{U}\right)\right)\neq \varnothing$, it is clear that both $S(f_1,\ldots,f_s)$ and especially $S(L_z(f_1),\ldots,L_z(f_s))$ are Zariski-dense in $\mathbb{R}^n$.\bew
\begin{Theo}\label{kxtt1}
Let $A=\bigoplus_{d\in \mathbb{Z}}A_d^{(z)}$ be a $z$-grading and $M$ a finitely generated quadratic module in $A$.\\
If for a set of generators $f_1,\ldots,f_s$ of $M$ the set $S(L_z(f_1),\ldots,L_z(f_s))\subseteq \mathbb{R}^n$ is Zariski dense, then $M$ is totally stable with respect to the $z$-grading. If $M$ is closed under multiplication, then total stability implies the Zariski denseness for any finite set of generators of $M$. 
\end{Theo}
\textbf{Proof}: See \cite[Theorem 4.3]{z12}.\bew
\begin{Theo}\label{kxtt32}
Let $f_1,\ldots,f_s$ be polynomials in the graded polynomial algebra $A=\bigoplus_{d\in \mathbb{Z}}A_d^{(z)}$, where $z\in \mathbb{Z}^n$. If the set $S(f_1,\ldots,f_s)\subseteq \mathbb{R}^n$
contains some rational tentacle $T_{K,\varphi}$, then the quadratic module $M=\mathrm{QM}(f_1,\ldots,f_s)$ is totally stable with respect to the $z$-grading. If $M$ is closed under multiplication, then $S(f_1,\ldots,f_s)$ must contain a tentacle $T_{K,z}$ for $M$ to be totally stable.
\end{Theo}
\textbf{Proof}: Combine Proposition \ref{kxpp1} and Theorem \ref{kxtt1}.\bew  
\begin{Rem}
If someone is interested in stability and the quadratic module $M=\mathrm{QM}(f_1,\ldots,f_s)$ is closed under multiplication, then there is no point using Theorem \ref{kxtt32} over Theorem \ref{kxtt11}. In general, however, this is not true as the next remark will illustrate it. Another advantage of the tentacle is that it is more flexible than an ordinary tentacle. A tentacle may loose its property of being a tentacle even by small manipulations, while it is harder doing so with respect to a rational tentacle.
\end{Rem}
\begin{Rem} \textbf{Stability under isomorphism}:
Let $\chi:\mathbb{R}^n\DistTo \mathbb{R}^n$ be given by a matrix in $\mathrm{GL}_n$. Consider a basic closed semi-algebraic set $S=S(f_1,\ldots,f_s)\subseteq \mathbb{R}^n$. Set $S'=\chi(S)$. Then we have $S'=S\left(f_1\circ \chi^{-1},\ldots,f_s\circ \chi^{-1}\right)$. Hence $S'$ is again a semi-algebraic set. In the following we write $\chi_i=\sum_{j}a_{ij}\mathrm{x}_j$, where $a_{ij}\in \mathbb{R}$.
Reprise that $\mathcal{T}_1:=\mathcal{T}(S)$ resp. $\mathcal{T}_2:=\mathcal{T}\left(S'\right)$ is the set of all rational tentacles $T$ such that $T\subseteq S$ resp. $T\subseteq S'$. Let $D_1$ resp. $D_2$ denote the set of all degrees of all tentacles in $\mathcal{T}_1$ resp. $\mathcal{T}_2$.
We are interested in the relationship between $D_1$ and $D_2$.
Take a rational tentacle $T:=T_{K,\varphi}\in \mathcal{T}_1
$ of degree $z\in \mathbb{Z}^n$. For the sake of simplicity let us assume that $\varphi$ is defined on the set $[1,\infty)$. Let the $i$-th component of $\chi\left(\varphi_1\mathrm{x}_1,\ldots,\varphi_n\mathrm{x}_n\right)$ be given by
$$\chi_i\left(\varphi_1\mathrm{x}_1,\ldots,\varphi_n\mathrm{x}_n\right)=\sum_{j}a_{ij}\varphi_j\mathrm{x}_j=\sum_{j}a_{ij}\mathrm{x}_i^{-1}\varphi_j\mathrm{x}_j\mathrm{x}_i.$$
Choose a point $x\in \mathrm{int}(K)$ that satisfies the following two conditions: 
\begin{itemize}
\item For each $i=1,\ldots,n$ the $i$-th component of $x$ and $\chi(x)$ does not vanish.
\item For each $i=1,\ldots,n$ the degree of $\chi_i\left(\varphi_1x_1,\ldots,\varphi_nx_n\right)\in \mathbb{R}[\uplambda]$
is equal to $\tilde{z}_i:=\max\left\{z_j:a_{ij}\neq 0,j=1,\ldots,n\right\}$. 
\end{itemize}
Set $\tilde{\varphi}_i=\sum_{j}a_{ij}x_i^{-1}\varphi_jx_j$. There is a small neighborhood $U\subseteq \mathrm{int}(K)$ of $x$ such that $\left(\tilde{\varphi}_1(\lambda)y_1,\ldots,\tilde{\varphi}_n(\lambda)y_n\right)$ lies in $S'$ for all $y\in \overline{U}$ and all $\lambda\geq 1$.
Let us prove this assertion. Set $\tilde{a}_{ij}=a_{ij}x_i^{-1}y_i$. Thus every $y\in U$ defines a small perturbation of the coefficients $a_{ij}$. Hence we can write $\tilde{a}_{ij}=a_{ij}+\varepsilon_{ij}$, where $\varepsilon_{ij}\in \mathbb{R}$. Note that $\varepsilon_{ij}=0$ if $a_{ij}=0$. Hence for small $\varepsilon_{ij}$ we get $\sum_j\tilde{a}_{ij}\varphi_j(\lambda)x_j=\sum_ja_{ij}\varphi_j(\lambda)x_j+\sum_j\varepsilon_{ij}\varphi_j(\lambda)x_j\geq 0$ for all $i=1,\ldots,n$ and all $\lambda\geq 1$. Therefore
$$\widetilde{T}=\left\{\left(\tilde{\varphi}_1(\lambda)y_1,\ldots,\tilde{\varphi}_n(\lambda)y_n\right):y\in \overline{U},\lambda\geq 1\right\}\subseteq S'.$$
By the continuity of $\chi$, an open neighborhood $U'\subseteq U$ of $x$ can be found such that 
$$T'=\left\{\left(\tilde{\varphi}_1(\lambda)\frac {x_1} {\chi_1(x)}\chi_1(y),\ldots,\tilde{\varphi}_n(\lambda)\frac {x_n} {\chi_n(x)}\chi_n(y)\right):y\in \overline{U'},\lambda\geq 1\right\}\subseteq \widetilde{T}.$$
Now, the identity $$T'=\left\{\left(\tilde{\varphi}_1(\lambda)\frac {x_1} {\chi_1(x)}y_1',\ldots,\tilde{\varphi}_n(\lambda)\frac {x_n} {\chi_n(x)}y_n'\right):y'\in \chi\left(\overline{U'}\right),\lambda\geq 1\right\}$$
proves that $T'$ is in $\mathcal{T}_2$.\par\smallskip  
The degree of $T'$ is given by $\left(\deg\left(\tilde{\varphi}_1\right),\ldots,\deg\left(\tilde{\varphi}_n\right)\right)$, which is nothing more than\\ 
$\tilde{z}:=\left(\tilde{z}_1,\ldots,\tilde{z}_n\right)$.
This gives us a map $\xi_1:D_1\rightarrow D_2,z\mapsto \tilde{z}$. 
On the other hand, we can start with a tentacle $T'\in \mathcal{T}_2$ and repeat the same argumentation done so far by replacing $\chi$ with $\chi^{-1}$. This gives us a map $\xi_2:D_2\rightarrow D_1$.\par\smallskip  
If $\chi=\mathrm{id}_{\mathbb{R}^n}$ then it is clear that $\xi_1$ and $\xi_2$ are the identity maps. Suppose $\chi\in \mathrm{GL}_n$ is not the identity map. Even under this circumstances neither $\xi_1$ nor $\xi_2$ need to be linear or inverse to each other.\par\smallskip        
By using Theorem \ref{kxtt32} we see that if $M=\mathrm{QM}(f_1,\ldots,f_s)$ is stable with respect to a $z$-grading, then $M'=\mathrm{QM}\left(f_1\circ \chi^{-1},\ldots,f_s\circ \chi^{-1}\right)$ is stable with respect to a $\xi_1(z)$-grading. 
On the other hand, if $M'$ is stable with respect to a $z'$-grading, then $M$ is stable with respect to a $\xi_2(z')$-grading. 
Note that this result is impossible by just using the ordinary tentacle defined in \ref{kxd1x} and Theorem \ref{kxtt11}.   
\end{Rem}
\section{Tentacles and the S4-conjecture} 
Let $n$ be a natural number.
For a subset $I\subseteq \{1,\ldots,n\}$ we define an involution $\pi_I:\mathbb{R}^n\rightarrow \mathbb{R}^n$ by $\pi_I(x)=x'$ where $x_j'=x_j$ for $j\notin I$ and $x_j'=-x_j$ for $j\in I$. This kind of maps form a group $G$. Furthermore every $\pi\in G$ maps a rational tentacle $T_{K,\varphi}$ to another rational tentacle $\pi_I(T_{K,\varphi})=T_{\pi_I(K),\varphi}$. Let $f$ be a polynomial in $A$. In the following we denote by $L_{\mathrm{lex}}(f)$ the leading term of $f$ with respect to the lexicographical ordering.   
Finally, set $\mathbb{N}^n_1:=\left\{z\in \mathbb{N}^n:\forall i\in \{1,\ldots,n\}:z_1\geq z_i\right\}$. Then we can state:
\begin{Theo}\label{tult1}
Let $S_1=S(q)$ and $S_2=S(p)$ be two semi-algebraic sets in $\mathbb{R}^n$. Suppose the following conditions are satisfied:
\begin{enumerate}[label=(\alph*)]
\item We have $L_{\mathrm{lex}}(p)\notin L_{\mathrm{lex}}(q)A$.
\item For every $z\in \mathbb{N}^n_1$ there exists a rational tentacle $T\in \mathcal{T}(S_1)$ 
of degree $z$ and an element $\pi\in G$ such that $\pi(T)\notin \mathcal{T}(S_2)$. Furthermore, all unbounded $T'\in \mathcal{T}(S_1)$ with $\pi(T')\subseteq \pi(T)$ satisfy $\pi(T')\notin \mathcal{T}(S_2)$.
\end{enumerate} 
Then there is no non-negative polynomial $t\in A$ such that $p(y)-t(y)q(y)\geq 0$ for all $y\in \mathbb{R}^n$.
\end{Theo}
\textbf{Proof}: Without loss of generality we can assume that $S_1\subseteqq S_2$. Suppose that there exists a non-negative polynomial $t\in A$ such that $p(y)-t(y)q(y)\geq 0$ for all $y\in \mathbb{R}^n$. The prove is divided in several steps:\par\smallskip  
(i): Let $L_{\mathrm{lex}}(q)=a_\alpha \mathrm{x}_1^{\alpha_1}\cdots \mathrm{x}_n^{\alpha_n}$, $L_{\mathrm{lex}}(p)=b_\beta \mathrm{x}_1^{\beta_1}\cdots \mathrm{x}_n^{\beta_n}$ and $L_{\mathrm{lex}}(t)=d_\gamma \mathrm{x}_1^{\gamma_1}\cdots \mathrm{x}_n^{\gamma_n}$ denote the leading terms of $q$, $p$, and $t$ with respect to the lexicographical ordering. Furthermore, we define $I(q)=\left\{\delta\in \mathbb{N}_0^n:a_\delta\neq 0,\delta\neq \alpha\right\}$, $I(p)=\left\{\delta\in \mathbb{N}_0^n:b_\delta\neq 0,\delta\neq \beta\right\}$ and $I(t)=\left\{\delta\in \mathbb{N}_0^n:d_\delta\neq 0,\delta\neq \gamma\right\}$. We are going to show that there is a tuple $z\in \mathbb{N}^n_1$ such that
\begin{itemize} 
\item $\langle \alpha,z\rangle>\langle \delta,z\rangle$ for all $\delta\in I(q)$.
\item $\langle \beta, z\rangle>\langle \delta,z\rangle$ for all $\delta\in I(p)$.
\item $\langle \gamma,z\rangle>\langle \delta,z\rangle$ for all $\delta\in I(t)$.  
\end{itemize}
Let us start with $z'=(1,\ldots,1)$. Consider the set $N(z')=\left\{\delta\in I(t):\langle \gamma,z'\rangle\leq \langle\delta,z'\rangle\right\}\cup \left\{\delta\in I(p):\langle \beta,z'\rangle\leq \langle\delta,z'\rangle\right\}\cup \left\{\delta\in I(q):\langle \alpha,z'\rangle\leq \langle\delta,z'\rangle\right\}$ and the numbers $$r_1(z')=\begin{cases}\max\{j\in \mathbb{N}_{\leq n}:\exists \delta\in N(z')\cap I(t)\forall i\leq j:\delta_i=\gamma_i\}\,\text{if}\,N(z')\cap I(t)\neq \varnothing\\
0,\,\text{otherwise}
\end{cases},$$
$$r_2(z')=\begin{cases}\max\{j\in \mathbb{N}_{\leq n}:\exists \delta\in N(z')\cap I(p)\forall i\leq j:\delta_i=\beta_i\}\,\text{if}\,N(z')\cap I(p)\neq \varnothing\\
0,\,\text{otherwise}
\end{cases},$$ 
$$r_3(z')=\begin{cases}\max\{j\in \mathbb{N}_{\leq n}:\exists \delta\in N(z')\cap I(q)\forall i\leq j:\delta_i=\alpha_i\}\,\text{if}\, N(z')\cap I(q)\neq \varnothing\\
0,\,\text{otherwise}
\end{cases},$$ and 
$r(z')=\max\{r_1(z'),r_2(z'),r_3(z')\}$. Suppose that $r(z')=r_1(z')$. We see that $n>r_1$,
 since $\gamma\succ_{\mathrm{lex}} \delta$ for all $\delta\in N(z')\cap I(t)$ with respect to the lexicographical ordering. Now take $\delta\in N(z')\cap I(t)$ whose components $1,\ldots,r_1(z'):=r_1$ are identical to those of $\gamma$. Since $\gamma\succ_{\mathrm{lex}} \delta$, the inequality $\gamma_{r_1+1}>\delta_{r_1+1}$ must hold. Now we can enlarge the $r_1+1$-th component of $z'$ in such a way that $\langle \gamma,z'\rangle>\langle \delta,z'\rangle$. In fact, we can achieve $\langle \gamma,z'\rangle>\langle \delta,z'\rangle$ for all $\delta\in N(z')\cap I(t)$ whose first $r_1$ components are identical to those of $\gamma$. If $r_2(z')=r_1(z')$ or $r_3(z')=r_1(z')$ we, if necessary, enlarge the $r_1+1$-th component of $z'$ further such that both inequalities $\langle \beta,z'\rangle>\langle \delta_1,z'\rangle$, $\langle \alpha,z'\rangle>\langle \delta_2,z'\rangle$ hold for all $\delta_1\in N(z')\cap I(p)$, $\delta_2\in N(z')\cap I(q)$ whose $1,\ldots,r_1$ components are identically to whose of $\beta$ resp. $\alpha$. If $r(z')=r_2(z')$ or $r(z')=r_3(z')$ just use the same argumentation again. That is, replace $r_1(z')$ by $r_2(z')$ or $r_3(z')$ and simply repeat the argumentation done in this matter.\par\smallskip     
Let $z''$ denote $z'$ with the enlarged $r_1+1$-th component and consider $N(z'')$ resp. $r(z'')$. Then it is clear that $r(z'')\leq r(z')-1$. 
Start over again with the new data $N(z'')$ and $r(z'')$ and notice that after each finished repetition the value $r(z'')$ will decrease at least by one.\par\smallskip     
Thus after $k\leq r(z')$ repetitions we finally get a tuple $z:=z^{(k)}$ that will satisfy $N(z)=\varnothing$ resp. $r(z)=0$, which is the same thing as saying that $z$ will satisfy all three inequalities mentioned above. It is obvious that we can choose $z$ in such a way that the first component is the largest. To be more precise, if the first component of $z$ is not the largest, then we can enlarge it without violating the three inequalities.\par\smallskip    
(ii): According to condition (b) in the theorem, there exists a rational tentacle $T\in \mathcal{T}(S_1)$ of degree $z$ and an involution $\pi\in G$ such that all statements in (b) are satisfied. Write $T_{K,\varphi}$ for $T$, where as usual $K\subseteq \mathbb{R}^n$ is a compact set with non-empty interior and $\varphi\in \mathbb{R}(\uplambda)^n$. For any polynomial $f\in A$, we define $\hat{f}(\uplambda,x')$ to be the rational fraction $\hat{f}(\uplambda,x')=f\left(\varphi_1(\uplambda)x_1',\ldots,\varphi_n(\uplambda)x_n'\right)\in \mathbb{R}(\uplambda)$, where $x'\in \mathrm{int}(K)$.
We know from condition (b) and (a) that there is a point $x'\in \mathrm{int}(K)$ with $x_i'\neq 0$ for all $i=1,\ldots,n$ such that
$$\deg_\uplambda\left(a_\alpha x^\alpha \varphi_1^{\alpha_1}\cdots \varphi_n^{\alpha_n}-d_\gamma b_\beta x^{\gamma+\beta}\varphi_1^{\gamma_1+\beta_1}\cdots \varphi_n^{\gamma_n+\beta_n}\right)=\max\left\{\langle \alpha,z\rangle,\langle \gamma+\beta,z\rangle\right\}$$
resp.
$$\deg_\uplambda\left(\hat{L}_{\mathrm{lex}}(p)(\uplambda,x')-\hat{L}_{\mathrm{lex}}(t)(\uplambda,x')\hat{L}_{\mathrm{lex}}(q)(\uplambda,x')\right)=\max\left\{\langle \alpha,z\rangle,\langle \gamma+\beta,z\rangle\right\}$$
holds.
Set $x=\pi(x')$ and $\varrho(\hat{q},x)=\left\{\lambda\in \mathbb{R}_{\geq 1}:\hat{q}(\lambda,x)\geq 0,\hat{q}(\lambda,x)\,\text{is defined}\right\}$.\\
\begin{figure}[h]
\begin{tikzpicture}

\path[name path=border1] (0,0) to[out=-10,in=150] (6,-2);
\path[name path=border2] (12,1) to[out=150,in=-10] (5.5,3.2);
\path[name path=redline] (0,-0.4) -- (12,1.5);
\path[name path=redlinex] (0,-0.4) -- (10,3);

\path[name intersections={of=border1 and redline,by={a}}];
\path[name intersections={of=border2 and redline,by={b}}]; 
\shade[left color=gray!10,right color=gray!80] 
  (0,0) to[out=-10,in=150] (6,-2) -- (12,1)to[out=150,in=-10] (5.5,3.7) -- cycle;
\draw[dashed,line width=2.8pt,shorten >= 3pt,shorten <= 3pt, opacity=0.3]
(0,0) to[bend left] (-1,4) to[bend right] (-0.3,3.5);
\draw[dashed,line width=2.8pt,shorten >= 3pt,shorten <= 3pt, opacity=0.3]
(-1,4) to[bend right] (6,7) to [bend left] (5.4,3.7);
\draw[dashed,line width=2.8pt,shorten >= 3pt,shorten <= 3pt, opacity=0.3]
(6,7) to[bend left] (9,7.5);
\shade [ball color=white] (5,2.5) circle [radius=1cm];
\draw[blue ,line width=1.8pt,shorten >= 3pt,shorten <= 3pt]
(5.4,1.5) to[bend right] (8.1,4.6);
\draw[blue ,line width=1.8pt,shorten >= 3pt,shorten <= 3pt]
(5,3.5) to[bend right] (7.7,5.1);
\draw[dashed, blue ,line width=1.8pt,shorten >= 3pt,shorten <= 3pt]
(5,3.5) to[bend right] (5.4,1.5);
\draw[blue ,line width=1.8pt,shorten >= 3pt,shorten <= 3pt]
(7.55,4.15) to[bend left] (5.7,3);
\draw [fill opacity=0.3,fill=blue!80!blue]
(5,3.5) to[bend right] (5.4,1.5) to[bend right] (8.1,4.5) to[bend left]  (7.5,4) to [bend left] (7.7,5) to[bend left] (5,3.5);
\draw [fill opacity=0.3,fill=blue!80!blue]
(8.1,4.5) to[bend left]  (7.5,4) to [bend left] (7.7,5) to[bend left] (8.1,4.5);
\shade[left color=red!20, right color=red!60]
(-0.3,3.5) to[out=-10,in=225] (10,7.5) to[bend right] (9,7.5) to[out=225,in=-10] (1,4);
\shade[left color=red!0, right color=red!50]
(b) to[out=172,in=-10]  (a) to[out=10,in=150] (b);
\draw[dashed,line width=2.8pt,shorten >= 3pt,shorten <= 3pt] 
(-0.3,3.5) to[out=-10,in=225]  
  coordinate[pos=0.27] (aux1) 
  coordinate[pos=0.52] (aux2) 
  coordinate[pos=0.75] (aux3) (10,7.5);
\draw[dashed,red,line width=2.8pt,shorten >= 3pt,shorten <= 3pt] 
(1,4) to[out=-10,in=225]  
  coordinate[pos=0.1] (auxx1) 
  coordinate[pos=0.6] (auxx2) 
  coordinate[pos=0.9] (auxx3) (9,7.5);
\draw[dashed,line width=2.8pt,shorten >= 3pt,shorten <= 3pt] 
  (b) to[out=172,in=-10] 
	coordinate[pos=0.8] (bux1) 
  coordinate[pos=0.5] (bux2) 
  coordinate[pos=0.2] (bux3) (a);
\draw[dashed,red,line width=2.8pt,shorten >= 3pt,shorten <= 3pt] 
  (b) to[out=150,in=10] 
	coordinate[pos=0.9] (buxx1) 
  coordinate[pos=0.4] (buxx2) 
  coordinate[pos=0.15] (buxx3) (a);
\foreach \coor in {1,2,3}
  \draw[dashed,red,line width=1.8pt,shorten >= 3pt,shorten <= 3pt] (auxx\coor)to[bend right] (buxx\coor);
\foreach \coor in {1,2,3}
  \draw[dashed,line width=1.8pt,shorten >= 3pt,shorten <= 3pt] (aux\coor)to[bend left] (bux\coor);
\shade[left color=black!0, right color=black!80]
 (a) to[out=-10,in=150] (6,-2) --  (12,1)to[out=150,in=-10] (b) to[out=172,in=-10] (a); 
\draw [fill opacity=0.2,fill=red!80!red]
(b) to[out=172,in=-10] (a) to[bend left] (-0.3,3.5)  to[out=-10,in=225]  (10,7.5) to[bend left] (b);
\draw [fill opacity=0.2,fill=red!80!red]
(b) to[out=150,in=10] (a)  to[bend left] (1,4) to[out=-10,in=225] (9,7.5) to[bend left] (b);
\shade [ball color=white] (7,0.5) circle [radius=1cm];
\draw[blue ,line width=1.8pt,shorten >= 3pt,shorten <= 3pt]
(7.4,-0.4) to[bend left] (12.1,1);
\draw[blue ,line width=1.8pt,shorten >= 3pt,shorten <= 3pt]
(6.7,1.35) to[bend left] (12.1,2);
\draw[dashed, blue ,line width=1.8pt,shorten >= 3pt,shorten <= 3pt]
(7.4,-0.4) to[bend left] (6.7,1.35);
\draw[blue ,line width=1.8pt,shorten >= 3pt,shorten <= 3pt]
(11.9,1.5) to[bend right] (7.9,0.5);
\draw [fill opacity=0.3,fill=blue!80!blue]
(7.4,-0.4) to[bend left] (6.7,1.35) to[bend left] (12,2) to[bend right] (11.8,1.5) to[bend right] (12,1) to[bend right] (7.4,-0.4); 
\draw [fill opacity=0.3,fill=blue!80!blue]
(12,2) to[bend right] (11.8,1.5)  to[bend right] (12,1) to[bend right] (12,2);
\draw[dashed,line width=2.8pt,shorten >= 3pt,shorten <= 3pt, opacity=0.3] 
(-0.3,3.5) to[bend right] (7,4);
\draw[dashed,line width=2.8pt,shorten >= 3pt,shorten <= 3pt, opacity=0.3] 
(6,-2) to[bend left] (7,4);
\draw[dashed,line width=2.8pt,shorten >= 3pt,shorten <= 3pt, opacity=0.3]
(7,4) to[bend right] (11.5,6.5);
\draw[dashed,line width=2.8pt,shorten >= 3pt,shorten <= 3pt, opacity=0.3]
(12,1) to[bend right] (11.5,6.5);
\draw[dashed,line width=2.8pt,shorten >= 3pt,shorten <= 3pt, opacity=0.3]
(11.5,6.5) to[bend right] (10,7.5);
\node[rotate=20] at (4.5,0.4) {\color{red}{$S_2\backslash S_1$}};
\node[rotate=30] at (7.7,0.8) {\color{green}{$\pi^{-1}\left(T_{\overline{U},\hat{\varphi}}\right)$}};
\node[rotate=30] at (5.5,2.5) {\color{green}{$T_{\overline{U},\hat{\varphi}}$}};
\end{tikzpicture}
\caption{This is how $T_{\overline{U},\hat{\varphi}}$ and $\pi^{-1}\left(T_{\overline{U},\hat{\varphi}}\right)$ must behave: The \textcolor{red}{red area} is the complement between $S_2$ and $S_1$. The rational tentacle \textcolor{green}{$T_{\overline{U},\hat{\varphi}}$} is contained in $\mathbb{R}^n\backslash S_2$ and \textcolor{green}{$\pi^{-1}\left(T_{\overline{U},\hat{\varphi}}\right)$} is contained in $S_1$.}  
\label{fig:fig1}
\end{figure}
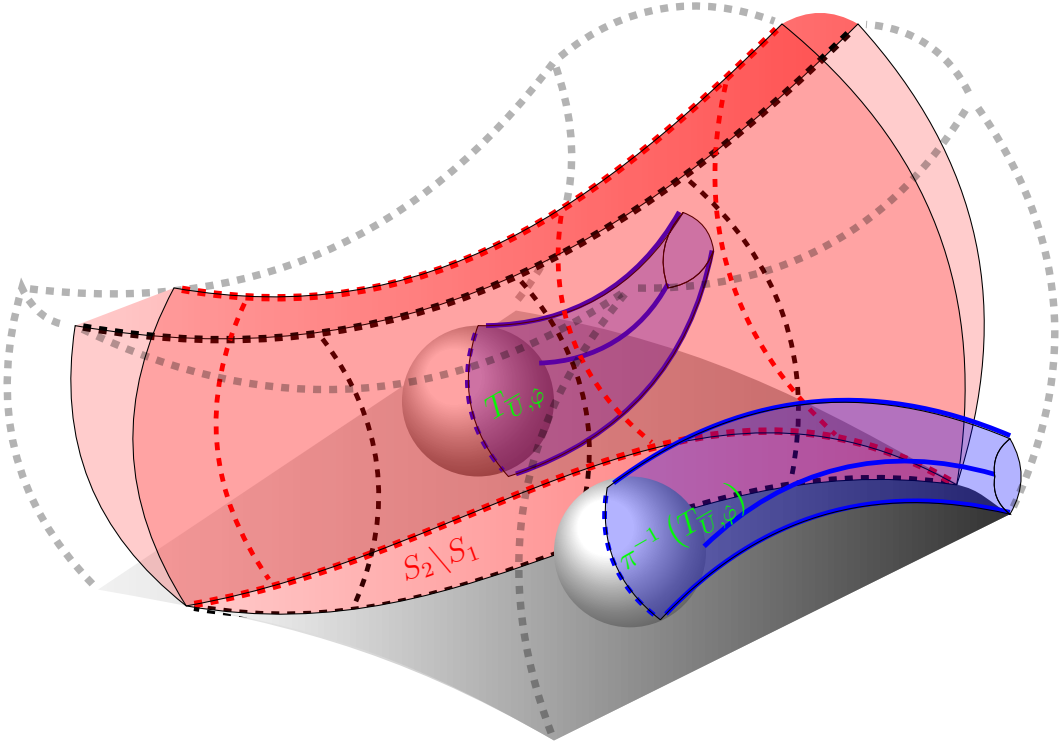\\
We are going to exclude that $\varrho(\hat{p},x)$ is unbounded. Suppose the opposite would be the case. 
Then there exists a $\lambda '\geq 1$ such that the rational fractions $\varphi_1,\ldots,\varphi_n$ have no poles and no roots, and that $\hat{p}(\lambda,x)$ is positive for all $\lambda\geq \lambda '$.\par\smallskip   
By taking a small neighborhood $U\subseteq \pi(K)$ of $x$ we can make sure that $\hat{p}(\lambda,y)$ will be positive for all $\lambda\geq \lambda'$ and all $y\in U$. Next, we define the rational fractions $\hat{\varphi}_i=\varphi_i(\uplambda+\uplambda'-1)$ for $i=1,\ldots,n$. Set $\hat{\varphi}=(\hat{\varphi}_1,\ldots,\hat{\varphi}_n)$. 
Then the rational tentacle $T_{\overline{U},\hat{\varphi}}$ lies in $\mathcal{T}(S_2)$.
According to our construction, $\pi^{-1}\left(T_{\overline{U},\hat{\varphi}}\right)$ is a subset of $T_{K,\varphi}$: This follows from $q(\hat{\varphi}_1(\lambda)y_1',\ldots,\hat{\varphi}_n(\lambda)y_n')=\hat{q}(\lambda+\lambda '-1,y')\geq 0$ for all $\lambda\geq 1$ and $y'\in \pi^{-1}\left(\overline{U}\right)\subseteq K$. 
Furthermore, this implies that $T_{\overline{U},\hat{\varphi}}$ is a subset of $\pi(T_{K,\varphi})$. Hence it is easy to see that while $T_{\overline{U},\hat{\varphi}}$ lies in $\mathcal{T}(S_2)$, the other rational tentacle $\pi^{-1}\left(T_{\overline{U},\hat{\varphi}}\right)$ lies in $\mathcal{T}(S_1)$.
Since the degree of $T_{K,\varphi}$ is $z$, we get $\lim_{\lambda\rightarrow \infty}|\varphi_1(\lambda)|=\infty$. Thus $\lim_{\lambda\rightarrow \infty}|\hat{\varphi}_1(\lambda)|=\infty$ and therefore $T_{\overline{U},\hat{\varphi}}$ is not bounded. But that contradicts (b). In fact, the rational tentacles $T_{\overline{U},\hat{\varphi}}$ and $\pi^{-1}\left(T_{\overline{U},\hat{\varphi}}\right)$ must behave like depicted in Figure \ref{fig:fig1} for a suitable $\lambda'\geq 1$ and neighborhood $U\subseteq \pi(K)$ of $x$. In other words, we just saw that $S_2\backslash S_1$ is just too small to contain a rational tentacle that would allow $\varrho(\hat{p},x)$ to be unbounded.\par\smallskip  
(iii): In (ii) we showed that $\varrho(\hat{p},x)$ is bounded resp. that $\varrho(-\hat{p},x)$ is unbounded. Now, we want the same thing for $\varrho(\hat{q},x)$ resp. $\varrho(-\hat{q},x)$. Without loss of generality, we can assume that $\pi(T)$ is not contained in $\mathcal{T}(S_1)$. Otherwise, we get $\hat{p}(\lambda,x)-\hat{t}(\lambda,x)\hat{q}(\lambda,x)<0$ for $\lambda>0$ big enough and therefore we are done. Thus $\varrho(-\hat{q},x)$ must be unbounded, since an infinite part of $\pi(T)$ must lie in $\mathbb{R}^n\backslash S_1$.\par\smallskip  
(iv): We know that $\varrho(-\hat{q},x)$ and $\varrho(-\hat{p},x)$ are unbounded. 
Thus there exists a real number $\lambda_0\geq 1$ such that $\hat{p}(\lambda,x)$ and $\hat{q}(\lambda,x)$ are defined for all $\lambda\geq \lambda_0$.
There is a positive real number $\hat{\lambda}\geq \lambda_0$ such that $\hat{L}_{\mathrm{lex}}(p)(\lambda,x)<0$, $\hat{L}_{\mathrm{lex}}(q)(\lambda,x)<0$ and $\hat{p}(\lambda,x)-\hat{t}(\lambda,x)\hat{q}(\lambda,x)>0$ for all $\lambda\geq \hat{\lambda}$. 
This implies $|\hat{L}_{\mathrm{lex}}(p)(\lambda,x)|<|\hat{L}_{\mathrm{lex}}(t)(\lambda,x)\hat{L}_{\mathrm{lex}}(q)(\lambda,x)|$ for all $\lambda\geq \hat{\lambda}$, if $\hat{\lambda}$ is large enough.\par\smallskip   
The same inequality $|\hat{L}_{\mathrm{lex}}(p)(\lambda,x')|<|\hat{L}_{\mathrm{lex}}(t)(\lambda,x')\hat{L}_{\mathrm{lex}}(q)(\lambda,x')|$ holds for $x'=\pi^{-1}(x)$ and all $\lambda\geq \hat{\lambda}$.
But here we have $\hat{L}_{\mathrm{lex}}(p)(\lambda,x')>0$, $\hat{L}_{\mathrm{lex}}(q)(\lambda,x')>0$ for all $\lambda\geq \hat{\lambda}$. Thus $\hat{L}_{\mathrm{lex}}(p)(\lambda,x')-\hat{L}_{\mathrm{lex}}(t)(\lambda,x')\hat{L}_{\mathrm{lex}}(q)(\lambda,x')<0$ for all $\lambda\geq \hat{\lambda}$. If we choose an appropriately large $\lambda\geq \lambda_0$, we will get $\hat{p}(\lambda,x)-\hat{t}(\lambda,x)\hat{q}(\lambda,x)<0$. However, this contradicts our assumption that $p(y)-t(y)q(y)\geq 0 $ for all $y\in \mathbb{R}^n$.\bew  
\begin{Prop}\label{scor}
Let $S_1=S(q)$ and $S_2=S(p)$ be two semi-algebraic sets in $\mathbb{R}^n$. Suppose the following conditions are satisfied:
\begin{enumerate}[label=(\alph*)]
\item We have $L_{\mathrm{lex}}(p)\notin L_{\mathrm{lex}}(q)A$.
\item The quadratic modules $\mathrm{QM}(q),\mathrm{QM}(-p)$ are totally stable with respect to any $z$-grading in $\mathbb{N}_1^n$ and neither $q=0$ nor $p=0$.
\end{enumerate} 
Then there is no non-negative polynomial $t\in A$ such that $p(y)-t(y)q(y)\geq 0$ for all $y\in \mathbb{R}^n$.
\end{Prop}
\textbf{Proof}: Let us start where part (i) in the proof of Theorem \ref{tult1} ended.
Unfortunately we need some new arguments, since condition (b) of Proposition \ref{scor} differs from that of Theorem \ref{tult1}.     
This is where the new part (ii') comes in. It serves as a link between part (i) and (ii) of Theorem \ref{tult1}, with the purpose that we can use the arguments already developed in the preceding theorem. For the sake of simplicity let us assume that $S_1\subseteqq S_2$.\par\smallskip  
(ii'): Let $z\in \mathbb{N}_1^n$ be same tuple we used in part (ii). According to Theorem \ref{kxtt11}, we can find two compact sets $K,K'\subseteq \mathbb{R}^n$ with non-empty interior such that $T_{K,z}\in \mathcal{T}_0(S(q))$ and $T_{K',z}\in \mathcal{T}_0(S(-p))$. Furthermore,  we can find two points $x\in \mathrm{int}(K)$ and $x'\in \mathrm{int}(K)$ with non-vanishing components. Note that $\hat{L}_{\mathrm{lex}}(q)(\uplambda,x)=a_\alpha x^\alpha \uplambda^{\langle \alpha,z\rangle}$, $\hat{L}_{\mathrm{lex}}(q)(\uplambda,x')=a_\alpha x'^\alpha \uplambda^{\langle \alpha,z\rangle}$, $\hat{L}_{\mathrm{lex}}(p)(\uplambda,x)=b_\beta x^\beta \uplambda^{\langle \beta,z\rangle}$ and $\hat{L}_{\mathrm{lex}}(p)(\uplambda,x')=b_\beta x'^\beta \uplambda^{\langle \beta,z\rangle}$.
It is obvious that there are two positive real numbers $\lambda_1$ and $\lambda_2$ such that 
$$\hat{q}(\lambda,x),\hat{p}(\lambda,x),\hat{L}_{\mathrm{lex}}(q)(\lambda,x),\hat{L}_{\mathrm{lex}}(p)(\lambda,x)>0$$ holds for all $\lambda\geq \lambda_1$ resp. 
$$\hat{q}(\lambda,x'),\hat{p}(\lambda,x'),\hat{L}_{\mathrm{lex}}(q)(\lambda,x')
,\hat{L}_{\mathrm{lex}}(p)(\lambda,x')<0$$ 
holds for all $\lambda\geq \lambda_2$.\\ 
Set $M(\alpha)=\left\{i:\mathrm{sgn}\left(x_i^{\alpha_i}\right)\neq \mathrm{sgn}\left(x_i'^{\alpha_i}\right)\right\}$ and $M(\beta)=\left\{i:\mathrm{sgn}\left(x_i^{\beta_i}\right)\neq \mathrm{sgn}\left(x_i'^{\beta_i}\right)\right\}$. Then both sets are not empty, because otherwise we would get $$\mathrm{sgn}\left(\hat{L}_{\mathrm{lex}}(q)(\lambda,x)\right)=\mathrm{sgn}\left(\hat{L}_{\mathrm{lex}}(q)(\lambda,x')\right)$$ or
$$\mathrm{sgn}\left(\hat{L}_{\mathrm{lex}}(p)(\lambda,x)\right)=\mathrm{sgn}\left(\hat{L}_{\mathrm{lex}}(p)(\lambda,x')\right)$$ for all $\lambda\geq 1$, which would result in a contradiction. The intersection $M(\alpha)\cap M(\beta)$ is not empty, since $S(q)\subseteq S(p)$:\par\smallskip 
Suppose that the intersection would be empty. Take an element $i\in M(\beta)$. Set $y_1=(1,\ldots,1,x_i,1,\ldots,1)$ and $y_2=(1,\ldots,1,x_i',1,\ldots,1)$. 
Then $\hat{L}_{\mathrm{lex}}(p)(\lambda,y_1)$ and $\hat{L}_{\mathrm{lex}}(p)(\lambda,y_2)$ have different signs for all $\lambda\geq 1$, while $\hat{L}_{\mathrm{lex}}(q)(\lambda,y_1)$ and $\hat{L}_{\mathrm{lex}}(q)(\lambda,y_2)$ have the same sign for all $\lambda\geq 1$. Thus $\hat{q}(\lambda,y_1)$ and $\hat{q}(\lambda,y_2)$ are both negative or positive, while $\hat{p}(\lambda,y_1)$ and $\hat{p}(\lambda,y_2)$ have different signs for all $\lambda$ large enough. This can only work if $\hat{q}(\lambda,y_1)$ and $\hat{q}(\lambda,y_2)$ are negative for all sufficiently large $\lambda$. Thus $a_\alpha<0$. On the other side, we get $b_\beta>0$ by repeating the same arguments with $j\in M(\alpha)$. It is not hard to see that this cannot work. Set $$\tilde{y}=\begin{cases}(1,\ldots,1,x_j',1,\ldots,1,x_i',1,\ldots,1)\,\text{if}\, x_j'<0\, \text{and}\, x_i'<0\\
(1,\ldots,1,x_j,1,\ldots,1,x_i',1,\ldots,1)\,\text{if}\, x_j<0\, \text{and}\, x_i'<0\\
(1,\ldots,1,x_j,1,\ldots,1,x_i,1,\ldots,1)\,\text{if}\, x_j<0\, \text{and}\, x_i<0\\
(1,\ldots,1,x_j',1,\ldots,1,x_i,1,\ldots,1)\,\text{if}\, x_j'<0\, \text{and}\, x_i<0
\end{cases}$$ 
In fact, $\hat{q}(\lambda,\tilde{y})$ is positive for all $\lambda>1$ large enough, while $\hat{p}(\lambda,\tilde{y})$ is negative, since $a_\alpha \tilde{y}^\alpha>0$ and $b_\beta\tilde{y}^\beta<0$. But this contradicts $S(q)\subseteq S(p)$.
Take a natural number $k$ out of the set $M(\alpha)\cap M(\beta)$ and let $\pi_{k}$ denote the map $\mathbb{R}^n\rightarrow \mathbb{R}^n,(x_1,\ldots,x_k,\ldots,x_n)\mapsto (x_1,\ldots,-x_k,\ldots,x_n)$.
Then we can find a two positive real numbers $c_1$ and $c_2$ such that the following two equations 
$c_1\hat{L}_{\mathrm{lex}}(q)(\uplambda,\pi_{k}(x'))=\hat{L}_{\mathrm{lex}}(q)(\uplambda,x)$
and 
$c_2\hat{L}_{\mathrm{lex}}(p)(\uplambda,\pi_{k}(x'))=\hat{L}_{\mathrm{lex}}(p)(\uplambda,x)$
hold. Hence it is easy to see that there is a positive real number $\lambda_1'$ such that 
$$c_1\hat{q}(\lambda,\pi_{k}(x')),c_2\hat{p}(\lambda,\pi_{k}(x'))>0$$ 
for all $\lambda\geq \lambda_1'$.
On the other side, we still have 
$$c_1\hat{q}(\lambda,x'),c_2\hat{p}(\lambda,x')<0$$ 
for all $\lambda\geq \lambda_2$. We are now ready to construct the rational tentacles needed for the second part of Theorem \ref{kxtt11}. For each $i=1,\ldots,n$ we define $\varphi_i=(\uplambda+\uplambda_3-1)^{z_i}\in \mathbb{R}[\uplambda]$, where $\lambda_3=\max\{\lambda_1',\lambda_2\}$. As usual set $\varphi=(\varphi_1,\ldots,\varphi_n)$. By what we have done so far and by using the fact that, if we substitute $x$ by some nearby other point $y$, none of those inequalities used in this proof will be affected (see Proposition \ref{rooters}), we see that there is an open neighborhood $U$ of $x$ such that $T_{\overline{U},\varphi}\in \mathcal{T}(S(q))$ and $\pi_{k}(T_{\overline{U},\varphi})\notin \mathcal{T}(S(p))$. 
Repeating part (ii)-(iv) of Theorem \ref{kxtt11} concludes the proof.\bew 
\begin{Rem}
Let $n$ be a natural number greater than $2$.
Suppose $f$ and $g$ are two different irreducible homogeneous polynomials 
in $\mathbb{R}[\mathrm{x}_1,\ldots,\mathrm{x}_n]$ of odd degree. If $f$ and $g$ satisfy the condition (a) of Theorem \ref{tult1}, then it has some interesting geometric consequences for $V_1=\mathcal{V}(f)$ and $V_2=\mathcal{V}(g)$. Let $\Lambda_1$ resp. $\Lambda_2$ denote all singular points of $V_1(\mathbb{R})$ resp. $V_2(\mathbb{R})$. According to \cite[Theorem 1, p. 239]{z2} the intersection $V_1(\mathbb{R})\cap V_2(\mathbb{R})$ is not empty. Consider a point $x\in V_1(\mathbb{R})\cap V_2(\mathbb{R})$. The following cases may occur:
\begin{itemize}
\item The point $x$ is in $\Lambda_1$: If $x$ is a local minimum of $g$, then is must also a local minimum of $f$. Thus $x$ is a point in $\Lambda_2$.
\item The point $x$ is not in $\Lambda_1$: It is easy to see that $x$ is a boundary point of $S(g)$. If $x$ is in $\Lambda_2$, then $x$ is either a local minimum of $f$ or it is a saddle point of $f$. If $x$ is not in $\Lambda_2$, then $V_1$ and $V_2$ intersect non-transversely at $x$.     
\end{itemize}  
\end{Rem}    
\begin{window}[3, r, \includegraphics[scale=0.3]{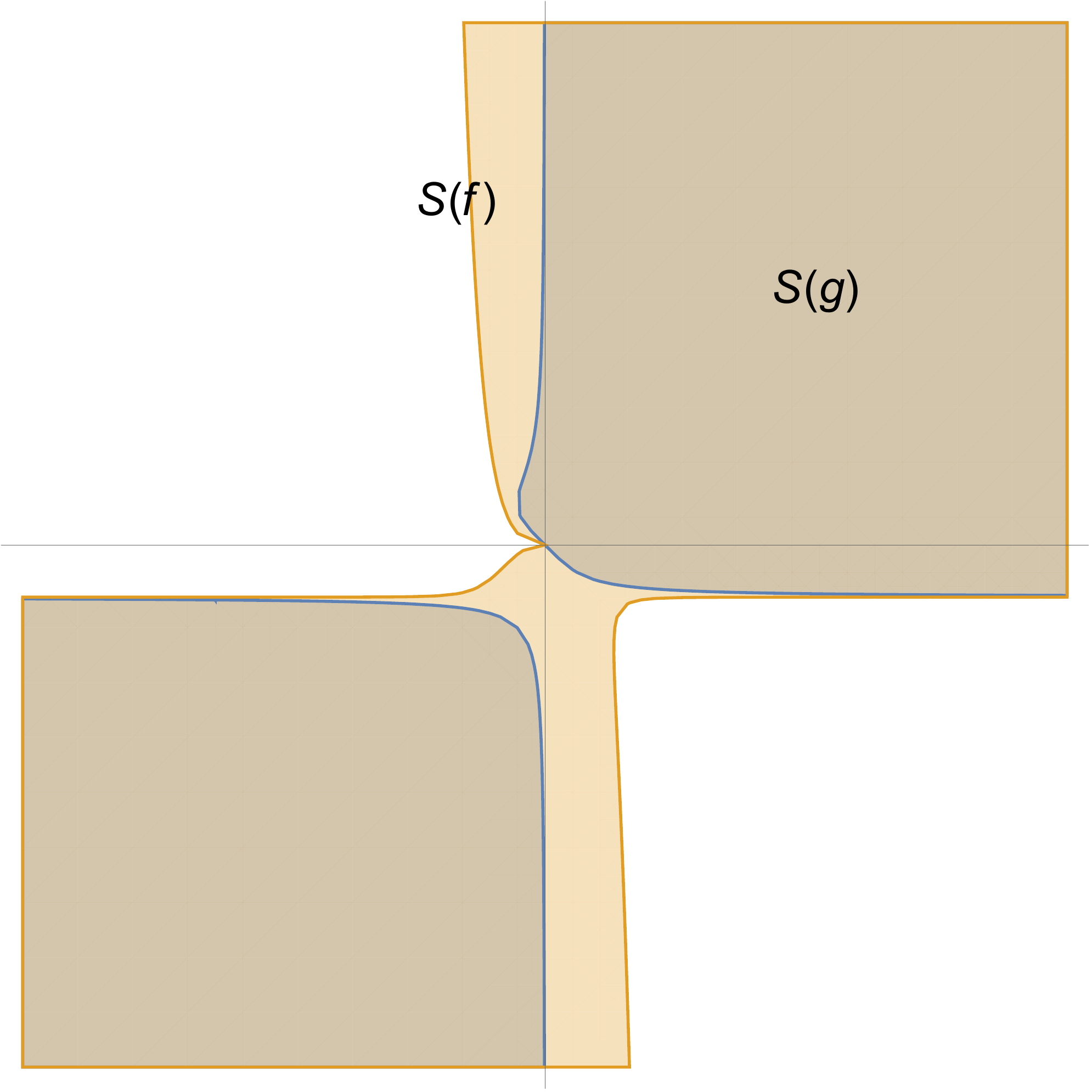},{}]
\begin{Exp}
Consider the polynomials $g=\mathrm{x}_1+\mathrm{x}_2+\mathrm{x}_1\mathrm{x}_2^3$ and $f=\mathrm{x}_1^5+\mathrm{x}_1^5\mathrm{x}_2+\mathrm{x}_2^2$. Let us check if $f$ and $g$ satisfy all conditions of Theorem \ref{tult1}.\\
By using \cite{z15} as we did in \ref{k6.e1x}, we see that $S(g)\subseteq S(f)$.\\
(a): Obviously, $L_{\mathrm{lex}}(g)=\mathrm{x}_1\mathrm{x}_2^3$ and $L_{\mathrm{lex}}(f)=\mathrm{x}_1^5\mathrm{x}_2$. Since $\deg_{\mathrm{x}_2}\left(L_{\mathrm{lex}}(g)\right)>\deg_{\mathrm{x}_2}\left(L_{\mathrm{lex}}(f)\right)$, the polynomial $L_{\mathrm{lex}}(f)$ cannot lie in $L_{\mathrm{lex}}(g)\mathbb{R}[\mathrm{x}_1,\mathrm{x}_2]$.\\
(b): Set $y'=(5,5)$, $y=(5,-5)$, and consider $\hat{g}(\uplambda,y')=g\left(\uplambda^{z_1}y_1',\uplambda^{z_2}y_2'\right)=5\uplambda^{z_1}+5\uplambda^{z_2}+625\uplambda^{z_1+3z_2}$, $\hat{f}(\uplambda,y')=f\left(\uplambda^{z_1}y_1',\uplambda^{z_2}y_2'\right)=3125\uplambda^{5z_1}+25\uplambda^{2z_2}+15625\uplambda^{5z_1+z_2}$, $\hat{g}(\uplambda,y)=g\left(\uplambda^{z_1}y_1,\uplambda^{z_2}y_2\right)=5\uplambda^{z_1}-5\uplambda^{z_2}-625\uplambda^{z_1+3z_2}$, $\hat{f}(\uplambda,y)=f\left(\uplambda^{z_1}y_1,\uplambda^{z_2}y_2\right)=3125\uplambda^{5z_1}+25\uplambda^{2z_2}-15625\uplambda^{5z_1+z_2}$, where $(z_1,z_2)\in \mathbb{Z}^2$. It is easy to see that if we take $(z_1,z_2)\in \mathbb{N}_1$, then $\hat{g}(\lambda,y'), \hat{f}(\lambda,y')$ are positive for all $\lambda\geq 1$, while $\hat{g}(\lambda,y'), \hat{f}(\lambda,y')$ are negative for all $\lambda\geq 1$.
By taking a small compact neighborhood $U'$ of $y'$ resp. $U$ of $y$, we get a tentacle $T_{U',z}$ belonging to $\mathcal{T}(S(g))$ and another one $T_{U,z}$ belonging to $\mathcal{T}(S(-f))$. 
Finally, Theorem \ref{tult1} implies that there is no non-negative polynomial $t\in \mathbb{R}[x_1,x_2]$ such that $f(x)-t(x)g(x)\geq 0$ for all $x\in \mathbb{R}^2$.  
\end{Exp}
\end{window} 
\begin{window}[8, r, \includegraphics[scale=0.3]{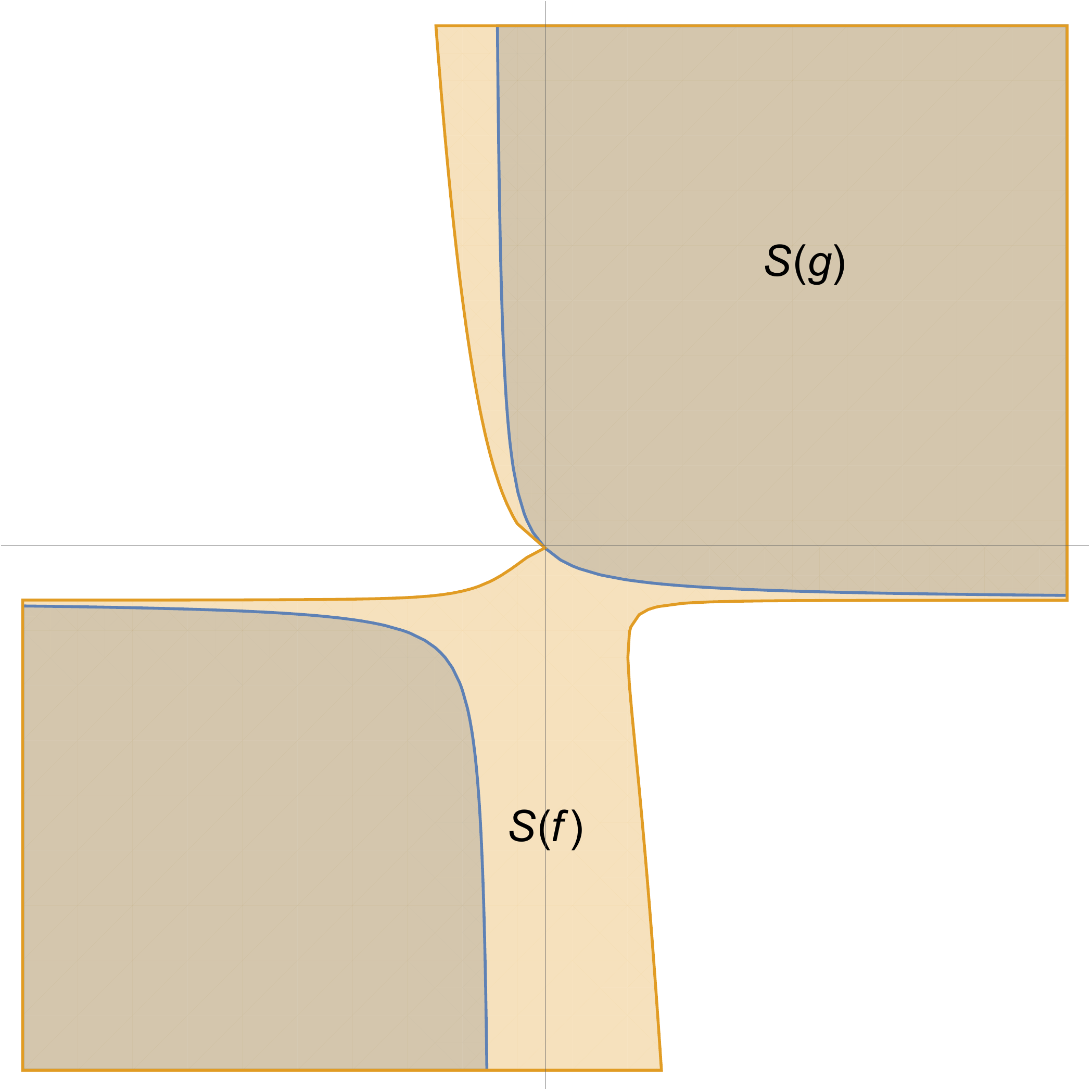},{}]
\begin{Exp}\label{Counter2}
Let us revisit the counterexample $f=\mathrm{x}_1^3+\mathrm{x}_1^3\mathrm{x}_2+\mathrm{x}_2^2$ and $g=\mathrm{x}_1+\mathrm{x}_2+\mathrm{x}_1\mathrm{x}_2$.
Unfortunately, condition (a) of Theorem \ref{tult1} is violated. We have $L_{\mathrm{lex}}(g)=\mathrm{x}_1\mathrm{x}_2$ and $L_{\mathrm{lex}}(f)=\mathrm{x}_1^3\mathrm{x}_2$. Thus $L_{\mathrm{lex}}(f)=\mathrm{x}_1^2L_{\mathrm{lex}}(g)$. On the other side, it is easy to see that the condition (b) of Theorem \ref{tult1} is satisfied. Can we still make use of Theorem \ref{tult1}?\\ 
Let $t\in \mathbb{R}[\mathrm{x}_1,\mathrm{x}_2]$ be a non-negative polynomial. If $L_{\mathrm{lex}}(t)\neq \mathrm{x}_1^2$ then we can just repeat the arguments of Theorem \ref{tult1}. However, the problematic case is if $L_{\mathrm{lex}}(t)=\mathrm{x}_1^2$ holds. Interchange the variables $\mathrm{x}_1$ and $\mathrm{x}_2$ in $f$ and $g$, giving $f=\mathrm{x}_2^3+\mathrm{x}_2^3\mathrm{x}_1+\mathrm{x}_1^2$ and $g=\mathrm{x}_1+\mathrm{x}_2+\mathrm{x}_1\mathrm{x}_2$. Now, we are out for some suitable $z$-grading such that $L_z(f)\notin L_z(g)\mathbb{R}[\mathrm{x}_1,\mathrm{x}_2]$. Choose $z\in \mathbb{N}_1^2$ with $3z_2=z_1$. Then we have $L_z(g)=\mathrm{x}_1\mathrm{x}_2$, $L_z(f)=\mathrm{x}_2^3\mathrm{x}_1+\mathrm{x}_1^2$ resp. $\hat{L}_z(g)(\uplambda,\mathrm{x})=\mathrm{x}_1\mathrm{x}_2\uplambda^{z_1+z_2}$, $\hat{L}_z(f)(\uplambda,\mathrm{x})=\left(\mathrm{x}_2^3\mathrm{x}_1+\mathrm{x}_1^2\right)\uplambda^{2z_1}$. Obviously, $L_z(f)$ does not lie in $L_z(g)\mathbb{R}[\mathrm{x}_1,\mathrm{x}_2]$. 
Let us check condition (b). Set $x=(5,5)$. Then $\hat{L}_z(g)(\uplambda,x)=25\uplambda^{z_1+z_2}$ and $\hat{L}_z(f)(\uplambda,x)=5^4\uplambda^{2z_1}+5^2\uplambda^{2z_1}$. On the other hand, if we set $x=(-5,5)$ we get $\hat{L}_z(g)(\uplambda,x)=-25\uplambda^{z_1+z_2}$ and $\hat{L}_z(f)(\uplambda,x)=-5^4\uplambda^{2z_1}+5^2\uplambda^{2z_1}$. This means that if we take an appropriate neighborhood $U$ of $x$ and define $\pi:\mathbb{R}^2\rightarrow \mathbb{R}^2,(x_1,x_2)\mapsto (-x_1,x_2)$ we get $T_{U,z}\in \mathcal{T}(S(g))$ and $\pi(T_{U,z})\in \mathcal{T}(S(-f))$.
By using the here defined $z$-grading, instead of the one defined in the proof of Theorem \ref{tult1}, and repeating the arguments in part (ii)-(iv) we get that there is no non-negative polynomial $t\in \mathbb{R}[\mathrm{x}_1,\mathrm{x}_2]$ such that $f(y)-t(y)g(y)\geq 0$ for all $y\in \mathbb{R}^2$. In fact we can state: 
\end{Exp}
\end{window}
\begin{Theo}\label{stxf2}
Let $S_1=S(q)$ and $S_2=S(p)$ be two semi-algebraic sets in $\mathbb{R}^n$. Suppose the following condition is satisfied:\par\smallskip 
There is a $z\in \mathbb{N}_1^n$ such that $L_z(p)\notin L_z(q)A$, a rational tentacle $T\in \mathcal{T}(S_1)$ 
of degree $z$ and an element $\pi\in G$ such that $\pi(T)\notin \mathcal{T}(S_2)$. Furthermore, all unbounded $T'\in \mathcal{T}(S_1)$ with $\pi(T')\subseteq \pi(T)$ satisfy $\pi(T')\notin \mathcal{T}(S_2)$.\par\smallskip 
Then there is no non-negative polynomial $t\in A$ such that $p(y)-t(y)q(y)\geq 0$ for all $y\in \mathbb{R}^n$.
\end{Theo}
\textbf{Proof}: Simply repeat the same arguments in part (ii)-(iv) of Theorem \ref{tult1}.\bew\\\\
In the same manner we can modify Proposition \ref{scor}:
\begin{Prop}\label{scor2}
Let $S_1=S(q)$ and $S_2=S(p)$ be two semi-algebraic sets in $\mathbb{R}^n$, where neither $q=0$ nor $p=0$. Suppose the following condition is satisfied:\par\smallskip 
There is a $z\in \mathbb{N}_1^n$ such that $L_z(p)\notin L_z(q)A$
and the quadratic modules $\mathrm{QM}(q)$ and $\mathrm{QM}(-p)$ are totally stable with respect to the $z$-grading.\par\smallskip 
Then there is no non-negative polynomial $t\in A$ such that $p(y)-t(y)q(y)\geq 0$ for all $y\in \mathbb{R}^n$.
\end{Prop}
\begin{Rem}
The difference between Theorem \ref{tult1} and Theorem \ref{stxf2} is simple. In Theorem \ref{tult1} we demanded that $L_{\mathrm{lex}}(p)\notin L_{\mathrm{lex}}(q)A$. Then we constructed a special $z$-grading, where $z\in \mathbb{N}_1^n$. Under this $z$-grading we had $L_z(p)=L_{\mathrm{lex}}(p)$ resp. $L_z(q)=L_{\mathrm{lex}}(q)$ and therefore $L_z(p)\notin L_z(q)A$. 
Let us refer to the $z$-gradings satisfy $z\in \mathbb{N}_1^n$ and $L_z(p)\notin L_z(q)A$ as special gradings.
The difference between Theorem \ref{stxf2} and Theorem \ref{tult1} is, that the special $z$-grading constructed in part (i) of the proof of Theorem \ref{tult1}, is already given in the prerequisites of Theorem \ref{stxf2}. The disadvantage of Theorem \ref{stxf2} compared to Theorem \ref{tult1} is, that one must find such a special $z$-grading for Theorem \ref{stxf2} to work, while Theorem \ref{tult1} does not require such a procedure. The advantage of Theorem \ref{stxf2} is, that it allows a wider range of gradings as Example \ref{Counter2} illustrates it.
In fact, we have an two explanations why we used $g=\mathrm{x}_1\mathrm{x}_3+\mathrm{x}_2\mathrm{x}_3+\mathrm{x}_1\mathrm{x}_2$. 
The first explanation is a geometric one. According to Example \ref{Counter2} the polynomial $g$ has all the necessary properties for Theorem \ref{stxf2} to work.
The second explanation is an algebraic one. 
According to Theorem \ref{kxt1} quadratic forms $q$ that have a vanishing diagonal part give rise to quadratic module $\mathrm{QM}(q)$ that is totally stably with respect to any $z$-grading. Thus it is (was) convenient to choose $g=\mathrm{x}_1\mathrm{x}_3+\mathrm{x}_2\mathrm{x}_3+\mathrm{x}_1\mathrm{x}_2$.   
\end{Rem}
\section{A non-geometric counterexample}
\begin{window}[2, r, \includegraphics[scale=0.20]{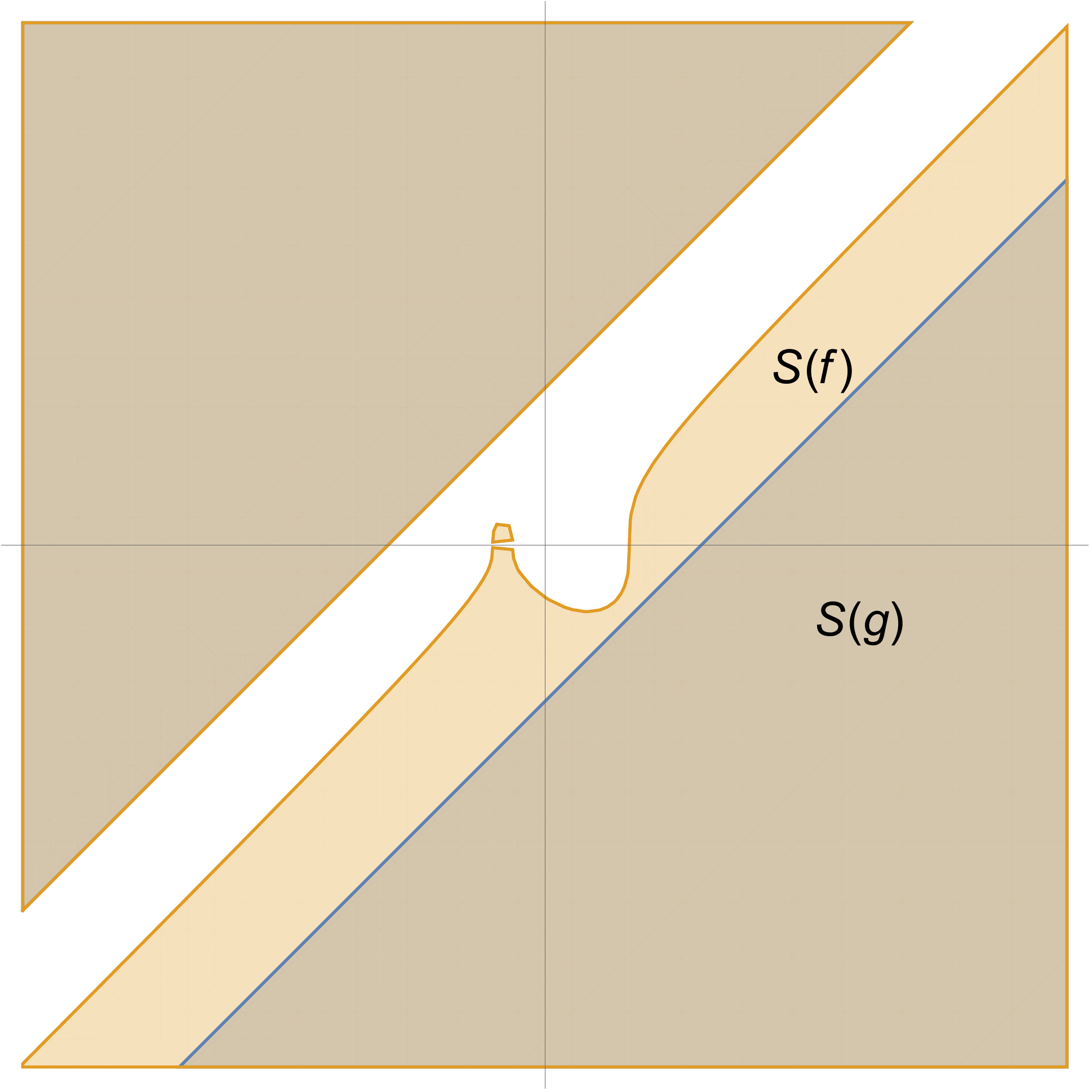},{}]
In the last chapter we saw, that Counterexample \ref{k6.e1} proved the S4-conjecture wrong because of geometric reasons. The straightforward question is obvious: Can we find for all counterexamples a geometric reason? At least,we will give a counterexample that does not work because of a arithmetical reason. From now on, we set $g=(-3+\mathrm{x}_1-\mathrm{x}_2)(3+\mathrm{x}_1-\mathrm{x}_2)$, $p=-\mathrm{x}_1^3+\mathrm{x}_2^3+2\mathrm{x}_1+1$, $l_1=-3+\mathrm{x}_1-\mathrm{x}_2$, $l_2=3+\mathrm{x}_1-\mathrm{x}_2$ and $f=-l_2p$. 
By using \cite{z15} (see \ref{k6.e1}) it is easy to verify that $S(g)\subseteq S(f)$. The next obvious step is:
\end{window}
\begin{Prop}\label{cc1x}
There is no non-negative polynomial $t\in \mathbb{R}[\mathrm{x}_1,\mathrm{x}_2]_2$ such that $f(y)-t(y)g(y)\geq 0$ for all $y\in \mathbb{R}^2$.
\end{Prop} 
\textbf{Proof}: Suppose there is a non-negative polynomial $t\in \mathbb{R}[\mathrm{x}_1,\mathrm{x}_2]_2$ contradicting the statement of this proposition. In the following we fix a real number $x_2\in \mathbb{R}$ and consider $f(\mathrm{x}_1,x_2),g(\mathrm{x}_1,x_2),t(\mathrm{x}_1,x_2)$ as polynomials in $\mathrm{x}_1$. 
The polynomial $f(\mathrm{x}_1,x_2)-t(\mathrm{x}_1,x_2)g(\mathrm{x}_1,x_2)\in \mathbb{R}[\mathrm{x}_1]$ has a root at $x_1=x_2-3$. Since $f(\mathrm{x}_1,x_2)-t(\mathrm{x}_1,x_2)g(\mathrm{x}_1,x_2)\in \mathbb{R}[\mathrm{x}_1]$ is non-negative for every $x_2\in \mathbb{R}$, it must be divided by $l_2^2(\mathrm{x}_1,x_2)\in \mathbb{R}[\mathrm{x}_1]$.\par\smallskip  
The remainder of the polynomial division $f(\mathrm{x}_1,x_2):l_2^2(\mathrm{x}_1,x_2)$, as polynomials in $\mathbb{R}[\mathrm{x}_1]$, is $r_1(\mathrm{x}_1,x_2)=9x_2^3-9\mathrm{x}_1x_2^2-52x_2^2-25\mathrm{x}_1x_2+97x_2^2-22\mathrm{x}_1-66\in \mathbb{R}[\mathrm{x}_1]$. And for $g(\mathrm{x}_1,x_2):l_2^2(\mathrm{x}_1,x_2)$ we have $r_2(\mathrm{x}_1,x_2)=6x_2-6\mathrm{x}_1-18\in \mathbb{R}[\mathrm{x}_1]$. Finally, let $r_3(\mathrm{x}_1,x_2)\in \mathbb{R}[\mathrm{x}_1]$ denote the remainder of $t(\mathrm{x}_1,x_2):l_2^2(\mathrm{x}_1,x_2)$. Since $f(\mathrm{x}_1,x_2)-t(\mathrm{x}_1,x_2)g(\mathrm{x}_1,x_2)$ is divided by $l_2^2(\mathrm{x}_1,x_2)$, we get the identity $r_1(\mathrm{x}_1,x_2)-r_3(\mathrm{x}_1,x_2)r_2(\mathrm{x}_1,x_2)=0$. This leads to $r_3(\mathrm{x}_1,x_2)=\frac {r_1(\mathrm{x}_1,x_2)} {r_2(\mathrm{x}_1,x_2)}=\frac {1} {6}\left(22-25x_2+9x_2^2\right)$.
Set $\tilde{t}=a_{x_2}l_2^2+\frac {1} {6}\left(22-25\mathrm{x}_2+9\mathrm{x}_2^2\right)$ and choose $a_{x_2}\in \mathbb{R}$ such that the equality $\tilde{t}(\mathrm{x}_1,x_2)=t(\mathrm{x}_1,x_2)$ holds for $x_2\in \mathbb{R}$.
It is easy to see that the leading term of $f-\tilde{t}g$ in $\mathrm{x}_2$ is $\left(-\frac {1} {2}-a_{x_2}\right)\mathrm{x}_2^4$ and that $f(0,0)-\tilde{t}(0,0)g(0,0)=30+81a_{x_2}$. For large $x_2\in \mathbb{R}$ we see that $a_{x_2}$ must satisfy $a_{x_2}\leq -\frac {1} {2}$ and $a_{x_2}\geq -\frac {30} {81}$, which is impossible.\bew       
\begin{Rem}
Under an appropriate change of coordinates the homogenization of $g$ can be written as $\overline{g}=-9\mathrm{x}_1^2+2\mathrm{x}_2^2$. Thus the signature of $\overline{g}$ is $0$. Note, that the signature in Counterexample \ref{k6.e1x} was $-1$. Furthermore, it is easy to see that neither Theorem \ref{tult1} nor Theorem \ref{stxf2} can be applied. In fact, $g$ violates condition (b) in Theorem \ref{tult1} and the condition in Theorem \ref{stxf2}.    
\end{Rem}
\section{Final thoughts}
Finally, we have come to the end of this article. Hence let us summarize what we have learned so far. First, we learned that the S4-conjecture is not true. Second, we learned that there are geometric reasons why the S4-conjecture cannot be true. Finally, we learned that there are arithmetic reasons why the S4-conjecture cannot work. 
Still there are many questions left. The signature of $g$ in \ref{k6.e1x} resp. \ref{cc1x} was $-1$ resp. $0$. So it is quite naturally to ask, if there is an counterexample, where $g$ has signature $1$.   
There is another obvious question: Under which conditions does the S4-conjecture work? Can these conditions be expressed in geometric or algebraic terms? Results in this direction can be found in \cite[Corollary 4.5]{z17}:
\begin{Prop}\label{sccco}
\cite[Corollary 4.5]{z17}: Let $h_1,\ldots,h_r\in \mathbb{R}[\mathrm{x}_0,\ldots,\mathrm{x}_n]$ be homogeneous polynomials of even degree, and let 
$$S=S(h_1,\ldots,h_r).$$
Assume there is $\xi\in \mathbb{R}^{n+1}$ with $h_i(\xi)>0$ for $i=1,\ldots,r$. If $p,q\in \mathbb{R}[\mathrm{x}_0,\ldots,\mathrm{x}_n]$ are homogeneous of even degree and positive on $S\backslash \{0\}$, then $pq^m$ lies in the preordering generated by $h_1,\ldots,h_r$, for all sufficiently large $m\geq 0$. 
\end{Prop}
\noindent In contrast to the S4-conjecture, we need a homogeneous polynomial $p$ of even degree that is positive on the set $S\backslash \{0\}$. And even then we can only conclude that there is a natural number $m\geq 1$ such that $p^m\in T(h_1,\ldots,h_r)$. If we could show that $p\in T(h_1,\ldots,h_r)$, we would still have considerable obstacles. For example, we do not know what kind of degree bounds the various representations of $p$ in $T(h_1,\ldots,h_r)$ have. Nevertheless, let us consider the polynomials $g=\mathrm{x}_1^2+\mathrm{x}_2^2-4\mathrm{x}_3^2$ and $f=\mathrm{x}_1^4+\mathrm{x}_2^4-\mathrm{x}_3^4$.
Then $f$ and $g$ satisfy the following conditions:
\begin{itemize}
\item There is a point $x'\in \mathbb{R}^3$ such that $g(x')>0$.
\item The polynomial $f$ is positive on the set $S(g)\backslash \{0\}$.
\item The polynomial $f-\mathrm{x}_3^2g$ is non-negative.
\end{itemize}
In other words, the two polynomials satisfy the S4-conjecture. 
It is easy to see that $V_1$ and $V_2$ are both non-singular curves. According to \cite[Proposition 11.6.2, p. 286]{z9} the set $V_1(\mathbb{R})\subseteq \mathbb{P}^2(\mathbb{R})$ 
decomposes into 1 or 0 ovals, and the set $V_2(\mathbb{R})\subseteq \mathbb{P}^2(\mathbb{R})$ into at most $4$ different ovals. For the definition of an oval see \cite[p. 286]{z9}. In our case $V_1(\mathbb{R})$ and $V_2(\mathbb{R})$ decompose into one oval that does not intersect the plane at infinity. Since $\mathrm{x}_1^2+\mathrm{x}_2^2$ and $\mathrm{x}_1^4+\mathrm{x}_2^4$ are positive on the set $\mathbb{R}^2\backslash \{0\}$, the two sets $V_1(\mathbb{R})$, $V_2(\mathbb{R})$ do not intersect the hyperplane at infinity. Another consequence is that we can tell something about the geometry of $S(g)$ resp. $S(f)$. For $y\in \mathbb{R}\backslash \{0\}$ set $H_y=\mathbb{R}^2\times \{y\}$ and interpret $V_1$ and $V_2$ as affine varieties in $\mathbb{A}^3$. Then we have    
$\partial S(g)\cap H_y=V_1(\mathbb{R})\cap H_y$ and $\partial S(f)\cap H_y=V_2(\mathbb{R})\cap H_y$ for all $y\in \mathbb{R}\backslash \{0\}$. The two inclusions $\partial S(g)\cap H_y\subseteq V_1(\mathbb{R})\cap H_y$ and $\partial S(f)\cap H_y\subseteq V_2(\mathbb{R})\cap H_y$ are obvious. The other two inclusions $\partial S(g)\cap H_y\supseteq V_1(\mathbb{R})\cap H_y$ and $\partial S(f)\cap H_y\supseteq V_2(\mathbb{R})\cap H_y$ hold, because $V_1(\mathbb{R})$ and $V_2(\mathbb{R})$ do not have singular points in $H_y$. 
Thus each slice $\partial S(f)\cap H_y$, $\partial S(g)\cap H_y$ 'looks' like a circle. 
If we replace $g$ by an arbitrary quadratic form $q\in \mathbb{R}[\mathrm{x}_1,\mathrm{x}_2,\mathrm{x}_3]$ and $f$ by an arbitrary ternary quartic $p$, then all geometric statements\footnote{Of course with an adjusted number of ovals} made so far in this matter remain true, if the ovals of $\mathcal{V}(q)(\mathbb{R})$ and $\mathcal{V}(p)(\mathbb{R})$ do not intersect the plane at infinity.
Interestingly, the ovals that do not intersect the plane at infinity have different topological properties than their counterparts that intersect the plane at infinity: It is a well known fact that the fundamental group of $\mathbb{P}^2(\mathbb{R})$ is exactly $\mathbb{Z}/2\mathbb{Z}$. By interpreting an oval as a loop, it turns out that all ovals that do not intersect the hyperplane at infinity, represent the identity element of the fundamental group. 
If all ovals of $\mathcal{V}(q)(\mathbb{R})$ and $\mathcal{V}(p)(\mathbb{R})$ do not intersect the hyper plane at infinity, then it is obvious that the topological situation compared to $\mathcal{V}(g)(\mathbb{R})$ and $\mathcal{V}(f)(\mathbb{R})$ has not changed much.\\  
Hence it is convenient to ask this final question:
\begin{Quest}
Let $q\in \mathbb{R}[\mathrm{x}_1,\mathrm{x}_2,\mathrm{x}_3]$ be a quadratic form and $p$ a ternary quartic. Set $V_1=\mathcal{V}(q)\subseteq \mathbb{P}^2$ and $V_2=\mathcal{V}(p)\subseteq \mathbb{P}^2$. Suppose that the following conditions are satisfied:
\begin{itemize}
\item There is a point $x'\in \mathbb{R}^3$ such that $q(x')>0$.
\item The ternary quartic $p$ is positive on the set $S(q)\backslash \{0\}$.
\item The projective varieties $V_1$ and $V_2$ are non-singular.
\item The set $V_1(\mathbb{R})$ is an oval and $V_2(\mathbb{R})$ decomposes into at least one oval. 
\item All ovals of $V_1(\mathbb{R})$ and $V_2(\mathbb{R})$ do not intersect the hyperplane at infinity.  
\end{itemize}
Can we find a non-negative homogeneous polynomial $t\in \mathbb{R}[\mathrm{x}_1,\mathrm{x}_2,\mathrm{x}_3]_2$ such that $p(y)-t(y)q(y)\geq 0$ for all $y\in \mathbb{R}^3$?
\end{Quest}

\end{document}